\documentclass[reqno,12pt,oneside]{amsart}
\usepackage{mathtext}
\usepackage{amsmath}
\usepackage{amsfonts}
\usepackage{amssymb}
\usepackage{amsthm}
\usepackage{verbatim}
\usepackage{enumerate}
\usepackage{version,etex}
\usepackage[top=3cm,bottom=2.5cm,left=2.5cm, right=3cm,marginparwidth=0cm]{geometry}
\usepackage[pdftex]{graphicx,color}
\usepackage{epsfig}
\usepackage[colorlinks]{hyperref}
\usepackage{amsmath}
\usepackage[displaymath, tightpage]{preview}

\usepackage{caption}

\thanks{Both authors were supported by ERC AdG grant no 339523 RGDD,
and Sof\'ia Trejo was additionally supported by FAPESP Grant No. 2014/09418-0}


\title{The boundary of chaos for interval mappings}
\date{}
\author{Trevor Clark} 
\address{Trevor Clark, Imperial College London, London, UK}
\email[]{t.clark@imperial.ac.uk}

\author{Sof\'ia Trejo}
\address{Sof\'ia Trejo, UNAM, Mexico City, Mexico}
\email[]{sofia.trejo.a@gmail.com}

\begin{document}
\setcounter{secnumdepth}{4}
\setcounter{tocdepth}{1}
\newcommand{\proclaim}[2]{\medbreak {\bf #1}{\sl #2} \medbreak}

\newcommand{\ntop}[2]{\genfrac{}{}{0pt}{1}{#1}{#2}}

\let\newpf\proof \let\proof\relax \let\endproof\relax
\newenvironment{pf}{\newpf[\proofname]}{\qed\endtrivlist}

\def\basins{\mathcal{B}}
\def\immbasins{\mathcal{B}_{0}}
\def\LL{\mathcal{L}}
\def\comp{\mathrm{Comp}}
\def\N{\mathbb{N}}
\def\I{\boldsymbol{I}}
\def\L{\boldsymbol{L}}
\def\K{\boldsymbol{K}}
\def\J{\boldsymbol{J}}
\def\U{\boldsymbol{U}}
\def\poin{\boldsymbol{Poin}}
\def\Y{\mathcal{Y}}

\def\cascend{I^{\hat{m}}}
\def\PC{\mathrm{PC}}
\def\crit{\mathrm{Crit}}
\def\R{\mathbb{R}}
\def\C{\mathbb{C}}
\def\dom{\mathrm{Dom}}

\newcommand{\parabolic}{\operatorname{Parabolic}}
\newcommand{\remark}{\noindent\textbf{Remark.}}

\newtheorem{thm}{Theorem}[section]
\newtheorem{cor}[thm]{Corollary}
\newtheorem{lem}[thm]{Lemma}
\newtheorem{prop}[thm]{Proposition}

\numberwithin{equation}{section}
\newcommand{\thmref}[1]{Theorem~\ref{#1}}
\newcommand{\propref}[1]{Proposition~\ref{#1}}
\newcommand{\secref}[1]{\S\ref{#1}}
\newcommand{\lemref}[1]{Lemma~\ref{#1}}
\newcommand{\corref}[1]{Corollary~\ref{#1}} 

\theoremstyle{definition}
\newtheorem{defn}{Definition}[section]
\newtheorem{definition}{Definition}[section]

\newcommand{\diam}{\operatorname{diam}}
\newcommand{\dist}{\operatorname{dist}}
\newcommand{\cl}{\operatorname{cl}}
\newcommand{\inter}{\operatorname{int}}
\renewcommand{\mod}{\operatorname{mod}}

\newcommand{\s}{\vspace{5mm}}

\maketitle

\begin{abstract}
A goal in the study of dynamics on the 
interval is to understand the
transition to positive topological entropy.
There is a conjecture from
the 1980's that the only
route to positive topological entropy is through a
cascade of period doubling bifurcations.
We prove this conjecture
in natural families of smooth interval maps,
and use it to study the structure of the boundary of mappings with
positive entropy. 
In particular, we show that in
families of mappings with a fixed number of critical points
the boundary is locally connected,
and for analytic mappings that it is a cellular set.
\end{abstract}

\newtheorem{thmx}{Theorem}
\renewcommand{\thethmx}{\Alph{thmx}}

\section{Introduction}
This paper is motivated by the following 
conjectures in one-dimensional dynamics
about the boundary of mappings with positive
topological entropy:

\medskip
Given a map $f$ of an interval, $I$, let 
$$\mathcal{P}er(f)=\{n\in \N : f^n(p)=p \text{  for some  } p\in I,\mbox{ and
}f^j(p)\neq p \mbox{ for } 1\leq j<n \}.$$
We refer to $\mathcal Per(f)$ as the {\em set of periods of} $f$.

\medskip
\noindent\textbf{Boundary of Chaos Conjecture I.}
{\em All endomorphisms of the interval,
$f\in \mathcal{C}^k(I), k=0,1,2,\dots,\infty,\omega,$
with $\mathcal{P}er(f)=\{2^n:n\in\mathbb N\cup\{0\}\},$
are on the boundary
of mappings with positive topological
entropy and on the  boundary of the set of mappings with
finitely many periods.}

\medskip

Interest in this conjecture is strongly 
motivated by its implications on the 
{\em routes to chaos}, that is, on the transition from zero to
positive entropy, for mappings of the
circle or the interval, see \cite{MaT-circle} and \cite{MaT-bimodal}.
Indeed, for $\mathcal C^1$ mappings,
this conjecture implies that the transition to positive entropy for
mappings on the interval occurs through successive period doubling
bifurcations.

In \cite{OT}, the following conjecture was made about the internal
structure
of the boundary of mappings with positive topological entropy:

\medskip
\noindent\textbf{Boundary of Chaos Conjecture II.}
{\em
An open and dense subset of the boundary of 
mappings with positive topological entropy
splits into disjoint cells such that each cell
is contained in the basin of the quadratic-like
fixed point of renormalization. See \cite{OT}
for a more precise statement.
}
\medskip
 
Conjecture I was first made for the space
$\mathcal C^1(I)$ in \cite{BH} and later for each 
$k=0,1,2,\dots,\infty,\omega,$ see
\cite{MaT-circle}, \cite{MaT-bimodal}, \cite{TH} and \cite{OT}.
It is known in $\mathcal C^0(I)$ and $\mathcal C^1(I)$.
In \cite{Kl}, it was proved that mappings with
positive topological entropy are dense in $f\in\mathcal{C}^0(I).$
In fact,
for any compact manifold $M,$
infinite topological entropy is a generic
property amongst endomorphisms of $M$ in 
the $\mathcal C^0$ topology \cite{Y}. In \cite{JS}, it was proved that
any $f\in\mathcal{C}^0(I)$
with $\mathcal{P}er(f)=\{2^n:n\in\mathbb N\cup\{0\}\}$
can be approximated by mappings with
finitely many periods. 
In \cite{J},
Conjecture I was proved for $\mathcal C^1(I)$.
The results in lower regularity are perturbative,
and this approach does not seem to work
in higher regularity.

\subsection{Main results}
To fix some notation, let $I=[-1,1]$ and
$\underline b=(\ell_1,\ldots,\ell_b)$ be
a vector of even integers greater than one.
Let $\mathcal A_{\underline b}(I)$
denote the space of 
analytic mappings of the interval,
with critical points $-1<c_1<c_2<\dots<c_b<1,$
such that the order of $c_i$ is $\ell_i.$ We describe this condition precisely
on page~\pageref{page:Ak}.
If $U\subset\mathbb C$ is open, we let
$\mathcal B_U$ denote the space of 
mappings that are holomorphic on $U$ and
continuous on $\overline U$. We consider
$\mathcal B_U$ with the supremum norm.
We prove the following result for
analytic mappings:
\begin{thmx}\label{thm:main}
\em{
All analytic endomorphisms 
$f\in\mathcal A_{\underline b}(I)$ with
$\mathcal{P}er(f)=\{2^n:n\in\mathbb N\cup\{0\}\}$
are on the boundary of mappings with positive topological
entropy and on the  boundary of
mappings with finitely many periods
in $\mathcal A_{\underline b}(I)$. 

More precisely, suppose that
$f\in\mathcal A_{\underline b}(I),$ with
$\mathcal{P}er(f)=\{2^n:n\in\mathbb N\cup\{0\}\}.$
Let $U\subset\mathbb C,$ $U\supset I$, be an open domain so that
$f\in\mathcal B_U$ and each critical point of 
$f|_U$ is in $I$. Then $f$ can be approximated in
$\mathcal A_{\underline b}(I)\cap\mathcal B_U$ by mappings with positive entropy
and by mappings with finitely many periods.
}
\end{thmx}

Recall that by Sharkovskii's Theorem, an interval mapping
$f$ has finitely many periods if and only if 
for some $N\in\mathbb N\cup\{0\},$
$\mathcal Per(f)=\{2^n:0\leq n\leq N\}.$
Let us point out that for $k$ at least one,
mappings with finitely many periods 
are in the interior of
mappings with zero topological entropy in $\mathcal C^k(I)$,
\cite[Proposition 2.1]{Mi1}. There, this result was proved with $k=1$, but the proof  goes through for any $k\geq 1$. Thus one may 
replace ``boundary of mappings with finitely
many periods'' with ``boundary of the interior of the set of mappings with
zero entropy'' in the statement of the Theorem~\ref{thm:main}.
Let us also recall that a mapping with with 
$\mathcal Per(f)=\{2^n:n\in\mathbb N\cup\{0\}\}$ 
has zero entropy \cite[Theorem 4]{M-horse}.

Theorem~\ref{thm:main} is closely related to the
\emph{Density of hyperbolicity}, \cite{KSS-density}, which
tells us that every mapping can be approximated by 
mappings where every critical point converges to 
a periodic attractor, but it does not specify the combinatorics of 
the mapping used to carry out the approximation. 
Conjecture I implies that for mappings $f$ with
$\mathcal{P}er(f)=\{2^n:n\in\mathbb N\cup\{0\}\},$
this approximation can be done in two combinatorially
different ways, and it specifies the combinatorics of the
approximating mappings with zero entropy precisely.

Our method used to prove Theorem A leads us to
the following result, which was inspired by Conjecture II:

\begin{thmx}\label{thm:II}
{\em
The boundary of mappings
with positive entropy, $\Gamma\subset \mathcal A_{\underline b}(I),$
admits a cellular decomposition.
Moreover, there exists an open and dense subset of 
$\Gamma$ consisting of disjoint 
cells, each contained in the basin of a
unimodal, polynomial-like fixed point of 
renormalization. 
}
\end{thmx}

A {\em cell} is a connected set of (finite) codimension-$k$
whose boundary contains a (relatively) open
and dense set of codimension-$(k+1)$.
A set $X$ {\em admits a cellular decomposition} 
if it can be expressed as a disjoint union of 
cells.
By {\em basin} we mean the set of all mappings $f$
which have a critical point $c$ (of order $\ell$)
at which $f$ is infinitely renormalizable with period doubling combinatorics,
and whose
renormalizations at $c$ converge to the unimodal
fixed point of renormalization whose critical point is of order $\ell$.

Theorem~\ref{thm:II} implies that there is an open and dense set $\Gamma'$ of mappings in
$\Gamma,$ so that each $f\in\Gamma'$ has a critical point $c_0$ at which $f$ is
infinitely renormalizable, and with the property that the symbolic dynamics on
the solenoidal attractor  $\omega(c_0)$  is the same as the 
symbolic dynamics on the solenoidal attractor for the
unimodal Feigenbaum mapping.

Using the complex bounds of \cite{ClvSTr},
we are able to extend Theorem~\ref{thm:main}
to spaces of smooth mappings with critical points
of even order.

\begin{thmx}\label{thm:main smooth}
{\em Let $k\geq 3$ and $b\in\mathbb N$.
If $\underline b$ is a $b$-tuple with only even entries,
then each $f\in\mathcal A_{\underline b}^k(I)$
with
$\mathcal{P}er(f)=\{2^n:n\in\mathbb N\cup\{0\}\}$
is on the boundary of mappings with positive topological
entropy and on the  boundary of the set of
mappings with finitely many periods in 
$\mathcal A_{\underline b}^{k}(I).$

% In greater generality,
% for each $b$-tuple with integer entries,
% $\underline b,$ there exists $M, b'\in\mathbb N$ so that
% the following holds:
% Each $f\in\mathcal A_{\underline b}^k(I)$,
% with
% $\mathcal{P}er(f)=\{2^n:n\in\mathbb N\}$
% is on the boundary of mappings with positive topological
% entropy and on the  boundary of the set of
% mappings with finitely many periods in 
% $\mathcal A_{M,b'}^{k}(I).$
}
\end{thmx}

See page~\pageref{page:Ak} for the definition of the
space $\mathcal A^k_{\underline b}(I)$. 

\medskip

To prove Conjecture II
for certain smooth mappings, we make use of the
hyperbolicity of renormalization of
$\mathcal C^{2+\alpha}$-unimodal mappings with period-doubling
combinatorics, \cite{Davie}.
See \cite{dFdMP} for the generalization of 
this result to all bounded combinatorics.
The hyperbolic structure at the 
quadratic-like fixed point of renormalization,
gives us a means to understand the 
structure of the set of mappings 
on the boundary of positive entropy in spaces of mappings
with several critical points. 
We let $\mathcal{A}^r_{even, b}(I)$ denote the space of
mappings with $b$ critical points all of even order,
see page~\pageref{page:Aeven}.

\begin{thmx}\label{thm:II smooth}
{\em 
There exists an open and dense set of mappings contained in the boundary of
positive entropy in $\mathcal A_{even, b}^r(I),$
$r> 3,$ which is 
a union of disjoint codimension-one submanifolds of  $\mathcal A_{even, b}^r(I),$
and each of these submanifolds is contained in  the basin of a unimodal, quadratic-like 
fixed point of renormalization.
}
\end{thmx}
Specifically, the dense set of mappings which can be decomposed
into codimension-one manifolds consists of mappings with all critical points
non-degenerate and with exactly one solenoidal attractor.
The boundaries of these manifolds 
contain mappings where the solenoidal
attractor contains more than one critical point.
Since we do not know that sets of such mappings
are manifolds, we are unable to obtain the 
cellular decomposition of the boundary of positive entropy for smooth mappings.

\medskip

The following is an
 interesting consequence of Theorem~\ref{thm:II smooth}.
\begin{thmx}\label{thm:III}
{\em Let $r>3$, and let $\underline b$ be a $b$-tuple
of even integers. The connected components of 
the boundary of mappings with 
positive topological entropy in $\mathcal{A}^r_{\underline b}(I)$
are locally connected.
}
\end{thmx}
This result should be contrasted with the theorem of \cite{FrT} that
the boundary of mappings with positive entropy in the family of
bimodal  mappings of the circle is not locally connected, and the
result of \cite{BvS2}, which shows that many isentropes in families of
polynomials are not locally connected. Let us point out that the
mechanisms used to produce non-local connectivity in these cases are
not present in our setting. The result of \cite{FrT} relies on there
being an accumulation of pieces of Arnold tongues in the boundaries of
phase locking regions with definite ``height" above the critical line
in the boundary of mappings with positive entropy. This phenomenon
creates a comb-like structure in the boundary. 
The families considered in
\cite{BvS2} do not have a constant number of critical points,
and its proof that certain isentropes are
non-locally connected requires that the entropy of the
isentrope is positive.

Mappings with 
$\mathcal Per(f)=\{2^n:n\in\mathbb N\cup\{0\}\}$ are
infinitely renormalizable \cite[Theorem 3]{TH},
see Section~\ref{sec:entropy and renormalization},
and such mappings have been the subject of intense
study over the past thirty years. 
Previous results in the direction of those in this paper
have been obtained via proofs of the Hyperbolicity of 
Renormalization Conjectures, \cite{TC,CT,F1,F2}
(or at least convergence of renormalization
together with certain rigidity results, \cite{Sm2}).
For unimodal mappings with critical points of even order,
the solution of the renormalization conjectures imply, roughly,
that the connected component of the
stable manifold containing the fixed point $f_*$
of the (period-doubling) renormalization
operator consists of mappings which are topologically conjugate
to $f_*,$ and the family 
$\{f_*+\lambda v\},$ where $v$ is the expanding direction for
renormalization and $\lambda\in(-\varepsilon,\varepsilon)$,
is transverse to the topological conjugacy class of
$f_*$.
Moreover, $f_*$ is a
polynomial-like mapping, which is hybrid conjugate to the Feigenbaum 
polynomial, and the family $z\mapsto z^2+c$ is transverse to the 
topological conjugacy class of $f_*$ too. Thus one obtains Conjecture I 
for such unimodal mappings from the solution of the renormalization
conjectures together with the solution of Conjecture I for 
unicritical, real, polynomials and Douady-Hubbard Straightening Theorem.
Theorem~\ref{thm:main} has been proved
for analytic unimodal mappings
\cite{Lanford, Sullivan, L-Feig}.
Renormalization results for smooth unimodal mappings with 
quadratic critical points
were obtained in \cite{Davie} and \cite{dFdMP}.
In \cite{dFdMP}
for $\gamma\in(0,1),$ sufficiently close to one, the authors proved
hyperbolicity of renormalization (with bounded combinatorics)
for $\mathcal C^{2+\gamma}$ mappings, and proved that the stable manifold
of the renormalization operator is a $\mathcal C^1$ codimension-one
submanifold of the space of $\mathcal C^{3+\gamma}$ mappings.
Thus proving Theorem~\ref{thm:main} for $\mathcal C^{3+\gamma}$
unimodal mappings with
non-degenerate critical points. 
In \cite{Sm2}, using convergence of renormalization and rigidity,
Smania proved Conjecture I for 
multimodal mappings with all critical points non-degenerate 
and with the same $\omega$-limit set (indeed, 
in \cite{Smania-shy} he goes beyond this to 
prove hyperbolicity of renormalization for these mappings).
In this paper, we remove these two conditions to prove 
Theorem~\ref{thm:main}.
We remove the condition that each critical point is non-degenerate
by using the complex bounds of \cite{ClvSTr}, see
Theorem~\ref{thm:complex bounds}. The condition on the number
of solenoidal attractors is removed through a technical perturbation 
argument, Lemma~\ref{lem:perturb}.

While we do not focus on renormalization in this paper,
let us point out that by now it is not difficult to
remove the condition that all critical points are non-degenerate
from \cite{Sm2}.
 McMullen, \cite{McM2}, proved exponential convergence of 
renormalization acting on quadratic-like mappings,
which are infinitely renormalizable of bounded type.
This was extended to multimodal mappings
with quadratic critical points by Smania, \cite{Sm2}.
From the complex bounds of \cite{ClvSTr} and
the quasiconformal rigidity of analytic mappings,
\cite{ClvS}, it is possible to extend this proof to 
infinitely renormalizable mappings of bounded type in 
$\mathcal{A}_{\underline b}(I)$.
Let us mention that using the {\em decomposition of 
a renormalization} and building on the exponential convergence of
renormalization of analytic mappings, exponential convergence of
renormalization for $\mathcal C^k,$ $k\geq 3,$ symmetric unimodal
mappings, in the $\mathcal C^k$-topology, was proved in
\cite{AMdM}.
Renormalization ideas figure heavily in our proof;
however, we leave the investigation of 
the rate of convergence of
renormalization (of, in particular, smooth mappings) 
to future work. 

We believe that the methods used in this paper can be improved on
to extend Theorem~\ref{thm:main smooth} to
$\mathcal C^2$ mappings with critical points of integer orders;
however, developing these tools (in particular proving the complex bounds for
these mappings) would take us far from
the goal of this paper. Since our proof of 
Theorem~\ref{thm:II smooth} depends on 
hyperbolicity of the quadratic fixed point of renormalization,
extending this result to mappings with lower regularity
would require a different approach. 
Let us also remark that our methods depend heavily on complex tools, 
so we do not obtain results for mappings with flat critical points
or with critical points of non-integer order.

\subsection{Outline of the paper}
In Section~\ref{sec:prelim}, we state some basic definitions
which will be used throughout this paper, and give the necessary 
background in real dynamics. 
In Section~\ref{sec:entropy and renormalization}, to make this paper
more self-contained,
we reduce Theorem~\ref{thm:main} to an equivalent statement
about infinitely renormalizable mappings with zero entropy, 
Theorem~\ref{thm:main bis}.
In Section~\ref{sec:spaces}, we introduce the different spaces of 
mappings in which we will work.

Of particular importance to us is the space
of stunted sawtooth mappings, $\mathcal S$, see Section~\ref{subsec:ssm}.
Stunted sawtooth mappings were introduced in \cite{MT}.
From a combinatorial point of view, they model 
mappings of the interval with finitely many critical points well.
Moreover, the space of stunted sawtooth mappings is a 
convenient space of mappings to work in since, 
in this space, entropy is monotone in each of the parameters
(the ``signed heights'' of the plateaus). Indeed the analogue of
Theorem~\ref{thm:main} is known in this space (see Section~\ref{subsec:ssm} for
the necessary terminology):

\begin{thm}\cite{TH}\label{thm:boundary ssm}
Let  $T_{\xi}\in \mathcal S$ be so that $\mathcal Per(T_{\xi})=\{2^n:n\in \N\cup\{0\}\}.$ Given $\nu>0$ there exist $\alpha,\beta \in [-e,e]^m$ so that $|\xi-i|<\nu$ for $i=\alpha, \beta,$ where $h(T_{\alpha})>0$ and $T_{\beta}$ has only finitely many periods.  
\end{thm}

This
result is the starting point for the results of this paper.
In Section~\ref{sec:boundary},
we will transfer it successively to the space of
polynomials using ideas from \cite{BruvS}, then via the 
Douady-Hubbard Straightening Theorem
to polynomial-like mappings,
and finally to analytic mappings with even critical points
via renormalization and specifically the complex bounds of
\cite{ClvSTr}. Using the transversal non-singularity of the derivative of
the renormalization operator acting from the space of 
analytic mappings to the space of polynomial-like germs,
we go on to prove Theorems~\ref{thm:main} and \ref{thm:II}.
We obtain 
Theorem~\ref{thm:main smooth}
from Theorem~\ref{thm:main}
via an approximation argument, which is similar to one used in
\cite{GMdM}. 
Once we have proved Theorem~\ref{thm:main smooth},
we use it together with results of \cite{Davie}
on the hyperbolicity of the period-doubling
renormalization operator acting on smooth unimodal mappings
to prove Theorem~\ref{thm:II smooth}. Finally we 
deduce Theorem~\ref{thm:III}.

\subsection{Standing assumptions}
Unless otherwise stated,  we will assume the following:
\begin{itemize}
  \item The vector of orders of critical points
$\underline b=(\ell_1, \dots ,\ell_b)$ is a vector positive even
integers.
\item All renormalizations are period-doubling, see \pageref{page:period2}.

\end{itemize}

\section{Preliminaries}
\label{sec:prelim}

\subsection{Notation and terminology}

Given a topological space $X$ and $A\subset X$ we denote the 
\emph{boundary} of $A$ by $\partial A$ and its 
\emph{closure} by $\text{cl}(A).$ 
If $X$ is a metric space, we denote the 
\emph{open ball of radius $\varepsilon$} centred around 
$x\in X$ by $ B_\varepsilon(x)=\{y\in X:\dist(x,y)<\varepsilon\}.$ 

As usual, $\mathbb R$ and $\mathbb C$ denote the
real line and the complex plane, respectively, and 
$I$ will always denote a compact interval in 
$\mathbb R.$ It will be convenient to assume that
$I=[-1,1].$ We denote the circle $\mathbb R\ \mod 1$
by $\mathbb T$. If $X$ is a set and $x\in X$, we let 
$\mathrm{Comp}_x(X)$ be the connected component
of $X$ containing $x$. For $\varepsilon>0,$ we let
$\mathbb D_{\varepsilon}$ denote the disk of radius
$\varepsilon$ centred at the origin.

Given a continuous
piecewise monotone map $f\colon I\to I,$ we call its local extrema
\emph{turning points.} If $f$ has finitely many
turning points %all of which are in the interior of $I$ 
and
$f(\partial I)\subset \partial I,$ then $f$ is called a
\emph{multimodal map.}
The images 
of the turning points of a multimodal mapping
are called \emph{critical values.}

\subsection{Background in dynamics}

Given a function $f\colon X\to X$ acting on a topological space 
$X,$ the \emph{orbit} of a point $x\in X$ 
is defined as the set $\mathcal O_f(x)=\{f^n(x): n\in \N\}.$ 
The set of accumulation points of $\mathcal O_f(x)$ is known as the 
\emph{$\omega$-limit set} of $x$ and is denoted by $\omega(x).$  
A point $x\in X$ is called \emph{non-wandering} 
if given any open set $U\ni x,$ there exists $n\in \N$ such that 
$f^n(U)\cap U\neq \emptyset.$ 
The set of non-wandering points of a map $f$ will be denoted by 
$\Omega(f).$ In particular, if $x\in \omega(x),$ 
then we say that $x$ is \emph{recurrent}.  

\begin{defn}
In a space of mappings $\mathcal X,$ we let 
$\Gamma_{\mathcal X}$ denote the subset of
$\mathcal X$ consisting of mappings $f$ with 
$\mathcal Per(f)=\{2^n:n\in\mathbb N\cup\{0\}\}$.
 When it will not cause confusion, we omit
$\mathcal X$ from the notation.
\end{defn}

Given a piecewise monotone map 
$f\colon I\to I$ with $m$-turning points $-1< c_1<\ldots < c_m<1,$ 
we denote by $i_f(x)$ the \emph{itinerary of x}
and the \emph{kneading sequence} of $c_j$ by
$$\nu_j:= \lim_{x\downarrow c_j} i_f(x),$$
where the sequence $\nu_j$ 
consists of the symbols $I_0,\ldots, I_m,$ 
where the $I_i$'s are the intervals from 
$I\setminus \{c_1,\ldots, c_m\}.$  Finally, we denote by
$$\nu(f)=(\nu_1, \ldots , \nu_m)$$
the \emph{kneading invariant of f}.
See \cite[Section II.3]{dMvS} for the definition of the
itinerary of a point.

% We say that a multimodal mapping is
% {\em critically finite} if the union of the set of orbits
% of its turning points is a finite set of points.

\begin{figure}[h]
\centering{
\resizebox{50mm}{!}{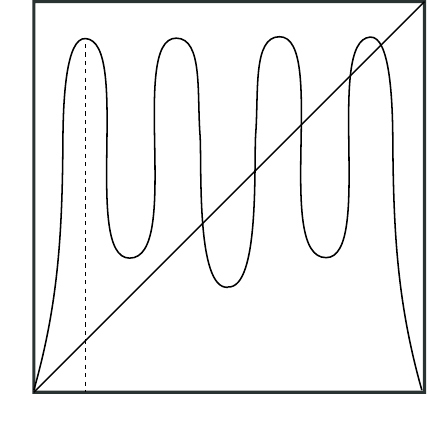}
\captionsetup{labelformat=empty}
\caption{Figure 1}
\label{fig:shape}
}
\end{figure}

\begin{defn}
Let $f\colon I\to I$ be a piecewise monotone map with turning points $-1<c_1<c_2<\ldots < c_{m}<1$  and critical values $v_1<v_2<\ldots <v_{s}.$ Then, we define its \emph{shape} as the set:
$$\tau=\{ (i,j_i): 1\leq i\leq m\},$$ 
where  $j_i\in \{1,\ldots, s\}$ is so that $f(c_i)=v_{j_i}.$ %and $k_i$ is the order of the critical point $c_i.$
\end{defn}
For example, the map in Figure \ref{fig:shape} has shape
$$\tau=\{ (1,3), (2,2), (3,3), (4,1), (5,3), (6,2), (7,3)\}.$$
The shape keeps track of the linear order of the critical values in the real line, which critical
points have which critical values and, in particular,
which critical points have the
same critical values. This notion of shape is useful in the study of
mappings arising as renormalizations. Since such mappings are
compositions of unimodal mappings, they have more ``symmetries''
than general polynomials.

 A \emph{shape}  is a set of ordered pairs 
$$\tau=\{(i,j_i)\in\{1,2,\dots, m\}\times\{1,2,\dots,s\}:i\in\{1,2,\dots,m\}\}$$
with the property that the mapping 
$i\mapsto j_i$ from $\{1,2,\dots,m\}$ to $\{1,2,\dots,s\}$
is onto.

% If $\tau$ is a shape, we let
% $\mathcal{PM}(\tau)$ denote the set of piecewise
% monotone mappings with shape $\tau$ and
% ${\mathcal Poly}(\tau)$ denote the space of
% polynomials with shape $\tau$. 

%We will need to refine the notions of 
%kneading invariant and
%shape slightly to accomodate mappings with critical points of 
%odd order. From the real point of view, these critical points are
%topologically ``invisible''. However, we are going to make use of
%complex tools and we will need to keep track of odd critical points.
%Given a differentiable mapping $f:I\to I$ with 
%$n$-critical points $-1< c_1<\ldots < c_n<1,$
%the $*$-\emph{kneading invariant of } $f$ is 
%ordered set of itineraries
%$\nu(f)=(\nu_1, \ldots , \nu_n),$
%defined as in the definition of the kneading invariant
%(only taking into account the turning points in the definition
%of each itinerary). The $*$-\emph{shape} of a 
%piecewise monotone mapping $f$ 
%is the shape $\tau(f)$ together with a finite list of marked 
%itineraries. In our applications to piecewise monotone mappings, 
%the marked itineraries will correspond to unique points.
%The $*$-\emph{shape of a polynomial} $p$ is the shape
%$\tau(p)$ together with the itineraries of its odd critical points.

Given a multimodal map $f$ and a forward invariant set $A\subset I,$ we say that $A$ is a \emph{topological attractor}
if its basin $\mathcal B(A)= \{x :\omega (x) \subset A\}$ satisfies the following properties:
\begin{itemize}
\item $\mathcal B(A)$ contains a residual subset of an open subset of $I$,
\item there exists no closed forward invariant set $A'\subset A$ with $A'\neq A$
  for which $\mathcal B(A')$ and $\mathcal B(A)$ coincide up to a meager set.
\end{itemize}

\subsubsection{Renormalization}

\begin{defn}
Let $f\colon I\to I$ be an interval map and let $n\in \N.$ A proper subinterval $J\subset I$ is called a \emph{restrictive interval} of period $n$ if:
\begin{itemize}
\item the interiors of $J,\ldots, f^{n-1}(J)$ are pairwise disjoint;
\item $f^n(J)\subset J$ and $f^n(\partial J)\subset \partial J$;
\item at least one of the intervals $J,\ldots f^{n-1}(J)$ contains a turning point of $f$;
\item $J$ is maximal with respect to these properties.
\end{itemize}

If $f$ has a restrictive interval $J$ of period $n\geq 2$, we say that $f$ is
\textit{renormalizable}. Let us assume further that $J$ is the largest such
restrictive interval (that is, the restrictive interval with the smallest period $\geq 2$) about a turning point $c$.
If $\Phi\colon J\to I$ is an affine 
surjection, the {\em renormalization operator},
$f \to \mathcal R_{c}(f),$ 
is defined by
$$\mathcal R_c(f) =  \Phi \circ f^n\circ \Phi^{-1}\colon I \to I,$$
and $\mathcal R_c(f)$
is known as a \emph{renormalization of $f$.} When it will not cause confusion,
we may omit $c$ from the notation.
\end{defn}

Assume $f$ possesses infinitely many restrictive intervals $J_n\ni c,$ then we
say that $f$ is \emph{infinitely renormalizable} at $c.$ Let  $q_n$ denote the
period of $J_n.$
Under these circumstances
% if the sequence $J_n$ contains all restrictive intervals around $c$ the following holds (see % \cite[Theorem 4.1] \cite{dMvS}). 
the set $\omega(c)$ is a solenoidal attractor $L$ with
$$ L=\bigcap_{n=1}^\infty L_n \, \,\, \,\text{ where } \, \, \, \,L_n=\bigcup_{k=0}^{q_n-1} f^k(J_n).$$

For a proof see Theorem III.4.1 in \cite{dMvS}.
%We will denote infinitely renormalizable maps with zero topological entropy as \emph{Feigenbaum maps}. 

Suppose that $f$ is infinitely renormalizable at a turning point $c$,
and let $\{q_n\}_{n=1}^\infty, q_n\in\mathbb N,$ be the strictly increasing 
sequence  such that for each $n\in\mathbb N,$
$f$ has a restrictive interval $J_n\owns c$ of period $q_n,$
and for any $q\in\mathbb N\setminus\{q_n\}_{n=1}^\infty$, there is no
restrictive interval of period $q$ about $c.$
We say that $f$ has {\em bounded combinatorics} at $c$ if there
exists $M\in\mathbb N$ so that $q_{n+1}/q_n\leq M$ for all $n$.
A mapping is infinitely renormalizable with
\label{page:period2} period-doubling
combinatorics at $c$ if and only if $q_{n+1}/q_n=2$ for all $n$.
Taking $\Phi_n\colon J_n\to I$ to be the affine surjection from $J_n$ onto $I,$
we define the  $n${\em -th renormalization of $f$ at $c$}:
$$\mathcal R^n_c(f)=\Phi_n\circ f^{q_n}\circ\Phi_n^{-1}\colon I\to I.$$

Any $\mathcal C^1$ mapping of the interval is semi-conjugate to a polynomial
mapping \cite[Theorem II.6.4]{dMvS}. The
semi-conjugacy collapses to points wandering intervals and attractors whose basins do not
contain turning points.
The mappings under consideration in this paper do not have wandering
intervals, but they may have such attractors. We say that two interval
mappings $f$ and $g$ have the \emph{same combinatorics} if there exists a
polynomial to which both $f$ and $g$ are semi-conjugate, and if corresponding
critical points of $f$ and $g$ have the same orders. Note that this definition
is more restrictive than just asking for $f$ and $g$ to have the same kneading
invariant, but we need to require corresponding critical points to have the same
orders to make use of complex extensions of $f$ and $g.$

The following result makes use of restrictive intervals to 
decompose the non-wandering set of $f$, denoted by $\Omega(f).$
\begin{thm}\cite[Theorem III.4.2]{dMvS}\label{omega limit}
Given a multimodal mapping $f$, there exists $N\in \N\cup\{\infty \}$ such that the following holds.
\begin{enumerate}
\item $\Omega(f)$ can be decomposed into closed forward invariant subsets $\Omega_n$:
$$\Omega(f)= \bigcup_{n\leq N} \Omega_n,$$
where the set $\Omega_n$ is defined as follows. Let $K_0=I$ and let $K_{n+1}$ be the union of all  maximal restrictive intervals of $f|_{K_n}.$ Then $K_n$ is a  decreasing  sequence of  nested sets, each consisting of a finite union of intervals for each finite $n\leq N.$ Then
$$\Omega_n:= \Omega(f) \cap cl(K_n\setminus K_{n+1})$$
for $n<N$ and $\Omega_N=\Omega(f)\cap K_N.$  If $N=\infty,$ we define $K_\infty= \cap_{n=0}^\infty K_n.$
\item For each finite $n\leq N,$ the set  $\Omega_n$ is a union of transitive sets. If $N=\infty,$ we have that  $\Omega_\infty=K_\infty$ is a union of solenoidal attractors.
\item The map $f$ has zero entropy if and only if $\Omega_n$ consists of periodic orbits of period $2^n$ for every finite $n\leq N.$
\end{enumerate}

\end{thm}

Theorem~\ref{omega limit} implies that the attractors of maps  in
$\Gamma$ can only be periodic or solenoidal. In the latter case the
attractor is equal to $\omega(c),$ where $c$ is some turning point
at which $f$ is infinitely renormalizable.

\subsection{Analytic and smooth mappings}
Given $a>0$, let
$\Omega_a=\{z\in\mathbb C:\dist(z,I)<a\}.$
We let $\mathcal{B}_{\Omega_a}$ denote the
space of complex-analytic mappings on
$\Omega_a$
which are continuous on $\text{cl}(\Omega_a).$
We endow $\mathcal{B}_{\Omega_a}$ with the sup-norm.
We let $\mathcal{B}_{\Omega_a}^{\mathbb R}$ denote the subspace of mappings 
in $\mathcal{B}_{\Omega_a}$
which commute with complex conjugation, and 
call such mappings {\em real}.

Given $k\in \N\cup\{0,\infty\}$ we let  $\mathcal C^{k}(I)$ denote the space of
$\mathcal C^k$ multimodal 
maps
of the compact interval $I;$ 
\emph{i.e.} continuous maps which are $k$-times differentiable with
continuous $k$-th derivative
on some small (real) neighbourhood of $I$.
We endow $\mathcal C^k(I)$ with the usual norm:
$$\|f\|_{\mathcal C^k(I)}=\max_{0\leq i\leq k}\;\sup_{x\in I}|f^{(i)}(x)|,$$
where $f^{(i)}$ denotes the $i$-th derivative of $f$.

We let $\mathcal C^\omega(I)$ denote the space of
real-analytic functions on $I$. 
We endow 
$\mathcal C^{\omega}(I)$ with a topology defined as follows:
We say that a net $\{f_\alpha\}$ converges 
to $f$ if all the $f_{\alpha}$ are analytic on some
fixed neighbourhood 
$\Omega$ of $I$ and $f_\alpha$ converges uniformly to 
$f$ on every compact subset of $\Omega$.

Given $k\in\{1,2,\dots, \infty,\omega\}$ and
$\underline b=(\ell_1,\ldots,\ell_b)$
a vector of positive integers we say that 
$f\in\mathcal C^k(I)$ belongs to $\mathcal A^k_{\underline b}(I)$
\label{page:Ak}
if the following holds. 
The map $f$ has finitely many parabolic cycles and
$b$ critical points
$c_i, 1\leq i\leq b,$ labelled so that 
$c_1<c_2<\dots<c_b$, and
each $c_i$ has a 
neighbourhood on which we can express $f$ as
$$f(x)=\sigma_i\cdot \phi_i(x-c_{i})^{\ell_i}+f(c_{i}),$$ where $\phi_i$ is a local
$\mathcal C^{k}$ diffeomorphism with
$\phi_i(0)=0,$ $\sigma_i\in\{1,-1\},$ and $\ell_i\in \N$ is at least two.  
We say that $\ell_i$ is the {\em degree} or {\em order} of $c_i$.
If $\ell_i$ is 
even, we say that the corresponding critical point $c_i$ has {\em even
order}. 
We let $\mathrm{Crit}(f)$ denote the set of critical points of $f$.

For many of the results in real dynamics that we will
recall later, the condition that the critical points 
have integer order is unnecessary. 
%and the condition on the Taylor
%series at parabolic points are unnecessary. 
The results which use complex analysis require this condition on
the local behaviour of the critical points. %and the condition on the Tayor
%series is required for Theorem~\ref{thm:qs rigidity}.

%Moreover, we assume that 
%at each parabolic cycle
%$\{p_0,p_1,\dots ,p_{s-1}\}$ of period $s$
%there exists $d$ so that
%the Taylor expansion
%of $f^s$ at $p$ is given by
%$$f^{s}(x)=p+\lambda(x-p)+a(z-p)^{d+1}+R(|x-p|),$$
%where $R=o(|x-p|^{d+1})$.
%This condition excludes the possibility of $f$ having
%infinitely many parabolic points, since any limit point of 
%such parabolic points would be tangent to the identity with infinite
%order.
We will denote by $\mathcal A_{\underline b}(I)$
the set of analytic maps $\mathcal A^{\omega}_{\underline b}(I).$
%Of course for analytic mappings, the condition on the 
%Taylor series is redundant.
When it will not cause confusion we will drop the 
subscript $\underline b$ from the notation.
%For $M, b\in\mathbb N$, with $M\geq 2,$ we let 
%$\mathcal{A}^k_{M,b}(I)=\cup_{\underline b\in\{2,3,\dots, M\}^b}
%\mathcal A_{\underline b}^k(I)$.
For $b\in\mathbb N$,
let \label{page:Aeven} $\mathcal A^k_{even,b}(I)=\cup_{\underline b}\mathcal
A^k_{\underline b}(I),$
where the union is taken over $b$-tuples,
$\underline b,$ with all entries even.

We say that a mapping $f$ is {\em critically finite} when its post-critical set
$$\{f^i(c):c\in\mathrm{Crit}(f), i\in\mathbb N\}$$
is a finite set. A mapping is critically finite if and only if
all of its critical points are periodic or 
pre-periodic.

\subsubsection{Real bounds}
Real {\em a priori} bounds, were first proved
for unimodal infinitely renormalizable mappings by Sullivan, \cite{Sullivan,dMvS}.
For multimodal mappings 
with all critical points even, 
real bounds were obtained in
\cite{Sm1} for
infinitely renormalizable mappings.
These were generalized in \cite{Shen-C^2}.
We have the following
 real bounds for infinitely renormalizable mappings.
\begin{thm}[Real Bounds, \cite{ClvS, ClvSTr}]
\label{thm:real bounds}
There exists $\delta>0$ so that the following holds.
Suppose that $f\in\mathcal A^3_{\underline b}(I)$ is infinitely renormalizable 
at a critical point $c$, suppose that 
$$J_1\supset J_2\supset J_3\supset\dots$$ is a sequence 
of restrictive intervals about $c$. Then for
$n$ sufficiently large,
if $f^s\colon J\to J_{n}$ is a diffeomorphism, we have that
there exists an interval $\hat J\supset J$ so that
$f^s\colon \hat J\to (1+\delta)J_{n}$ is a diffeomorphism.
Moreover, $(1+\delta)J_{n+1}\subset J_{n}$.
\end{thm}

\subsubsection{Asymptotically holomorphic mappings}
Asymptotically holomorphic mappings have proved to be
vital in extending known results for analytic mappings to the 
case of smooth mappings. One of their first uses in 
dynamical systems was in a proof of 
rigidity of quadratic Fibonacci mappings, \cite{L-Fib}.
We will make use of a particularly effective 
asymptotically holomorphic extension given by
\cite{GSS}.
These extensions have been used to study smooth 
mappings of the interval, \cite{ClvSTr, CdFvS,ClvS},
and on the circle, \cite{GdM,GMdM}.

Suppose that $J\subset\mathbb R,$ and that $U$
is an open subset of $\mathbb C$ containing $J$.
We say that a $\mathcal C^1$  mapping $f\colon U\to\mathbb C$ is 
{\em asymptotically holomorphic of order} $k$ on $J$
if 
$$ \frac{\partial f}{\partial \bar z}(x)=0 \mbox{ for } x\in J,\quad
\mbox{and}\quad
 \frac{\frac{\partial f}{\partial \bar z}(x+iy)}{|y|^{k-1}}\rightarrow
 0\;
\mbox{as}\;
y\rightarrow 0.$$

Let $\kappa\geq 1$ and let $U\in\mathbb C$ be an open set.
We say that a mapping $f\colon U\to\mathbb C$ is 
$\kappa$-{\em quasiregular} if it is orientation preserving, with
locally square integrable distributional derivatives, $f_z$ and $f_{\bar z}$,
which satisfy 
$$\max_{\alpha}|\partial_{\alpha}f(z)|\leq 
\kappa\min_{\alpha}|\partial_{\alpha}f(z)|,$$
for almost every $z\in U,$ where
$$\partial_{\alpha}f(z)=\cos(\alpha)f_x(z)+\sin(\alpha)f_y(z),\;\mbox{for}\;
\alpha\in[0,2\pi).$$
We say that
$f$ is {\em quasiregular} if it is $\kappa$-quasiregular for
some $\kappa\geq 1$.

\begin{thm}\cite{GSS}\label{thm:GSS}
Suppose that $f\in \mathcal C^k(I)$, then 
there is an asymptotically holomorphic extension of order $k$ of $f$ to a
neighbourhood of the interval in the complex plane.
\end{thm}

\subsubsection{Smooth polynomial-like mappings}

A {\em polynomial-like mapping} is 
a proper holomorphic branched covering map $f\colon U\to V,$ 
where $U\Subset V\neq\mathbb C$ are two simply connected 
complex domains. 
We will consider polynomial-like mappings up to
affine conjugacy.
We define the 
{\em filled Julia set} for a polynomial-like map $f$ as:
$$K(f)= \bigcap_{n\in \N} f^{-n}(V).$$
The {\em Julia set} of $f$, denoted by $J(f)$, is the boundary of $K(f).$
We say that a polynomial-like mapping $f\colon U\to V$ is {\em real} if 
$U$ and $V$ are real-symmetric 
and $f$ commutes with complex conjugation.

We say that two polynomial-like mappings $f\colon U_f\to V_f$ and 
$g\colon U_g\to V_g$ are \emph{quasiconformally equivalent} 
if there exists a qc-mapping $H$ defined on a neighbourhood $W$ 
of $K(f)$ to a neighbourhood of $K(g)$ such that 
$H\circ f(z)=g\circ H(z),$ $z\in W.$ 
If additionally we have that
$\bar \partial H=0$ on $K(f)$, then we say that
$f$ and $g$ are {\em hybrid equivalent}\label{def:hybrid eq}.

An \emph{asymptotically holomorphic polynomial-like mapping},
abbreviated \emph{AHPL-mapping}, of order $k$
is a proper $\mathcal C^k$ branched covering map $f\colon U\to V,$ 
where $U\Subset V\neq\mathbb C$ are two simply connected 
complex domains, which is asymptotically holomorphic of order
$k$ on $U\cap\mathbb R.$ Every such asymptotically holomorphic
polynomial-like mapping in this paper has the properties that
$U$ and $V$ are real-symmetric and $f$ commutes with
complex conjugation.

\subsubsection{Complex bounds}

Complex bounds for real mappings have a long history,
see the introduction of \cite{ClvSTr}. In the classes of mappings
most relevant to us, they were first proved for 
real-analytic, infinitely renormalizable unimodal mappings with
bounded combinatorics by Sullivan, \cite{Sullivan}.
This result was extended to analytic multimodal mappings with all
critical points even by Smania \cite{Sm1}. The authors together with
van Strien proved the following, which built on work of
\cite{Shen-C^2, KSS-rigidity}.

\begin{thm}\cite{ClvSTr}\label{thm:complex bounds}
Suppose that $f\in\mathcal A^k(I),$ $k\geq 3,$
is infinitely renormalizable
at an even critical point $c_0$. Let 
$J_1\supset J_2\supset\dots$ denote the sequence of restrictive
intervals for $f$ about $c_0,$ where the period of $J_i$ is $q_i$.
Then for all $i$ sufficiently large, there exists an 
asymptotically holomorphic polynomial-like
mapping of order $k$,
$F\colon U\to V$ with $F=f^{q_i}|_U$,
 $U\supset J_i,$ $\mathrm{diam}(V)\asymp J_i,$ and $\mathrm{mod}
(V\setminus U)$ bounded away from zero.
Moreover, the dilatation of $F$ tends to zero as
$i$ tends to infinity.

If $f$ is analytic, then $F$ is a polynomial-like mapping.
\end{thm}
When $f$ is analytic, the polynomial-like mapping $F$ is 
is constructed using the holomorphic 
extension of $f$ to a neighbourhood of the interval.
If $f\in\mathcal A^k(I)$, the extension is constructed for any
$\mathcal C^k$ asymptotically holomorphic extension of order $k$ of $f$
to a neighbourhood of $I$. By Theorem~\ref{thm:GSS} at least one
such extension exists.

The following lemma is useful for
working with asymptotically holomorphic mappings.
\begin{lem}[Stoilow Factorization]
\label{lem:Stoilow}
If $f\colon U\to V$ is a quasiregular mapping, then we can factor $f$
as $f=h\circ \phi$ where $\phi\colon U\to U$ is quasiconformal and
$h\colon U\to V$ is analytic.
\end{lem}

\subsubsection{Quasisymmetric rigidity}
Quasisymmetric (quasiconformal) 
techniques were introduced into one-dimensional
dynamics by Dennis Sullivan, who observed that 
quasisymmetric rigidity of unimodal mappings
could be used to prove density of hyperbolicity.
Quasiconformal rigidity was first proved in
\cite{L-acta, GS} for quadratic polynomials. 
It was later proved for real polynomials 
with all critical points even and real
in \cite{KSS-rigidity}. The first author together with van
Strien proved the following:

\begin{thm}\cite{ClvS}\label{thm:qs rigidity} For mappings
of the interval we have the following:
\begin{enumerate}
\item {\em Rigidity for analytic mappings.}
Suppose that $f,\tilde f\in\mathcal A_{\underline b}(I)$, 
are topologically conjugate mappings 
by a conjugacy
which is a bijection on
\begin{itemize}
\item the sets of parabolic points,
\item the sets of critical points and corresponding critical points 
have the same order.
\end{itemize}
Then $f$ and $\tilde f$ are quasisymmetrically conjugate.

\item {\em Rigidity for smooth mappings.} Suppose that $f,\tilde f\in\mathcal A^k_{\underline b}(I)$, 
$k\geq 3,$ do not have parabolic cycles, and that they
are topologically conjugate mappings 
by a conjugacy
which is a bijection on
the sets of critical points and corresponding critical points 
have the same order.
Then $f$ and $\tilde f$ are quasisymmetrically conjugate.
\end{enumerate}
\end{thm}

\medskip
\noindent\textit{Remark.} If $f$ and $\tilde f$ are 
$\mathcal C^k$ we can allow for parabolic points as in the 
theorem for analytic mappings
under some additional regularity assumptions, see \cite{ClvS}.
We will only apply rigidity of smooth mappings to deep
renormalizations, which do not have parabolic cycles by 
\cite[Theorem IV.B]{dMvS}.

\medskip

It is worth observing that since we are only concerned here with
infinitely renormalizable mappings with bounded geometry, we do not require the
full result of  Theorem~\ref{thm:qs rigidity}. Indeed, the following result is
sufficient:

For any $b\in\mathbb N,$ let $\mathcal G_{\underline b}$
be the collection of $\mathcal C^3$ maps
$$f\colon \big(\bigcup_{j=0}^m J_j \big)\cup \big(\bigcup_{j=0}^{b-1}I_i\big)\to\bigcup_{j=0}^{b-1}I_i$$
with the following properties:
\begin{itemize}
\item $I_i$’s are open intervals with pairwise disjoint closures;
\item $m$ is a non-negative integer;
\item each $J_j$ is an open interval contained in $I_0$, and the $J_j$'s have
  pairwise disjoint closures contained in $I_0,$ unless $m=0$ in which case we also allow $J_0=I_0;$
\item $f$ is a proper map;
\item $f$ extends to a $\mathcal C^3$ map defined on the closure of its domain such that $f'$ does not vanish at the boundary;
\item for each $1\leq j\leq m$ $f|_{J_j}$ is a diffeomorphism;
\item for any $U\in\{J_0, I_1,\dots,I_{b-1}\},$ $f|_U$ has a unique critical point $c_U$;
\item all the critical points do not escape under forward iterates of $f$ ,
\item all the critical points are non-periodic and recurrent, and they have the
same $\omega$-limit set which is a minimal set.
\item the extension of $f$ to the closure of its domain has only hyperbolic repelling periodic points.
\end{itemize}
The following proposition follows from the much more general \cite[Theorem 2]{Shen-C^2}.
\begin{thm}
  Suppose that $f, \tilde f\in \mathcal G_{\underline b}$ are two combinatorially
  equivalent maps with the property that each of their critical points is
  infinitely renormalizable with period doubling combinatorics.  Then they are quasisymmetrically conjugate on the postcritical sets, i.e., the combinatorial equivalence can be realized by a quasisymmetric map.
\end{thm}

For real polynomials we have the following:

\begin{thm}\cite{KSS-rigidity, ClvSTr}\label{thm:rigidity}
Suppose that $f$ and $\tilde f$ are
two real polynomials, with real critical points.
Assume that $f$ and $\tilde f$
are topologically conjugate as dynamical systems on the real line,
that corresponding critical points for $f$ and
 $\tilde f$ have the same order and that parabolic
points correspond to parabolic points, then $f$ and
 $\tilde f$ are quasiconformally conjugate as
dynamical systems on the complex plane.
\end{thm}
This result was proved for mappings with all critical points of
even order in \cite{KSS-rigidity}, and this restriction on the
degrees of the critical points was removed in \cite{ClvSTr}.

\subsubsection{Absence of invariant line fields}
A {\em line field} on a subset $E$ of $\mathbb C$ is a choice of
a line through the origin in the tangent space
$T_eX$ at each point $e\in E$.
For a polynomial, 
absence of invariant line fields on the Julia set
is an ergodic property of the dynamics, which 
is closely related to rigidity \cite{McM}.
Complex bounds are a key tool in the proof of quasisymmetric
rigidity, and they play a crucial role in establishing the 
absence of invariant line fields for polynomials.
Absence of invariant line fields were first proved in
\cite{McM}. Building on this, they were proved for real
infinitely renormalizable polynomial-like mappings in
\cite{Sm2}, and for real rational maps with all critical points real
and with even degrees in
\cite{Shen}.

\medskip
\noindent\textit{Remark \cite{McM}.}
\label{page:bdlf}
A line field may be identified with 
a Beltrami differential $\mu=\mu(z)\frac{d\bar z}{dz}$ with
$|\mu|=1$: The real line through $v=a(z)\frac{\partial}{\partial z}$ 
corresponds to the Beltrami differential $\frac{a}{\bar a}\frac{d\bar
  z}{dz}.$ Conversely a Beltrami differential determines a function
$\mu(v)=\mu(z)\frac{\bar a(z)}{a(z)}$, where 
$v=a(z)\frac{\partial}{\partial z}$ is a tangent vector;
the line field consists of those tangent vectors $v$ for which
$\mu(v)=1$.

\medskip

We will make use of the following theorems about
polynomials:

\begin{thm}\cite{McM, Shen, ClvSTr}\label{thm:nilf}
Suppose that $f$ is a real polynomial with real critical points.
Then $f$ 
supports no measurable invariant line field on
its Julia set.
\end{thm}

Theorems \ref{thm:nilf} and \ref{thm:rigidity} together with the
B\"ottcher Theorem imply the following:
\begin{cor}\label{cor:rigidity}
Suppose that $f$ and $\tilde f$ are topologically conjugate mappings
as in the statement of the Theorem~\ref{thm:rigidity}
with connected Julia sets and  all periodic points
repelling. Then $f$ and $\tilde f$ are affinely conjugate.
\end{cor}

\section{Entropy and renormalization}
\label{sec:entropy and renormalization}
In this section, we study maps $f$ with $\mathcal Per(f)=\{2^n: n\in
\N\cup\{0\}\}.$ 
Our goal is 
to show that to prove Theorem \ref{thm:main},
it is enough to prove:
\begin{thmx}\label{thm:main bis}
{\em Every map $f\in \mathcal A^3_{\underline b}(I),$
which is infinitely renormalizable with
entropy zero can be approximated by mappings with positive
topological entropy and by mappings with finitely many periods.}
\end{thmx}
The equivalence of Theorems~\ref{thm:main} and 
\ref{thm:main bis} is not new, but we include it
to help make the
paper more self contained.

\medskip

There are many equivalent definitions of topological entropy. For
simplicity of exposition, we use the one introduced in
\cite{MiSz1}. Given a continuous piecewise monotone map  $f\colon I\to I$ we
define the \emph{lap number} of $f,$ denoted by $\ell(f),$ as the number of maximal intervals on which $f$ is monotone. 
%For example, if $f$ is a $\ell$-modal map, the number of turning points of $f$ is equal to $\ell(f) -1.$ 
The topological entropy is defined as the rate of exponential growth of $\ell(f^n).$ 

\begin{defn}
Given a continuous piecewise monotone map $f\colon I\to I$ we define its \emph{topological entropy}
as

$$h(f)= \lim_{n\to\infty}\frac{1}{n} \log(\ell(f^n)).$$

\end{defn}

For simplicity, we will refer to topological entropy as {\em entropy}. 
\medskip

The following classical result, see \cite{MiLl}, relates the entropy of a map with the periods of its periodic orbits.

\begin{prop}\label{prop:zero entropy per 2}
 A map $f\in \mathcal C^0(I)$ has positive entropy if and only if $f$ has a periodic orbit of a period which is not a power of two.
\end{prop}

\medskip

To get a characterization of the boundary of chaos we will take a closer look at
the level sets of the entropy map. Recall that
$$\Gamma_{\mathcal C^k(I)}= \{f\in \mathcal C^k(I):  \mathcal Per(f)=\{2^n :
n\in \N\cup\{0\}\},$$ and that we may omit the subscript on $\Gamma$ when it is clear
in what space we are working.

\begin{prop}\label{prop:positive open}
We have the following:
\begin{enumerate}
\item the set of maps with positive entropy is open
in the space $\mathcal C^k(I)$ for $k\geq 2,$ and
\item $\Gamma$ is closed in $\mathcal C^k(I)$, for $k\geq 1$.
\end{enumerate}
\end{prop}
\begin{pf}
The first statement follows from the fact that the topological entropy
is continuous on $\mathcal C^k(I),$ for $k\geq 2,$
see Theorem 6 in \cite{MiSz1}.
The second statement corresponds to Proposition 2.1 in \cite{Mi1}. 

\end{pf}

\noindent \textbf{Remark:} To show that  $\Gamma$ is a closed set the author in
\cite{Mi1} proves the following result which we alluded to before:
The set of maps  for which the set $\mathcal Per(f)$ is bounded is
$\mathcal C^k$-open for $k\geq 1$. 

\medskip

We have the following corollary.

 \begin{cor}\label{cor:boundary}
If $f\in \mathcal C^k(I)$ for $k\geq 2$ is  on the boundary of the set of maps with positive entropy, then  $f\in \Gamma.$ The same holds for maps which lie on the boundary of the interior of the set of maps with zero entropy.
\end{cor}

This corollary, also proved in \cite{BH}, provides a characterization of maps on
the boundary of chaos in $\mathcal C^k(I)$ for $k\geq 2$ which remains true for maps in $\mathcal A(I).$
% This means that to show Theorem~\ref{thm:main} we need to prove that
% all maps in $\Gamma$ can be approximated by maps in the interior of
% the set of mappings with
% zero entropy and by maps with positive entropy.

%\subsection{Renormalization}In this section we study the renormalization of interval maps.  Let us start by introducing some terminology.

\medskip

The next result will help us determine the combinatorics of renormalizable maps with zero entropy.

%\begin{thm}\label{entropy per}
%If $f$ is a multimodal map, then $$h(f)\leq p(f),$$
%
%where $p(f)= \limsup_{n\to \infty} \log P_n(f)/n$ and $P_n(f)$ is equal to the
%cardinality of the set of points of period $n$.
%\end{thm}
%
%For a proof of this refer to \cite{MiSz1,MiSz2}. As a Corollary of the previous Theorem we have the following.
%
%\begin{cor}\label{coro zero entropy}
%If $f$ is a multimodal map and the cardinality of $\mathcal P(f)$ is finite, then $h(f)=0$.
%\end{cor}

\begin{prop}\cite[Proposition III.4.2]{dMvS}\label{prop:restrictive intervals}
If $f\in \mathcal A^k(I)$ and $h(f)=0,$ then each restrictive interval is contained in a restrictive interval of period $2.$ Furthermore, every point in $I$ is either eventually mapped into a restrictive interval of period $2$, or is asymptotic to a fixed point.
\end{prop}

\begin{lem}\cite[Theorem 2]{TH}\label{lem:infinitely renorm}
If $f\in \Gamma,$ then $f$ is infinitely renormalizable. Furthermore,
if $J_n$ and $J_{n+1}$ are consecutive restrictive intervals (meaning
that $J_{n+1}$ is a maximal, with respect  to containment, proper restrictive interval in $J_n$), then the period of $J_{n+1}$ inside of $J_{n}$ is two.
\end{lem}
\begin{pf}
Consider $\Delta_j(f),$ the set of accumulation points of periodic orbits of periods greater or equal to $2^j,$ and let $$\Delta(f)=\bigcap_{j\in \N}\Delta_j(f).$$
It is clear from the definition that $\Delta(f)$ is closed and $f$-invariant. In addition, by Lemma 1 in \cite{TH} we know that no point in $\Delta(f)$ is periodic. %So $\Delta(f)$ contains no periodic points.
Proposition \ref{prop:zero entropy per 2}  and Proposition \ref{prop:restrictive
  intervals}, imply that  every point which is not eventually mapped into a
restrictive interval of period two is asymptotic to a fixed point. Given $p\in
\Delta(f)$ we have that the orbit of $p$ enters $J_0,$ a restrictive interval of
period two. By definition, there exists a turning point $c$ contained either in $J_0$ or in $f(J_0).$
Repeating the argument, substituting $f$ by $f^{2^n}$ for $n\in \N$ we
can find a nested sequence of restrictive intervals  $J_n,$  with
$J_{n+1}$ of period two under $f^{2^n}$ inside $J_n$ and such that $c\in J_n$. 
 %Since  $h(f)=0,$ $f$ is Feigenbaum.
\end{pf}

\begin{cor}
If $f\in \mathcal A^k(I)$ finitely renormalizable and $h(f)=0,$ then the period of its periodic orbits is bounded.
\end{cor}

\medskip
\noindent\textit{Proof that Theorem~\ref{thm:main bis} is equivalent to
  Theorem~\ref{thm:main}.}
Let us first suppose that Theorem~\ref{thm:main} holds,
and assume that $f\in\mathcal A_{\underline b}^3(I)$ is infinitely
renormalizable and has entropy zero. Since $f$ has entropy zero, by
Proposition~\ref{prop:restrictive intervals}, we have that each restrictive
interval of $f$ has period a power of $2,$ and since $f$ is infinitely
renormalizable, we have that $f$ has restrictive interval of each period $2^n,$
for $n\in\mathbb N.$ Thus, $\mathcal{P}er(f)=\{2^n:n\in\mathbb N\cup\{0\}\}.$ By
Theorem~\ref{thm:main}, $f$ can be approximated by mappings with positive
topological entropy and by mappings with finitely many periods.

Now, let us assume that Theorem~\ref{thm:main bis} holds, and
$\mathcal{P}er(f)=\{2^n:n\in\mathbb N\cup\{0\}\}.$
By Proposition~\ref{prop:zero entropy per 2}, we have that $f$ has zero entropy
and by Lemma~\ref{lem:infinitely renorm}, we have that $f$ is infinitely
renormalizable, so by Theorem~\ref{thm:main bis}, $f$ can be approximated by
mappings with positive topological entropy and by mappings with finitely many
periods. \qed

% From Corollary~\ref{cor:boundary} 
% we know that maps which can be approximated by 
% mappings with finitely many periods and by mappings
% with positive entropy 
% belong to $\Gamma,$ \emph{i.e.}
% that the set of periods of their periodic orbits 
% corresponds to the set $\{2^n: n\in \N\}.$
% By  Proposition~\ref{prop:zero entropy per 2} and
% Lemma~\ref{lem:infinitely renorm} we know that
% $\Gamma$ corresponds to the set of
% infinitely renormalizable maps with entropy zero.
% Hence Theorem~\ref{thm:main bis} implies
% Theorem~\ref{thm:main}. \qed

% The proof of this Theorem~\ref{thm:main bis}
% occupies the remaining part of the paper
% and will be divided into several steps.

\section{Spaces of mappings}
\label{sec:spaces}

\subsection{Stunted sawtooth mappings}\label{subsec:ssm}
In this section, we recall the definition of
stunted sawtooth mappings and collect some 
useful facts about these mappings.
% In particular, 
% we prove:
% Any 
% stunted sawtooth mapping $T$
% with ${\mathcal Per}(T)=\{2^n:n\in\mathbb N\cup\{0\}\}$ 
% can be approximated by 
% \begin{itemize}
% \item a stunted sawtooth mapping with all plateaus periodic and entropy zero, and 
% \item by a stunted sawtooth mapping with all plateaus periodic 
% and with positive entropy. 
% \end{itemize}
% see  Proposition~\ref{prop:bcssm}  below.

\subsubsection{The definition of the space of
stunted sawtooth mappings}
We start by defining an auxiliary piecewise linear mapping $S_0$,
which will be used in the definition of stunted sawtooth mappings.
The \emph{basic shape} of a piecewise linear mapping $S$
is defined as
$$\epsilon(S)=\left\{
\begin{array}{ll}
1 & \mbox{if } S \mbox{ is increasing at the left endpoint of } I,\\
-1 & \mbox{otherwise}.\\
\end{array}\right.
$$
 
\begin{figure}[h]
\centering{
\resizebox{40mm}{!}{\Large{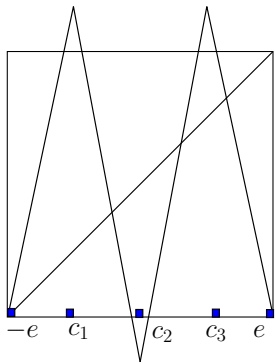}}
%\captionsetup{labelformat=empty}
\caption{The map $S_0$}
\label{fig:s_0}
}
\end{figure}

Let $\epsilon\in\{-1,1\}$.
Fix a constant $m\in\N,$ to be the number of
turning points,
and set $\lambda=m+2$. The slopes of the 
piecewise monotone mapping are either $\lambda$ or $-\lambda.$
Let $e=m\lambda/(\lambda -1),$ and set 
$A=[-e,e].$ 
One easily sees that there exists a unique $m$-modal
piecewise linear mapping $S_0$ with 
$\epsilon(S_0)=\epsilon$,
$m$ turning points,
$c_1,\ldots ,c_m$ at $-m+1, -m+3, \ldots, m-3, m-1$ with the following properties:
\begin{itemize}
\item $m+1$ intervals of monotonicity $I_0=[-e,c_1], I_1=[c_1,c_2], \ldots, I_m=[c_m,e].$
\item slopes $\pm \lambda$, extremal values $\pm \lambda;$ and
\item $S_0(\{-e,e\})\subset\{-e,e\}$.
\end{itemize}
See Figure~\ref{fig:s_0}.

\begin{figure}[h]
\centering{
\resizebox{40mm}{!}{\Large{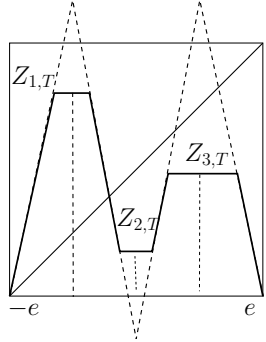}}
%\captionsetup{labelformat=empty}
\caption{The stunted sawtooth mapping with shape $\tau=\{(1,3),(2,1),(3,2)\}$ obtained from $S_0$.}
\label{fig:sstm1}
}
\end{figure}

The space of $\mathcal S=\mathcal S_{\epsilon,m}$
of \emph{stunted sawtooth maps} with
$m$ turning points consists of continuous maps 
$T$ with plateaus $Z_{i,T}$ with $i\in \{1,\ldots, m\}$ which are
obtained from the map 
$S_0$ (see Figure~\ref{fig:sstm1}) 
and satisfy the following:
 \begin{itemize}
 \item $Z_{i,T}$ is a closed symmetric interval around $c_i;$
 \item $T$ and $S_0$ agree outside $\cup_{i} Z_{i,T};$
 \item $T|_{Z_{i,T}}$ is constant and $T(Z_{i,T})\in [-e,e];$
 \item $Z_{i,T}$ have pairwise disjoint interiors.
 \end{itemize}

\begin{figure}[h]
\centering{
\resizebox{40mm}{!}{\Large{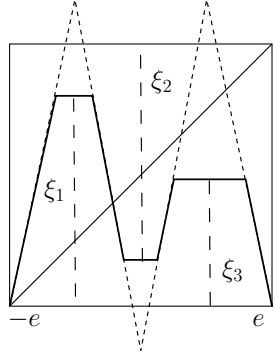}}
%\captionsetup{labelformat=empty}
\caption{The stunted sawtooth mapping parameterized by
$\xi=(\xi_1,\xi_2,\xi_3)$.}
\label{fig:sstm2}
}
\end{figure}

It is important to remark that a map $T\in\mathcal S$ could have touching plateaus. In other words, two of its plateaus could have one endpoint in common. In this case, we say that $T$ is $m$-{\em modal in the degenerate sense}. 
We use the $m$-signed extremal values $\xi\in [-e,e]^m$ to parameterize $\mathcal S$ in the following way.
\[
  \xi_i = \left\{
     \begin{array}{@{}l@{\thinspace}l}
       \, \, \,\, \, T(Z_{i,T})  & \text{if } c_i \text{ is a maximum of } S_0,
       \\
       -T(Z_{i,T})  & \text{if } c_i \text{ is a minimum of } S_0.
       \\
           \end{array}
   \right.
\]
Figure \ref{fig:sstm2} illustrates the parametrization.
We denote by $T_\xi$ the map in 
$\mathcal S$ with parameters $\xi=(\xi_1, \ldots \xi_m).$  

\medskip
We can identify the space $\mathcal S$ with the set
$$\{\xi=(\xi_1,\ldots \xi_m): \xi_i\in [-e,e] \text{ and } \xi_i\geq -\xi_{i+1}\}.$$
We define $T_{\xi}<T_{\xi'}$ if for the corresponding parameters
$\xi_i\leq \xi_i'$ for all $1\leq i\leq m$ with at least one strict inequality.

\medskip 

The definition of the shape of a stunted sawtooth mapping
is the same as for piecewise monotone mappings if we replace
turning points with plateaus.

\begin{defn} 
Given a map $T\in \mathcal S$ we will define its \emph{shape} in the
following way. 
Let $\ell\leq m$ be the number of distinct values of $T$ on the
plateaus $Z_{i,T}$, $1\leq i\leq m,$ and
label these values by $v_j,$ $1\leq j\leq \ell,$ so that 
$v_1<\ldots < v_{\ell}.$ 
The shape of $T$ is defined as the set of ordered pairs:
$$\tau(T)=\{ (i,j_i): 1\leq i\leq m\}$$
where $j_i$ is so that $T(Z_{i,T})=v_{j_i}.$ 
%The only difference between the shape of a map $P\in \mathcal P_{\underline n}$ with $m$ turning points and the shape of a map in $\mathcal S_m$ is that the vector associated to maps in $\mathcal P_{\underline n}$ has the extra information of the order of the critical points, since maps in $\mathcal S_m$ have plateaus, instead of turning points, this information is irrelevant in this case. 
For example, the shape of the map in Figure \ref{fig:sstm1} is 
$\tau(T)=\{(1,3),(2,1),(3,2)\}.$

Given a map $T\in \mathcal S$ with shape $\tau$ we define
$$\mathcal S({\tau})=\{ T'\in \mathcal S: \text{the shape of $T'$ is equal to } \tau\}.$$  
\end{defn}

\medskip

Let us recall two useful facts related to the entropy of
stunted sawtooth mappings.

\begin{prop}\cite[Proposition 4.1]{BruvS}\label{entropy monotone}
The map $\xi=(\xi_1,\ldots, \xi_m)\to h(T_\xi)$ is non-decreasing in each coordinate.
\end{prop}

The following is a slight variation of the
Theorem of Hu-Tresser.

\begin{prop}\cite[{\it c.f.} Theorem 1]{TH}\label{prop:bcssm}
Suppose that $T_{\xi}$ is a stunted sawtooth mapping
with ${\mathcal Per}(T_{\xi})=\{2^n:n\in\mathbb N\cup\{0\}\}.$
Let $\tau=\tau(T_{\xi}).$
Then for any $\varepsilon>0,$ there exist,
$\xi',\xi''$ with $|\xi-\xi'|<\varepsilon,$
 $|\xi-\xi''|<\varepsilon,$ such that 
\begin{enumerate}
\item $\tau(T_{\xi})=\tau(T_{\xi'})=\tau(T_{\xi''}),$
\item both $T_{\xi'}$ and $T_{\xi''}$ have all plateaus periodic or
pre-periodic,
\item $T_{\xi'}$ has positive entropy, and
\item $T_{\xi''}$ has finitely many periods.
\end{enumerate} 
\end{prop}
Since Proposition~\ref{prop:bcssm} 
can be proved using the same perturbative argument used to obtain 
Theorem 1 of \cite{TH} (to additionally obtain conclusions (1) and (2) in the
statement), we omit the proof.
% \begin{pf}
% We can apply Theorem~\ref{thm:boundary ssm} to find $\alpha$ and
% $\beta$ with $T_{\alpha}, T_\beta\in \mathcal S_{\epsilon,m},$ 
% $h(T_{\alpha})>0$
% and so that $T_{\beta}$ has only finitely many periods.
% %Since none of
% %the plateaus is periodic or 
% %pre-periodic and their orbits are disjoint from 
% %$W(T)= \bigcup T^{-k}(\mathrm{int}(\cup_{i=1}^mZ_{i,T}))$ we get the
% %following. 
% %All of the periodic orbits of $T$ are contained in $I\setminus W(T).$
% %Furthermore, 
% Since the entropy of $T_{\xi}$ is zero, either
% the right endpoint of each plateau of $T_{\xi}$ is periodic or it
% can be approximated by periodic orbits, and since 
% ${\mathcal Per}(S)=\{2^n:n\in\mathbb N\cup\{0\}\},$ at least one right endpoint
% of a plateau is not periodic.
% Hence, there is at least one plateau whose
% right end point is in $ \Delta(T_{\xi}),$
% where 
% $\Delta(T_{\xi})$ is as in Lemma~\ref{lem:infinitely renorm}. 
% Using this information we can use the proof of
% Theorem~\ref{thm:boundary ssm} 
% to ensure two extra properties: 
% First, that $T_{\xi'}$ and $T_{\xi''}$  have the same shape as $T_{\xi}.$
% Second, that the right end points of all plateaus of $T_{\xi'}$
% and $T_{\xi''}$ are periodic.
% \end{pf}

\subsection{Multimodal mappings of type $\underline b$}
In the next two subsections, we introduce two types of mappings
which arise naturally when one studies renormalization of 
multimodal mappings: multimodal mappings of type $\underline b$
and polynomials of type $\underline b.$ These sorts of mappings
were considered by Smania in \cite{Sm2} under the additional 
assumption that all critical points of the mappings
have order two, which
simplifies the description of the spaces a little.

\begin{defn}\label{def:multimodal map}
Given a vector $ \underline b=(\ell_1, \dots ,\ell_b)$
of positive even integers, we say that $f$ is a
\emph{multimodal map of type \underline b}
if it can be written as a decomposition of
$b$ maps $f_i\in \mathcal A(I)$  
(or more generally 
in $\mathcal A^k(I)$, $k=0,1,2,\dots,\infty,\omega)$ 
as follows:
\begin{itemize}
\item $f= f_b\circ \dots \circ f_1;$
\item $f_i$ has a unique critical point $c_i,$
which is a maximum
and has order $\ell_i;$ 
\item $f_i(\partial I)\subset \partial I;$
\item for $1\leq i\leq b-1$, $f(c_i)\geq c_{i+1}$ and
$f(c_b)\geq c_1$.
\end{itemize}

The vector $(f_1,\dots, f_b)$ gives a {\em decomposition} of $f.$
%To each decomposition of $f$ we associated an {\em extended map}
%$F: I\times\{1,2,\dots b\}\to I\times\{1,2,\dots b\},$
%with $I^n=\{(x,i): x\in I \text{ and }1\leq i \leq
%n\},$
%defined as 
%$$F(x,i)= (f_i(x), (i\ \mathrm{mod}\ b)+1).$$
%We define the \emph{critical set of $F$} as 
%$\crit(F)=\{(c_i,i):1\leq i \leq b\},$ where $c_i$
%is the critical point of $f_i.$
\end{defn}

% All critical points even:
% \begin{defn}\label{def:multimodal map}
% Given a vector $ \underline n=(\ell_1, \dots ,\ell_n)$ of positive integers, we say that $f$ is a \emph{multimodal map of type \underline n} if it can be written as a decomposition of $n$ maps $f_i\in \mathcal A(I)$  (or in $\mathcal A^k(I)$) as follows:
% \begin{itemize}
% \item $f= f_n\circ \dots \circ f_1;$
% \item $f_i$ has a unique turning point $c_i,$ which is a maximum and has order $\ell_i.$ 
% \item $f_i(\partial I)\subset \partial I;$
% \item $f(c_i)\geq c_{i+1}$ \emph{mod} n.
% \end{itemize}

% The vector $(f_1,\dots, f_n)$ gives a {\em decomposition} of $f.$ Each decomposition of $f$ has associated an extended map $F: I^n\to I^n,$ with $I^n=\{(x,i): x\in I \text{ and }1\leq i \leq n\},$ defined as 
% $$F(x,i)= (f_i(x), i+1 \textit{ mod n}).$$
% We define the \emph{critical set of $F$} as  $\crit(F)=\{(c_i,i):1\leq i \leq n\},$ where $c_i$ is the critical point of $f_i.$ 
% \end{defn}
Multimodal mappings of type $\underline b$ arise naturally
as renormalizations of multimodal mappings in $\mathcal A^k(I)$.
% They were introduced in \cite{Sm2} under the additional 
% assumption that all critical points have order two, which
% simplifies the description of the space a little.

%As for multimodal maps, we define renormalizable maps on this space using restrictive intervals. More formally, let $f$ be a multimodal map of type $\underline b$ and consider an extended map $F$ induced by a decomposition $(f_1,\ldots,f_b).$ We say that $J$ is a \emph{$k$-periodic interval} of $F$ if:
%\begin{itemize}
%\item $c_1\in J,$ where $c_1$ is the critical point of $f_1.$
%\item The interval $J, F(J) \ldots F^{k-1}(J)$ are pairwise disjoint.
%\item $F^k(J)\subset J$
%\item The orbit of $J$ contains $\crit(F).$
%\end{itemize}

%Under these circumstances we call $k$ the period of $J.$
%Furthermore, if $J$ is the maximal periodic interval of period $k,$
%then $F^k(\partial J)\subset \partial J$
%and we say that $J$ is a \emph{restrictive interval}
%of $F$ \emph{of period} $k$. 
%If a map $F$ possesses a restrictive interval
%we will say that it is \emph{renormalizable}.
It is easy to see that the renormalization of a 
multimodal mapping of type $\underline b$ is a 
multimodal mapping of type $\underline b'$,
where $\underline b'$ depends on $\underline b$ and 
the combinatorics of $f$.
%We will only renormalize about the critical point 
%$c_1$. 
%For us this is not a restriction, since the critical points
%of the multimodal mappings of type $\underline b$ we consider all have the
%same $\omega$-limit set.
%If $f$ possesses an infinite sequence of
%restrictive intervals, we say that the mapping $f$ is
%\emph{infinitely renormalizable}. 

\subsection{Polynomials of type $\underline b$}
\begin{defn} \label{P_n}
Given a vector  $ \underline b=(\ell_1, \dots ,\ell_b)$  of
positive even 
integers we define the space  $\mathcal P_{\underline b}$ of \emph{polynomials if type $\underline b$} as follows. A polynomial $p\colon I\to I$ belongs to $\mathcal P_{\underline b}$ if
$$p= q_b\circ\dots\circ q_1,$$   where $q_i\colon I\to I$  has the following properties for $i\in \{1,\ldots, b\}$: 
$q_i(-1)=q_i(1)=-1,$ $q_i(0)>0$ for $i\neq b$,
 and $q_i=A_i^{-1}\circ p_i\circ A_i,$ where  
$p_i\colon \R \to \R$ is a polynomial of the form 
$z^{\ell_i}+a_i$ which has an invariant interval $J_i$ and 
$A_i\colon  I \to J_i$ is an affine bijection. 
The vector given by $(q_1,\ldots, q_b)$ will be called a 
\emph{decomposition of $p.$}
We identify affinely conjugate polynomials.

\end{defn}
It is important to note that the number of turning points of maps in $\mathcal
P$ is not constant, but it is uniformly bounded by a constant depending only on
$b$. Also,
observe that, if in addition we have that $q_b(0)>0,$ then $p\in \mathcal
P_{\underline b}$ is a multimodal map of type $\underline b.$

% The families $\mathcal P_{\underline b}$ are the natural generalization of the polynomials of type $b$ used in \cite{Sm2}; Smania works under the assumption that all entries of the vector $\underline b$ are equal to two, while we consider all even integers. 

% We will denote by $N(n)\in \N$ the maximal number of turning points of maps in $\mathcal P_{\underline n}.$  

%Given $m\in \N$ we will denote by $\mathcal P^{m}$\todo{} the subset
%of maps in $\mathcal P$ with $m$-turning points. 

Given a polynomial $p\in \mathcal P $ with shape $\tau$ we let$$\mathcal P(\tau)=\{ q\in \mathcal P: \text{the shape of \emph{q} is equal to } \tau\}.$$  
%The shape will play a fundamental role in the proof of Theorem \ref{main thm}, but before we get to that we will show some important properties of the space $\mathcal P_{\underline n}.$
We say that a shape $\tau$ is \emph{admissible for polynomials of
  type} $\underline b$ if $\mathcal P(\tau)\neq\emptyset.$

Let us fix
$\underline b$ and denote the family $\mathcal P_{\underline b}$
simply as $\mathcal P.$

\begin{lem}\label{lem:uniqueness}
Each map $p\in \mathcal P$ has a unique decomposition.
\end{lem}
\begin{pf}
Let $p= q_b\circ\dots\circ q_1$ and let $q_i, p_i, A_i$ and $J_i$ be as in the definition of $\mathcal P.$ Since $p_i=z^{\ell_i}+a_i$ the interval $J_i$ is symmetric with respect to the origin and its right end point is the point $b_i>0$ so that $p_i(b_i)=b_i.$  Since  the  $\ell_i$ is fixed, the value of $b_i$ depends only on $a_i.$  There exist only two affine maps which map $[-1,1]$ to $[-b_i,b_i]$ bijectively, which are $z\to b_iz$ and $z\to -b_iz.$ Since $q_i(-1)=q_i(1)=1$ we get that $A_i(z) =-b_i z.$  Hence, $q_i$ depends only on $a_i.$ 
The uniqueness of the decomposition of $p$ can be proved by induction on $b,$  the length of $\underline b.$ If $b=1$ then $q_i=q_i'$ if $$\frac{1}{b_i}[(b_i z)^d+a_i]=\frac{1}{b'_i}[(b'_i z)^d +a'_i]$$
Hence $a_i=a_i'.$ By definition of $b_i$ and $b_i'$
this implies that $b_i=b_i'.$ So the decomposition is unique.
Assume that the decomposition is unique when the length of the decomposition is $b.$
To prove the result when the length of the decomposition is $b+1$ we observe the following. 
We have $p=q_{b+1}\circ \dots \circ q_1.$ 
The polynomial $q_{b+1}$
has exactly one turning point, and the critical value of
$p$ corresponding to the critical point of $q_{b+1}$
is determined by the
critical value of 
$q_{b+1},$ hence it depends only on $a_{b+1}.$
Since the map   $p'=q_b\circ\dots q_1$ has a unique decomposition the result follows. 
\end{pf}

\begin{lem}\label{lem:critical val}
Suppose that $\tau$ is  
admissible for polynomials of type $\underline b,$ and 
that $g\colon I\to I$ is a piecewise monotone map with $\tau(g)=\tau$,
then there exists a unique map $q\in \mathcal P(\tau)$
with the same critical values as $g$.
\end{lem}

% \begin{lem}\label{lem:critical val}
% Suppose that $\tau$ is a shape and that $X$ is 
% set of marked itineraries. Assume that the 
% data $(\tau, X)$ is
% admissible for polynomials of type $\underline b,$ and 
% that $g:I\to I$ a piecewise monotone map with $\tau(g)=\tau$,
% then there exists a unique map $q\in \mathcal P(\tau)$
% with the same critical values at turning points
% as $g$ and with odd critical points
% with itineraries in $X$.
% \end{lem}

\begin{pf}
This result can be shown by induction on the length of $\underline b,$ in a similar fashion as the proof of Lemma~\ref{lem:uniqueness}.
\end{pf}

We say that two mappings are \emph{essentially conjugate}
if they are topologically conjugate outside of their
basins of attraction.

\begin{prop}\label{prop:essentially conjugate}
Given a shape $\tau$ that is admissible for polynomials of type
$\underline b$
and a piecewise monotone map $g\colon I\to I$ with $\tau(g)=\tau$,
there exists a map $q\in \mathcal P(\tau),$
which is essentially conjugate to $g.$
%and has the same critical values.
\end{prop}

% \begin{prop}\label{prop:essentially conjugate}
% Given a shape $\tau$ and marked kneading data $X$
% so that $(\tau, X)$ is admissible for polynomials of type
% $\underline n$
% and a piecewise monotone map $g:I\to I$ with $\tau(g)=\tau$,
% there exists a map $q\in \mathcal P(\tau),$
% which is essentially conjugate to $g$ and has the same critical
% values at turning points as $g$,
% and with odd critical points with itineraries in $X$.
% \end{prop}

\begin{pf}
This result follows from the previous lemma and  the proof of Step 1 in Theorem II.4.1 in \cite{dMvS}. 
\end{pf}

The following two results of \cite{Sm2}
generalize immediately 
to polynomials with critical points not a power of two, 
so we have not included their proofs. Before we can state them, we must introduce some notation. Let  $\mathcal P^*$  denote the set of maps of the form $p=p_b\circ \dots \circ p_1,$ where $p_i=z^{\ell_i}+a_ i$ and $(\ell_1,\ldots,\ell_b)=\underline b$ and let $Poly(\hat b)$ denote the set of monic polynomials of degree $\hat b=\ell_1\ldots \ell_b.$

\begin{prop}\cite[Proposition 3.1]{Sm2}
The space $\mathcal P^*$ 
is a complex submanifold of $Poly(\hat b)$ with parametrization
$$(a_1,\ldots, a_b) \to P_{a_b}\circ\cdots\circ P_{a_1}.$$
\end{prop}

The \emph{connectedness locus} of $Poly(\hat b)$ is
the set of all mappings in $Poly(\hat b)$ with connected
Julia set. 
\begin{prop}\cite[Proposition 3.2]{Sm2}
The connectedness locus
of $Poly(\hat b)$ is compact.
\end{prop}

\subsubsection{Stunted sawtooth mappings and polynomials}
In this section we present the results from \cite{BruvS} which we will 
use in later sections. 
%It is important to note that  \cite{BruvS} works under the assumption that all turning points of polynomials are of order $2$. 
For a fixed vector $\underline m= (\ell_1,\ldots, \ell_m)$ of even integers  let  $\mathcal Q$ denote the space of polynomials $q\colon I\to I$ with $q(-1)=q(1)=-1,$ with $m$ turning points $-1<c_1<\ldots c_m<1$ where the order of $c_i$ is $\ell_i.$ 
%vectors $\underline n$ of positive even integers and the proofs remain unchanged. So without loss of generality we will state results from \cite{BruvS} in this general setting.

\medskip

\medskip
Given $p\in \mathcal Q$ the following holds. Let $\mathcal S=\mathcal S_m$ and $S_0$ be defined as in Section \ref{subsec:ssm}. Let $\nu(p)=(\nu_1,\ldots, \nu_m)$ be the kneading invariant of $p$ and let $s_i$ be the unique point in the $(i+1)$-th lap of $S_0$ such that
$$ \lim_{y\downarrow s_i} i_{S_0}(y)=\nu_i:=  \lim_{x\downarrow c_i} i_f(x).$$
Let $Z_i$ be the symmetric interval around the $i$-th turning point of
$S_0$ with right end point $s_i.$ Then we can define a map 
\label{def:psi}
$$\Psi\colon \mathcal Q\to \mathcal S\quad \text{ by }\quad p\to \Psi(p),$$
where $\Psi(p)$ is the unique map in $\mathcal S$ which agrees with $S_0$ outside $\cup^m_{i=1} Z_{i}$ and which is constant on $Z_i$ with value $S_0(s_i).$

The following result summarizes some of the key properties of $\Psi.$
\begin{lem}\cite[Lemma 5.1]{BruvS}\label{Psi}
The map $\Psi\colon  \mathcal Q\to \mathcal S$
\begin{itemize}
\item is well-defined;
\item the kneading invariant of $p$ and of $T=\Psi(p)$
are the same in
  the sense that $\lim_{y\downarrow Z_i} i_{T}(y)=\nu_i;$
\item $p$ and $\Psi(p)$ have the same topological entropy;
\item $\Psi(p)$ is non-degenerate (see the next paragraph). 
\end{itemize}
\end{lem}

Recall that given a map $T\in \mathcal S$ a pair of plateaus $(Z_i,Z_j)$ is called \emph{wandering} if there exists $n\in \N$ such that $T^n$ of the set $[Z_i,Z_j]$ (the convex hull of $Z_i$ and $Z_j$) is a point. We say that a map $T\in \mathcal S$ is \emph{non-degenerate} if for every wandering pair $(Z_i,Z_j)$ its convex hull belongs to the closure of a component of the basin of a periodic plateau. We will denote by $\mathcal S^*$\label{def:S*} the set of non-degenerate maps in $\mathcal S.$ In particular, Lemma  4.16 in \cite{BruvS} tells us the following.

\begin{lem}\label{admissible ssm}
If we take a map $T\in \mathcal S^*$ with no periodic attractors, then there exists $\nu_0>0$ so that $B_{\nu_0}(T)\subset  \mathcal S^*.$
\end{lem}

\subsection{Polynomial-like mappings and germs}
\subsubsection{Polynomial-like mappings of type $\underline b$}

\begin{defn}
Given a vector $\underline b=(\ell_1,\ell_2,\ldots, \ell_b)$ of
positive even integers we say that a polynomial-like mapping $f\colon U\to V$ is a
\emph{polynomial-like map of type $\underline b$} if there exist
simply connected domains $U=U_1,\ldots, U_b,U_{b+1}=V$
and holomorphic maps $f_i\colon U_i\to U_{i+1}$ with $i\in \{1,\ldots, b\}$ satisfying:
\begin{itemize}
\item For $1\leq i\leq b$, $f_i\colon U_i\to U_{i+1}$ is a branched covering 
of degree $\ell_i$ with exactly one ramification point.
\item $f=f_b\circ\cdots\circ f_1.$
\end{itemize}
\end{defn}
We denote the space of polynomial-like mappings of type 
$\underline b$ by $\mathcal{PL}_{\underline b}$. 
For future reference, we define the {\em type} of an AHPL-mapping in the same way.

% \medskip
% \noindent\textit{Remark.}
% \emph{Polynomial-like maps of type $b$} were
% introduced in \cite{Sm2}. The difference is that here
% we do not require 
% each $F_i$ to be quadratic.

\medskip

The following result is an analogue of the Douady-Hubbard 
Straightening Theorem for polynomial-like mappings of type $\underline b$.
\begin{prop}\cite[Proposition 4.1]{Sm2}
Let $\underline b=(\ell_1,\ell_2,\ldots, \ell_b)$ be a vector of
non-negative even integers. Assume $f\colon U\to V$ is a polynomial-like map
of type $\underline b$ and that the critical values of $f$ are
contained in $U.$ Then $f$ is 
hybrid conjugate
to a polynomial $P=\chi(f)$
in $\mathcal P^*.$
\end{prop}
We call the polynomial $P$ the {\em straightening of} $f$, and we 
refer to the mapping $\chi$ as the {\em straightening map}.
See page~\pageref{def:hybrid eq} for the definition of hybrid conjugate.

Following \cite{McM}, we endow the space of 
polynomial-like mappings with the {\em Carath\'edory
topology.} 
A {\em pointed disk} is a topological disk $U\subset \mathbb C$ with
a marked point $u\in U$.
Let $\mathcal D$ denote the set of pointed
disks $(U,u).$
We first define the {\em Carath\'eodory topology} on
$\mathcal D$. We say that $(U_n,u_n)\rightarrow (U,u)$ in
$\mathcal D$ if 
\begin{itemize}
\item $u_n\rightarrow u;$
\item for any compact set $K\subset U,$ $K\subset U_n$ for
all $n$ sufficiently large; and
\item for any connected $N\owns u,$ if $N\subset U_n$
for infinitely many $n$, then $N\subset U$.
\end{itemize}
Now, we define the Carath\'eodory topology on the
space of all holomorphic mappings $f\colon (U,u)\rightarrow\mathbb C$,
where $(U,u)$ is a pointed disk. We say that $f_n\colon (U_n,u_n)\rightarrow
\mathbb C$ {\em converges to} $f\colon (U,u)\rightarrow\mathbb C$ if:
\begin{itemize}
\item $(U_n,u_n)\rightarrow (U,u)$ in $\mathcal D$, and
\item for all $n$ sufficiently large $f_n$ converges to $f$ uniformly
  on compact subsets of $U$.
\end{itemize}
We endow the space of polynomial-like mappings 
$f\colon U\to V$ with the Carath\'eodory topology by choosing the marked
point in the filled Julia set.

\subsubsection{Polynomial-like germs of type $\underline b$}

% Two analytic maps $f$ and $g$ defined in a neighbourhood of 
% 0 \emph{represent the same germ at 0} if they coincide in a 
% neighbourhood of 0.

We have the following equivalence
relation on the space $\mathcal{PL}_{\underline b}$: 
Suppose that $f:U\rightarrow V$ and 
$\tilde f\colon \tilde U\rightarrow \tilde V$ are polynomial-like mappings of type
$\underline b$.
We say that
$f\sim \tilde f$ if
$f$ and $\tilde f$ have a common 
polynomial-like restriction of the same degree.
By \cite{McM}, Theorem 5.11, we have that if $f\sim \tilde f,$
then $K_f=K_{\tilde f},$ and for mappings with 
connected Julia set, this is an equivalence relation.
Classes of this equivalence
relation are called \emph{polynomial-like germs}
and we denote the equivalence class of a polynomial-like
mapping $f$ by $[f]$.
Let $\mathcal{PG}$ represent the space of polynomial-like germs,
up to affine conjugacy, and let 
 $\mathcal{PG}^{\R}$ be the subset of real polynomial-like germs.
The space of polynomials is naturally embedded in 
the space of polynomial-like germs.
We let $\mathcal{C}$ denote the connectedness locus 
in $\mathcal{PG}$, and let $\mathcal{C}^{\R}
=\mathcal{C}\cap\mathcal{PG}^{\R}$.

We say that a polynomial-like germ $f\colon U\to\mathbb C$
is {\em renormalizable} at a point $c\in\mathrm{Crit}(f),$
if there exists a neighbourhood $U_1\subset U$ of $c$
and an $s\in\mathbb N\setminus\{1\}$ so that
$f^s\colon U_1\to V:=f^s(U_1)$ is a polynomial-like mapping with connected
Julia set.

The definitions of quasiconformal equivalence and hybrid equivalence
for polynomial-like germs are the same as for polynomial-like mappings.
We denote the hybrid class of a polynomial-like mapping or germ $f$
by $\mathcal{H}_f$.
Any two polynomial-like germs
$[f]$ and $[g]$ in the same hybrid class
$\mathcal H$ can be included in a
\emph{Beltrami disk}: Let $h$
be a hybrid conjugacy between representatives
$f$ and $g$
and let $\mu=\bar\partial h/\partial h$ be its
Beltrami differential. Let $\varepsilon>0$ be so small that
$(1+\varepsilon)\|\mu\|_{\infty}<1$.
Define $\mu_{\lambda}=\lambda\mu,\lambda\in \mathbb D_{1+\varepsilon}.$
By the Measurable Riemann Mapping Theorem,
we obtain a family $h_{\lambda},$
$\lambda\in\mathbb D_{1+\varepsilon},$
of quasiconformal mappings, the solutions of the 
associated Beltrami equations.
A Beltrami disk through 
$f$ and $g$ is a family of mappings
$\{f_{\lambda}=h_{\lambda}\circ f\circ h_{\lambda}^{-1}\colon
\lambda\in\mathbb D_{1+\varepsilon}\}$. If $U$ is a neighbourhood of
$K_f$ to which $f$ has a polynomial-like restriction, then,
since $\mu$ is an invariant Beltrami differential, so is $\lambda\mu,$
and so
each $f_\lambda$ is holomorphic on $h_\lambda(U)$. 
We call a real
one parameter family 
$\{f_{\lambda}: -1-\varepsilon<
\lambda< 1+\varepsilon\}$
a \emph{Beltrami path} through
$f$ and $g$.

To define the topology on $\mathcal{PG}$, we push down the 
Carath\'eodory topology which we defined on the space of
polynomial-like mappings, see \cite{L-Feig}.
We say that a sequence of polynomial-like germs
$[f_n]\rightarrow [f]$
if the sequence of $[f_n]$ and be split into finitely
many subsequences $[f_m^i]$ which admit representatives
$f_m^i$ which converge to representatives $f^i$ of $f.$
In the case when $J(f)$ is connected, which is the one that is
important to us, we do not need to split the sequence
$[f_n]$ into subsequences.

\subsubsection{External mappings and matings}
Let us fix $d\in\mathbb N,$ $d\geq 2$.
Let $g\colon\mathbb T\to\mathbb T$ be a degree $d$ real-analytic
endomorphism of the unit circle. We say that 
$g$ is \emph{expanding} if it admits an extension to 
a degree $d$ covering map
$g\colon U\rightarrow V$ between annular neighbourhoods
of $\mathbb T$ with $U\Subset V$. We normalize
$g$ by the condition that $g(0)=0$.
Let $\mathcal{E}$ denote this space of normalized
expanding endomorphisms of the circle. We endow $\mathcal E$ with a topology as
follows, see \cite{AL}. Since $\mathbb R$ is the universal covering of $\mathbb T,$ any map
$g\in\mathcal E$ lifts to a mapping $\tilde g\colon \mathbb R\to \mathbb R,$
so that $\tilde g(x)=dx+\phi(x),$ where $\phi$ is a real-analytic function with
period one and $\phi(0)=0.$ Let $\mathcal A$ denote the space of all such
functions, and let $\mathcal A_n\subset \mathcal A$ denote the subset of functions which admit an extension to the strip $|\mathrm{Im}
z|<1/n,$ for $n\in\N,$ which are continuous up to the boundary. Since the
$\mathcal A_n$ are Banach spaces, $\mathcal A$ is realized an inductive limit of
Banach spaces, and we may  endow it with the inductive limit topology. This
topology on $\mathcal A$ yields a topology on $\mathcal E,$ and in this topology
a sequence $g_n\in\mathcal E$ converges to $g\in\mathcal E$ if there is a
neighbourhood $W$ of $\partial T$ such that all the $g_n$ admit a holomorphic
extension to $W,$ and $g_n$ converge uniformly to $g$ on $W.$
Let $\mathcal E^{\mathbb R}$ denote the space of 
real-symmetric expanding maps of the circle.

For $g\in\mathcal E$, let $\mod(g)=
\sup\mod(V\setminus(U\cup\mathbb D))$,
where the supremum is taken over all
extensions $g\colon U\rightarrow V$ as above.

To each $f\in\mathcal{PG}$
of degree $d$,
we associate its external mapping
$\pi(f)=g$: we let $f$ be a polynomial-like representative
of the germ  $f$, and then use the construction of
\cite{L-Feig}, Section 3.2 to obtain $g\colon\mathbb T\to \mathbb T$. 
See also, \cite{DH}, Section 2. 
We say that two polynomial-like germs, $f$ and $g$, are 
externally equivalent if $\pi(f)=\pi(g)$.

\begin{lem}
We have the following:
\begin{itemize}
\item If $f$ and $g$ are externally equivalent polynomial-like germs
with connected Julia sets,
then there is a conformal mapping
$h\colon \C\setminus K_f\rightarrow\C\setminus K_g$
which conjugates $f$ and $g$ near their Julia sets.

\item The external mapping $\pi(F)(z)=z^d$ if and only if
$F$ is a polynomial of degree $d$.
\end{itemize}
\end{lem}

\begin{thm}[Mating Theorem]\label{thm:DHstraightening}
If $P$ is a real polynomial of degree $d$ with connected Julia set and
$g\in\mathcal{E}$, then there exists, up to affine conjugacy,
a unique germ $f=i_P(g)\in\mathcal{PG}$ such that
$\chi(f)=P$ and $\pi(f)=g$. 
\end{thm}

The following theorem can be obtained in exactly the same
way as for unicritical mappings,
see \cite{AL}.
Let $\mathcal M\subset\mathcal P_{\underline b}$ denote
the subset of $\mathcal P_{\underline b}$
of mappings $f$ such that the Julia set, $J(f),$ is connected.
We let $\mathcal M^{\mathbb R}$ denote the 
real slice of $\mathcal M.$
To simplify matters, we restrict to the
real slices of these complex spaces.
\begin{thm}[\it{c.f.} \cite{AL}, Theorem 2.2, \cite{Sm2}, Proposition 4.1] 
\label{thm:lamination}
There is a canonical choice of the straightening
$\chi(f)\in\mathcal M^{\mathbb R},$ and an external mapping
$\pi(f)\in\mathcal E^{\mathbb R}$ associated to each germ
$f\in\mathcal C^{\mathbb R}$. %which depends continuously on $f$.
It has the following properties:
\begin{enumerate}
\item For each $P\in\mathcal M^{\mathbb R}$, the hybrid leaf
$\mathcal H_p^{\mathbb R}$ is the fiber $\chi^{-1}(P)
\cap\mathcal C^{\mathbb R}$ 
and the external
map $\pi$ restricts to a homeomorphism 
$\mathcal H_P^{\mathbb R}\to\mathcal E^{\mathbb R},$
whose inverse is denoted by $i_P$ and is called the
(canonical mating).
%\item $(\pi,\chi):\mathcal C^{\mathbb R}\to
%\mathcal E^{\mathbb R}\times\mathcal M^{\mathbb R}$
%is a homeomorphism.
\item For $P,\tilde P\in\mathcal M^{\mathbb R},$ if $f_{\lambda}$ is a 
Beltrami path in $\mathcal{H}_P$, then $i_{\tilde P}\circ
i_{P}^{-1}(f_{\lambda})$ is a Beltrami path in $\mathcal H_{\tilde P}^{\mathbb R}$.
\item the external map, straightening and mating are 
equivariant with respect to complex conjugation.
\end{enumerate}
\end{thm}

\begin{prop}
  \label{prop:contstr}
Suppose that $f\colon U\to V$ is a real polynomial-like mapping of type $\underline b$
with a single solenoidal attractor which contains $\mathrm{Crit}(f).$
Then the straightening map is continuous at $f.$
\end{prop}
\begin{pf}

We will use the following lemma:

\begin{lem}[{\it c.f.} \cite{DH}, Lemma page 313]\label{lem:DH-313}
Suppose that $\{f_{\lambda}\colon U_{\lambda}\rightarrow
V_{\lambda}\}_{\lambda\in\Lambda}$
is an analytic family of polynomial-like mappings,
 where $\Lambda$ is a complex analytic manifold.
Let $\lambda_0\in\Lambda$, and suppose that
$\{\lambda_n\}_{n=1}^\infty\subset\Lambda$ is a sequence so that
$\lambda_n\rightarrow\lambda_0$ as $n\rightarrow\infty$.
Then there exists a sequence $\{\lambda_{n_k}\}_{k=1}^\infty$ of
 $\{\lambda_n\}_{n=1}^\infty$ such that the sequence of polynomials
$P_{n_k}=\chi(f_{\lambda_{n_k}})$ converges to a polynomial $\tilde
P$, and $\tilde P$ is quasiconformally equivalent to $f_{\lambda_0}$.
\end{lem}

Thus for any sequence $f_{n}\colon U_n\to V_n$ of real polynomial-like mappings
converging to $f,$ we have that there exists a subsequence $f_{n_k}$ so that
$\chi(f_{n_k})$ converges to a polynomial $\tilde P,$ which is quasiconformally
conjugate to $f.$ If $\tilde P$ is hybrid conjugate to $f,$ then since the Julia
set of $f$ is connected, we have that the straightening of $f$ is unique, and so
$\chi(f)=\tilde P,$ and the straightening is continuous. So we may assume that
$\tilde P=h^{-1}\circ f\circ h $ is not hybrid conjugate to $f,$ but then the
Beltrami differential $\mu=\bar \partial h/\partial h$ gives an invariant line
field supported on the filled Julia set of $f$ (see page~\pageref{page:bdlf}), which contradicts Theorem~\ref{thm:nilf}.
 \end{pf}

\begin{prop}
\label{prop:windows shrink}
Let $b\in\mathbb N$, and let 
$\underline b$ be a $b$-tuple of even integers. Assume that
$\mathcal X$ is a compact subset of 
$\mathcal{PG}_{\underline b}.$ Suppose that
$\mathcal V_n \subset \mathcal{PG}_{\underline b},$
is the set of mappings that are at least $n$-times
renormalizable (of period 2). Then if $f_n\in\mathcal V_n\cap \mathcal X$ 
is any sequence,
$f_n\rightarrow\Gamma_{\mathcal{PG}_{\underline b}}$
in the Carath\'eodory topology.
\end{prop}
\begin{pf}
Since the $f_n$ are contained in compact set, 
any subsequence of the $f_n$ must have a 
subsequence which converges. By 
Theorem~\ref{thm:qs rigidity},
this limit is infinitely renormalizable,
so it is in $\Gamma$.
\end{pf}

\subsubsection{Infinitesimal structure of the space of polynomial-like germs}
The description of the tangent space of the space of polynomial-like germs was
first given in \cite[Section 4]{L-Feig} in the context of unicritical mappings, and
\cite[Sections 3 and 4]{ Smania-shy} treated polynomial-like germs with  several critical points. We refer to those papers for the details.

% Suppose that $f\in\mathcal{PG}.$ Let $T_f\mathcal{PG}$ denote the tangent space
% to $f$. It can be identified with the inductive limit, $\lim_U \mathcal B_U,$ over
% the directed set of domains $U$ on which $f$ has a quadratic-like restriction.

\begin{prop}
  \label{prop:transdim} Suppose that $f\in\mathcal{PG}_{\underline b}.$
Let $\mathcal H_f$ denote the hybrid class of $f.$ Then $\mathcal H_f$ is a
connected, codimension-$b$ complex-analytic submanifold of $\mathcal{PG}.$
\end{prop}

Since the renormalization operator is transversally non-singular, we can transfer some of this structure to
infinitely renormalizable analytic mappings:

\begin{thm}[\cite{ALdM}, Lemma 4.8, and\cite{Smania-shy}, Theorem 3]
  \label{thm:R is trans nonsing}
Suppose that $f\in\mathcal A_{\underline b}(I)$ is an analytic map and that $c$
is a critical point of $f,$ at which $f$ is infinitely renormalizable. Suppose
that $\mathcal R(f)=F\colon U\to V$ is a polynomial-like renormalization of $f$ at $c,$ and that
$v\in T_F\mathcal{PL}$ is transverse to the topological conjugacy class of $f$.
Then there exist vectors $w_i\in T_f\mathcal A_{\underline b}(I),$ so that
$D\mathcal R(f)w_i\to v.$ 
\end{thm}

\subsubsection{Convergence of renormalization for 
analytic mappings}
McMullen proved exponential convergence of renormalization of
quadratic-like mappings with bounded combinatorics, \cite{McM2}. 
These results were generalized by Smania to multimodal mappings
with all critical points of degree 2 \cite{Sm2}, \cite{Sm3}. 
By Theorem~\ref{thm:complex bounds}
and the quasiconformal rigidity of analytic mappings,
Theorem~\ref{thm:qs rigidity}, we have exponential convergence of
renormalization for infinitely renormalizable analytic mappings
with bounded combinatorics.

\begin{thm}
\label{thm:ren limit set}
For any $\underline b$,
there exists $\lambda\in(0,1), \delta>0$ so that the following holds.
Suppose that $f,g\in\mathcal{A}_{\underline b}(I)$ are topologically conjugate
mappings that are infinitely renormalizable at corresponding critical points
$c_0$ and $\tilde c_0,$ respectively.
Let $\mathcal R^n_{c_0}$ denote the $n$-th renormalization at $c_0$.
There exists $C>0$, depending also on the
combinatorics of $f$, so that 
$$\|\mathcal R^n_{c_0}(f)-\mathcal{R}^n_{\tilde c_0}(g)\|_{\mathcal C^0(U_\delta)}\leq
C\lambda^n,$$
where $U_\delta$ is a $\delta$-neighbourhood of the filled Julia set of
$\mathcal R^n_{c_0}(f).$
Moreover the limit set of $\mathcal R^n(f)$ is
contained in a Cantor set $\mathcal K$.
\end{thm}

\noindent\textit{Remark.}
The attractor of the period doubling renormalization
operator for multimodal mappings, with more than one
critical point, is a horseshoe, \cite{Sm2}, \cite{OT}.
This is in contrast with the unimodal case where the 
period-doubling renormalization operator has
stationary combinatorics and the attractor is a fixed point.

\section{The boundary of chaos}
\label{sec:boundary}
\subsection{Boundary of chaos for polynomials of type $\underline b$}
Let $p\in\mathcal P$ have decomposition $(q_1\ldots,q_b).$

% with $q_b(0)>0$, let 
%$$\hat p :I\times\{1,\dots,b\}\to I\times\{1,\dots,b\}$$ 
%denote the extension of $p$
%given in 
%Definition~\ref{def:multimodal map}.

%Observe that for a map $p\in\Gamma_{\mathcal P},$
%all turning points have the same
%$\omega-$limit set.

%\begin{thm}\label{boundary unif}
%Given $P\in \mathbb B$  with shape $\tau$ it is on the boundary of chaos of $\mathcal P_{\underline n}^{\tau}.$
%\end{thm} 

%To prove this result we will make use of the equivalent result for the space of stunted sawtooth maps, Theorem \ref{boundary ssm}. This will be done by relating stunted sawtooth maps with interval maps via kneading theory. For the definition and properties of kneading invariants we refer the reader to \cite{dMvS}.

\begin{thm}\label{thm:boundary poly}
Let $p\in \Gamma=\Gamma_{\mathcal P}\subset\mathcal P$ and let $(q_1\ldots,q_b)$ and  $\tau$
denote the decomposition and shape, respectively, of $p.$
Given
$\nu>0$ there exist $q,r\in \mathcal P(\tau)\cap B_\nu(p)$
both with all turning points periodic or pre-periodic such that
\begin{itemize}
\item $\mathcal Per(q)$ is finite, and
\item $r$ has positive entropy.
\end{itemize}
\end{thm}

\begin{pf}
% Recall that $\mathcal P(\tau),$ $\mathcal S(\tau)$ and 
% $\mathcal{PM}(\tau)$ respectively 
% denote the spaces of polynomials of type
% $\underline b$ 
% with shape $\tau$,
% stunted sawtooth maps with shape $\tau$, 
% and piecewise monotone maps with shape $\tau$.
 Since
$p\in\Gamma_{\mathcal P}$, $p$ has at least one solenoidal attractor and zero
entropy.
Apply Lemma~\ref{Psi} to obtain  a map $T=\Psi(p)\in \mathcal S$  with
the same kneading invariant as $p$ and with the same entropy, $h(T)=0$. 
Each critical point of $p,$ corresponds to a plateau of $T.$
Thus the plateaus of $T$ which correspond to critical points of $p$
that are not contained in basins of periodic attractors are neither
periodic nor pre-periodic and there exist such plateaus. Let us label these
plateaus by $Z_1,\dots, Z_m.$ 
All points with periodic itineraries for $T$ are contained in $I\setminus (\mathrm{int}(\cup_{i=1}^mZ_{i,T})).$ Hence, each periodic itinerary corresponds to a unique periodic orbit of $T$ and $\mathcal Per(T)=\{2^n:n\in \N \cup\{0\}\}.$

By Proposition~\ref{prop:bcssm},  we can find maps $\{T_{k}\}_{k\in \N}$
with all plateaus periodic or pre-periodic and with finitely many
periodic orbits. Furthermore, we can guarantee that  $T_{k}\in
B_{1/k}(T)$ for each $k\in \N.$ In addition, by 
Lemma~\ref{admissible ssm}  we may assume that $T_{k}\in \mathcal
S^*,$ see page~\pageref{def:S*}.  Now we follow a procedure used in Proposition 5.9 in
\cite{BruvS}.
For each 
$k\in\mathbb N$ define $x \sim_k y$ if there exists $i>0$ so that $T_k^i(x)$ maps the convex hull $[x,y]$ into one of the plateaus of $T_k.$ Then collapse each of these intervals $[x,y]$ to a point and let $\hat T_k$ be the corresponding map. By definition we have that  $\hat T_k$ is continuous and since $T\in \mathcal S^*$ has 
$m$-turning points so does  $\hat T_k.$  
Furthermore, by construction each $\hat T_k\in \mathcal S^*(\tau),$ has no wandering intervals, no inessential attractors and its kneading invariant corresponds to the one of $T_k$. 
%We apply the fullness of families result, the Corollary of Theorem 4.1 \cite{dMvS}, to find a polynomial $p_k$ with $m$-turning points. Furthermore, the proof of such result guaranteed that $p_k$ has the same critical values as $\hat T_k,$  $p_k$ has shape $\tau.$ 
By Proposition~\ref{prop:essentially conjugate}
there exists  $p_k\in \mathcal P(\tau),$ 
which is essentially conjugate to $\hat T_k.$ So $p_k$ and
$\hat T_k$ are conjugate. Hence $p_k$ has finitely many periodic orbits and entropy zero.
Since the connectedness locus of $\mathcal P^*$ is compact, there
exists a subsequence $\{p_{k_j}\}_{j\in \N}$ which converges to a map
$p'.$ Without loss of generality assume $p_k\to p'.$ 
Corollary~\ref{cor:boundary} and 
Lemma~\ref{lem:infinitely renorm}
imply that  $p'$ is an infinitely renormalizable map with entropy zero.
Finally, since the kneading invariant of the maps $T_k$ converges to
the kneading invariant of $T,$ we have that the kneading invariants of
the maps $p_k$ converge to the kneading invariant of
$p'.$ Hence, $p'$
has the same kneading invariant as $p.$ 
By Theorem~\ref{thm:rigidity}, 
$p$ and $p'$ are conjugate by an affine map.

By an analogous argument as the one used to
construct the maps $p_k$,
we can 
find a sequence of maps $q_k\in \mathcal P(\tau)$,
which have all critical points periodic or pre-periodic and
positive entropy so that $q_k\to q'.$
% By construction we have that $p'$ and $q'$ have the same critical values, so $p'=q'.$  
\end{pf}

\subsection{Boundary of chaos for polynomial-like germs}

\begin{prop}\label{prop:bcpl}
Suppose that $f$ is a real polynomial-like germ which is
infinitely renormalizable at a critical point $c$. Assume that
$h(f|_{K(f)\cap \mathbb R})=0$ and that all critical points of $f$ are
even and have 
the same $\omega$-limit set. Then there exist polynomial-like germs
$g$ and $\tilde g$ arbitrarily close to $f$ in the
Carath\'eodory topology such that
\begin{itemize} 
\item $g$ has finitely many periods, and 
\item $h(\tilde g|_{K(\tilde g)\cap\mathbb R})>0.$
\end{itemize} 
\end{prop}

\medskip
\noindent\textit{Remark.}
By infinitely renormalizable, we mean that 
the restriction of $f$ to its real trace is infinitely
renormalizable about $c$.

\medskip

\begin{pf}
Let $P=\chi(f)$ denote the straightening of $f$. 
By Theorem~\ref{thm:boundary poly},
there exists a sequence of 
polynomials $P_k$ converging to 
$P$ such that each $P_k$ is critically finite and $h(P_k)=0$.
By Theorem~\ref{thm:lamination}, the hybrid classes 
of the $P_k$ are connected submanifolds
in the space of polynomial-like germs and
laminate $\mathcal{PG}.$ Hence for any neighbourhood 
$B\subset\mathcal{PG}$ of $f$, there exists $k$ so that 
$\mathcal{H}_{P_k}\cap B\neq\emptyset$. Hence by Proposition~\ref{prop:contstr}
there exists  
a sequence of critically finite
polynomial-like germs $g_k$ with $h(g_k)=0$ converging to $f$.

Similarly, there exists a sequence of 
polynomials $P_k$ converging to 
$P$ such that each $P_k$ is critically finite and $h(P_k)>0$,
and the same argument implies that there is a sequence of
critically finite
polynomial-like germs $\tilde g_k$ with positive entropy
converging to $f$.
\end{pf}

\subsection{Boundary of chaos for analytic mappings}

Suppose that $\Lambda$ is a $\mathcal C^k$-manifold with
base point $\lambda_0\in\Lambda.$
Suppose that $X\subset \mathbb C,$ we say
that $h_{\lambda}\colon X\to\mathbb C,\lambda\in\Lambda,$
is a
$\mathcal C^k$-\emph{motion} if 
\begin{itemize}
\item $h_{\lambda_0}=\mathrm{id},$
\item $h_{\lambda}$ is an injection for each $\lambda\in\Lambda,$ and
\item $\lambda\mapsto h_{\lambda}(z)$ is $\mathcal C^k$ in $\lambda$.
\end{itemize}
We say that $h_\lambda\colon X\to\mathbb C,\lambda\in\Lambda$ is
a \emph{holomorphic motion} if additionally we
require $\Lambda$ to be a 
complex Banach manifold, and the mapping 
$\lambda\mapsto h_{\lambda}(z)$ to be holomorphic.

Suppose that $f\in\mathcal{B}_{\Omega_a}^{\mathbb R}.$
Let $\mathcal W$ be a neighbourhood of $f$ in
$\mathcal{B}_{\Omega_a}^{\mathbb R}.$ We say that 
a periodic interval $K=K^f\subset [-1,1]$ 
of period $s$
\emph{persists} in $\mathcal W$
if for each $g\in\mathcal W$, there is a $\mathcal C^k$-motion
$h^g\colon K^f\to K^g$, and $K^g$ is a restrictive interval of period $s$
and $h^g\circ f^s(z)=g^s\circ h^g(z)$ for $z\in\partial K$.
We call $K^g$ the \emph{continuation of} $K^f$ to $g.$
Similarly, if  $\mathcal W$ is a neighbourhood of $f$ in
$\mathcal{B}_{\Omega_a}$ and
$f^s|_{U}=F\colon U\to V$ is a polynomial-like mapping,
we say that 
$F\colon U\to V$ persists over $\mathcal W$ if for each 
$g$ in $\mathcal W$ there is a holomorphic motion
$h^{g}\colon (U,V)\to (U_G,V_G),$ 
a polynomial-like mapping, $g^s|_{U_G}=G\colon U_G\to V_G,$
and
$h^g\circ F(z)=G\circ h^g(z)$ for $z\in\partial U.$

\begin{lem}\label{renor window}
Let $f\in \Gamma_{\mathcal A^k(I)}$ and let $c$ be a turning point at which $f$ is infinitely renormalizable. Let $\{J_n\}_{n\in \N}$ be a sequence of restrictive intervals containing $c.$ For $n$ large enough there exists $\epsilon_n>0$ so that $J_n$ persists on $B_{\epsilon_n}(f).$

\end{lem} 
%For $n$ large enough there exists $\epsilon_n>0$ so that $J_n$ persists on $B_{\epsilon_n}(f).$\end{lem}

\begin{pf}
%Let $c$ be a turning point at which $f$ is infinitely renormalizable. Denote by $J_n$ the restrictive interval of period $2^n$ containing $c.$ 
Let $J_n$ be a sequence of restrictive intervals 
containing a turning point $c.$ 
By definition, the boundary points of $J_n$ are:
a periodic point $p_n$ of period $2^n$
and a preimage of $p_n$ under $f^{2^n}.$
By Theorem IV.B in \cite{dMvS} there exists
$M\in \N$ so that all periodic orbits of prime period greater than 
$M$ are repelling. 
Since $p_n$ is hyperbolic for  all $n>M$ 
we have that the interval $J_n$ persists on a 
$\mathcal C^1$-neighbourhood of $f$. 
In other words, there exist a neighbourhood 
$U_n \ni f$ so that the interval $J_n$ 
has a continuation on $U_n$. 
Given a map $g\in U_n,$ we will denote by $J_n^g$ 
its corresponding continuation.
Let us show that $J_n^g$ is a restrictive interval for $g$. 
%Taking $M$ larger if necessary, we can apply Proposition 3.20 in \cite{ClvSTr} to find $\delta>0$ so that the intervals $J_n$ are $\delta$-free for all $n>M;$ this means that the set $\omega(c)$ is not close to the boundary of $J_n$. Hence, for each $n>M$ there exists $\epsilon_n>0,$ such that for any $g\in B_{\epsilon_n}(f)$ the interval $J_n'$ is $\delta/4$-free.  Hence, $J_n'$ is a restrictive interval of period $2^n$ for $g$.
For all $n$ sufficiently large, we can guarantee that the results from
\cite{ClvSTr}  hold for $f$. Since $f\in \Gamma$ we get  that
$I_n=J_n,$ where $I_n$ is an interval form the generalized enhanced
nest. By Theorem 3.1 (a) in  \cite{ClvSTr} there exists
$\delta>0$ so that $V_{n+1}=(1+\delta)J_{n+1}\subset J_n.$
Make $\epsilon_n>0$ small enough so that $\| f^{2^n}-g^{2^n}\|<
\delta/4|J_n|$
and $J_n$ persists on $B_{\epsilon_n}(f)$.
Then, if $g\in B_{\epsilon_n}(f)$ all turning points of
$g^{2^n}|_{J_{n}^g}$
are contained in $V_{n+1}\subset J_n^g.$ 
Hence $J_n^g$ is a restrictive interval for $g$ and the result follows.
\end{pf}

\begin{lem}\label{K_n}
Let  $f\in \Gamma_{\mathcal A^k(I)}$ and let $K_n$ be as in Theorem \ref{omega limit}. For $n\in \N$ large enough, there exists $\nu_n>0$ so that $K_n$ persists on $B_{\nu_n}(f).$
\end{lem}
\begin{pf}
By Theorem \ref{omega limit}, we know that $f$ has a finite number of solenoidal attractors $ C_i.$ Furthermore, $C_i=\omega(c_i)$ for a turning point $c_i$
 at which $f$ is infinitely renormalizable. Each $c_i$  has associated a sequence of restrictive intervals $J_n^i\ni c_i$ of period $2^n$. If we let $K_n$ be as in Theorem \ref{omega limit}, then for $n$ large
$$K_n= \bigcup_i \bigcup_{k=0}^{2^n-1}f^k(J^i_n)$$
The persistence of $K_n$ follows directly from Lemma \ref{renor
  window} by taking $\nu_n>0$ equal to the minimum of the constants
$\epsilon_n$ associated to the intervals $J_n^i$. In addition, if
$g\in B_{\epsilon_n}(f)$ then the continuation 
$J^i_n(g)$
of $J^i_n$ associated to $g$, is a restrictive interval of period $2^n$ and 
$$K_n(g)=\bigcup_i \bigcup_{k=0}^{2^n-1}g^k(J^i_n(g)).$$
\end{pf}

\begin{lem}\label{lem:perturb}
Let  $f\in \Gamma_{\mathcal A^k(I)}$ and let $K_n$ be as in 
Theorem \ref{omega limit}. 
Given $n$ large enough, there exists $\epsilon_n>0$ 
so that $K_n$ and $K_{n+1}$ persist on $B_{\epsilon_n}(f).$ 
Furthermore, let  $K_i^g$ be the continuation of $K_i,$ 
$i=n,n+1,$
associated to $g\in B_{\epsilon_n}(f)$. 
Then for $0\leq j\leq n$,
 any $x\in \Omega'_j:= \Omega(g) \cap cl(K^g_j\setminus K^g_{j+1})$ 
is a periodic point of period $2^m$ for some $0\leq m\leq n$.

\end{lem}
\begin{pf}
In Lemma~\ref{K_n}
we proved that for $n$ sufficiently large, 
there exists $\nu_n>0$ so that 
$K_n$ persists on $B_{\nu_n}(f)$.
Taking $\nu_n$ smaller if necessary,
we can assume that all hyperbolic 
attracting basins for $f$
and all repelling periodic points with period less
than $2^n$
persist over $B_{\nu_n}(f)$

\medskip

\noindent
\textit{Claim 1.} 
Let $K'_0,K'_1$ be the intervals associated to $g$ by
Theorem~\ref{omega limit}.
There exists $\epsilon_0>0$
so that for $g\in B_{\epsilon_0}(f),$
$\Omega^g_0=\Omega(g)\cap\text{cl}(K'_0\setminus K'_1)$
consists of fixed points of $g$.

\medskip

The lemma follows inductively from the claim: Let $J$ be a component
of $K_{n}$ and consider $f^{2^n}|J$.
By the claim there exists $\epsilon_n\leq\epsilon_{n-1}$ so that if
$g\in B_{\epsilon_n}(f),$
$\Omega_n(g)=\Omega(g)\cap\text{cl}(K_n(g)\setminus K_{n+1}(g))$
consists of fixed points of $g^{2^n}$.

\medskip

\noindent
\textit{Proof of Claim 1.}
To conclude the proof of the lemma, we now prove Claim 1.
Let us start by describing how parabolic fixed points bifurcate over
small $\mathcal C^3$ neighbourhoods of $f$. 
Suppose that $p$ is a parabolic periodic point with multiplier 1.
We say that $p$ is of \emph{crossing type},
if on one side of $p$  the graph of $f$ is above the diagonal and on the other it is below.
Parabolic fixed points 
with multiplier -1 always cross the diagonal.

There exists a neighbourhood $Q$ of the set of parabolic points of $f$
such that if $g$ is sufficiently close to $f$, 
every fixed point of $g$ is either in $Q$ or is a continuation of 
a hyperbolic fixed point of $f$,
and each component of $Q$ contains a 
parabolic fixed point of $f.$ 
We denote the component of $Q$ that contains $p$ 
by $Q_p.$ We will show that close the boundary of
$Q_p$, the behaviour of $g$
is similar to the behaviour of $f$
and that in $Q_p$ either there are no periodic 
cycles, a periodic cycle or an invariant interval
for $g$.

% \medskip
% \noindent\textit{Claim.}
% If $p$ is a parabolic fixed point of $f$ with multiplier one, which is not of 
% crossing type, then there exists a neighbourhood $Q_p$ of 
% $f$ so that either $g$ has no fixed points in $Q_p$ or
% or there exists an invariant interval $[q,q']$ for $g$, such that
% the orbit of every point in $[q,q']$ tends to an attracting cycle.

% \medskip
% \noindent\textit{Proof of claim.}

\medskip
\noindent\textbf{Case 1: $p$ is a parabolic 
fixed point of $f$ with multiplier 1, which is not of crossing type.}
Assume that the graph of $f$ is above the diagonal in a neighbourhood of $p$.
Then $p$ is attracting from the left and repelling from the right.
If $Q_p$ contains no fixed point of $g$,
we say that a {\em gate} opens between the graph of $f$
and the diagonal. In this case,
locally, orbits under the perturbed mapping travel from the left of
$p$ to the right, and $g$ has no fixed points in $Q_p$.
So suppose that there is a fixed point of $g$ in $Q_p$.
If there is a non-parabolic fixed point, then since $g$ is close to
$f,$
there are
at least two fixed points for the perturbed map.  Assume
this is the case and let $q$ denote the fixed point in $Q_p$ furthest to the
left and
$q'$ the fixed point in $Q_p$ 
furthest to the right. We have that $q$ is attracting from the 
left, $q'$ is repelling from the right, and $[q,q']$ is an invariant interval 
(if it was not invariant, it would contain a critical point,  but then $g$ would not be close to $f$ in the $\mathcal C^1$ topology). 
The dynamics in the invariant interval are simple, each orbit converges to
a fixed point.
Similar analysis holds when the graph of $f$ is below the diagonal.

\medskip
\noindent\textbf{Case 2: $p$ is a parabolic fixed point of 
crossing type and multiplier 1.}
Either $p$ is attracting
or repelling, and the periodic point persists under small perturbations.
We have that if $p$ is an attracting parabolic fixed point of crossing type for $f$,
then either there is an attracting (not necessarily hyperbolic) fixed point for
$g$ close to $p$, or $g$ has an invariant interval containing no turning points
near $p$ that is attracting from the
left and the right. 
A similar analysis holds when $p$ is a repelling parabolic fixed point
of $f$ with multiplier 1, which is of crossing type:
either there is a repelling (not necessarily hyperbolic) fixed point for
$g$ close to $p$, or $g$ has an invariant interval containing no turning points
near $p$ that is repelling from the
left and the right. 

\medskip
\noindent\textbf{Case 3:
$p$ is
a parabolic fixed point with multiplier $-1$.} 
Then $p$ is of crossing type and
$p$ is a parabolic fixed point with multiplier
1 of crossing type for $f^{2}$, and we can apply the above analysis to
$f^2$
in a small
neighbourhood of $p$.

\medskip

Suppose that there are parabolic fixed points
$p_0, p_1,\dots, p_{k-1}$ 
of $f$ each with multiplier one and not of crossing type
such that for each $i\in\{0,\dots,k-1\}$ 
there is a point $x_i$ such that the following holds:
\begin{itemize}
\item $f^j(x_i)$ converges to $p_{i+1\mod k},$
\item $(f|_{\mathcal{Y}_p})^{-j}(x_i)$ converges to $p_i,$
where $\mathcal{Y}_p$ is the monotone branch of $f$
containing $p$.
\end{itemize}
We will call such a sequence a 
\emph{pseudo-cycle of orbits}.

\medskip
\noindent\textit{Claim 2.} If the entropy of $f$ is zero, then
no such pseudo-cycle of orbits exists.

\medskip

% We will show that We will need the following
% lemma:
% \begin{lem}\cite{M-horse}
% If $f$ is a continuous mapping of the interval, then 
% $f$ has positive topological entropy if and only if some
% iterate of $f$ has a horseshoe.
% \end{lem}
\noindent\textit{Proof of Claim 2.}
Let us recall that if a return mapping to an interval has two full
branches, then it has positive entropy
\cite{M-horse}.
Suppose that $p_j$ is the parabolic fixed point
that is furthest to the right in $I$, and $p_i$ is
furthest to the left.  
If the graph of $f$ is above the diagonal near $p_j$,
then $p_j$ must be repelling from the right.
By assumption, there is a pseudo-cycle of orbits which
enters
$(p_j-\lambda,p_j)$ for any $\lambda>0$. So
the closest turning point to the right of $p_j$, 
$c_1,$  is a local maximum. Furthermore, there are 
no fixed points between $p_j$ and $c_1.$

Let $\alpha$ be the orientation 
reversing fixed point closest to $c_1$. Let $J$
be the interval in $I\setminus\{f^{-1}(\alpha)\}$
that contains $c_1$, then since 
$J$ is not invariant under $f^2$ as there is
pseudo-cycle,
we have that the dynamics of $f^2$ on $J$ 
has positive entropy (the return map has two full
branches).
So we can assume that the graph of $f$ is below the diagonal 
at $p_j$. But then $p_j$ is attracting from the left, and we have that
there is a turning point $c_2$ contained in $[p_i,p_j]$, with
$f(c_2)>p_j.$ But now, the graph of $f$ must cross the diagonal between
$c_2$ and $p_j$, and the point where it crosses cannot be attracting,
since that would violate the condition on orbits near the parabolic point,
so it must be repelling, but now we can argue as before to see that
$f$ must have positive entropy.
So we can assume that there is a turning point between $c_2$ and $p_j$,
this point must correspond to a minimum of $f$, and it must be less that
the parabolic point closest to, and on the left of $c_2$. 
Again we have that $f$ has positive entropy.
So no such parabolic fixed points exist. \checkmark

\medskip

For each repelling periodic point $p$ of $f$, 
let $\varepsilon>0$ be chosen so small that
$B_{\varepsilon}(p)$ is contained in a
neighbourhood of $p$ where $f$ is conjugate
to $x\mapsto f'(p)x$.
Let $U$ be the union of $Q$ 
together with $\cup B_{\varepsilon}(p)$,
where the last union is taken over all repelling 
fixed points of $f$ in the complement of $\mathrm{int}(K_1)$.

Let $\mathcal B$ denote the union of $K_1$,  
basins of hyperbolic attractors,
small neighbourhoods of attracting parabolic points of crossing type
for $f$.
From Proposition~\ref{prop:restrictive intervals},
any point $x$ which is accumulated by 
$f^{-n}(\mathcal B)$, but which in 
not in $f^{-n}(\mathcal B)$ for any $n$,
is a (pre)fixed point of $f$ or is contained in the basin
of a one-sided parabolic attractor, so
we have that for $M$ large enough
$K_0\setminus (\cup_{n=0}^Mf^{-n}(\mathcal B))$
together with
$\cup_{n=0}^{\infty} f^{-n}(U)$
contains all but countably many points of 
$I,$ each whose forward orbit is eventually fixed,
and the complement of
$K_0\setminus (\cup_{n=0}^Mf^{-n}(\mathcal B))$
consists of points that are eventually mapped to
small neighbourhoods (possibly one-sided)
of repelling (not necessarily hyperbolic and 
possibly one-sided) points and 
points that 
converge to a one-sided parabolic attractor.

Suppose first that $K_1$ persists.
Then for any $x\in \cup_{n=0}^M f^{-n}(\mathcal B)$
under $g$ one of the following holds:
\begin{itemize}
\item the orbit of $x$ eventually lands in $K_1(g)$;
\item the orbit of $x$ converges to a hyperbolic attractor;
\item the orbit of $x$ is eventually contained in some $Q_p$ where $p$ is
a parabolic point of $f$ and converges to a fixed point 
of $g$ in $Q_p.$
\end{itemize}
If $x\in \cup_{n=0}^{\infty} f^{-n}(U)$ then either the orbit of $x$
eventually enters $\cup_{n=0}^M f^{-n}(\mathcal B),$ in which case we 
know the possibilities for its forward orbit, or the orbit of
$x$ enters $U$. In this case, either
\begin{itemize}
\item the orbit of $x$ is eventually contained in some $Q_p$ where $p$ is
a parabolic point of $f$ and converges to a fixed point or
\item the orbit of $x$ enters $\cup_{n=0}^M f^{-n}(\mathcal B).$ 
\end{itemize}

So let us assume now that $\partial K_1$
contains a parabolic point $p$ with multiplier 1
and that this point cannot be continued to all nearby mappings.
Then for some nearby map $g$ a gate opens up at the
boundary of $K_1$.
Let $K_1'$ be the union of maximal restrictive intervals of
$g$ and let
$\mathcal B'$ denote the
union of $K_1'$ together with 
\begin{itemize}
\item the corresponding basins of hyperbolic attractors and
\item neighbourhoods of the corresponding attracting parabolic points
  of crossing type.
\end{itemize}
By the analysis above there are no pseudo-cycles of orbits
outside $K_1(g)$,
so we have that an orbit travels through a 
bounded number of gates,
and eventually passes through one that has the property
that any fundamental domain for the dynamics is covered (except for
possibly finitely many points) by
$\cup_{n=0}^M g^{-n}(\mathcal B')\cup\cup_{n=0}^\infty f^{-n}(U)$.
In particular,
every point eventually converges to a fixed point 
for $g$ or enters $K_{1}'$. Thus Claim 1 follows.
\end{pf}

Now we prove Theorem~\ref{thm:main bis} for analytic mappings,
and thus obtain Theorem~\ref{thm:main}, see the end of
Section~\ref{sec:entropy and renormalization}.

\begin{thm}
\label{thm:AE}
Suppose that $f\in\mathcal A_{\underline b}(I)\cap \mathcal{B}_{\Omega_a},$ for some $a>0,$ with all critical points of 
even order,
which is
infinitely renormalizable
at some critical point $c$ and that 
$h(f)=0$. Then there exist mappings
$g,\tilde g\in\mathcal{A}_{\underline b}(I)\cap \mathcal{B}_{\Omega_a},$
arbitrarily close to $f,$
such that
\begin{itemize} 
\item $g$ has finitely many periods, and 
\item $h(\tilde g)>0.$
\end{itemize} 
\end{thm}
\begin{pf}
 Let $J_n\owns c$ be
the sequence of restrictive intervals with periods 
$2^n$ about $c$. By Theorem~\ref{thm:complex bounds},
for all $n$ sufficiently large
there exists a polynomial-like mapping of type $\underline b_1$,
$F\colon U\to V$, $U\owns c,$
$J_n\subset U$ and $F=f^{2^n}|_U$,
where $\underline b_1$ is depends on $\underline b$ and the combinatorics of 
the renormalization.
Moreover, there exists a neighbourhood $\mathcal U\subset \mathcal
B_{\Omega_a}\cap\mathcal A_{\underline b}(I)$
of $f$, so that the polynomial-like mapping
$F\colon U\to V$ persists over $\mathcal U$.
Observe that if $\mathcal U$ is sufficiently small, then
for each $g\in\mathcal U, g^s=:G\colon U_G\to V,$ where
$U_G=\comp_{c_1}(G^{-1}(V))$ is a polynomial-like mapping.
Let $\mathcal R\colon \mathcal U\to\mathcal{PG}_{\underline b_1}$
be the renormalization operator from $\mathcal U$
to the space of polynomial-like germs of 
type $\underline b_1$,
mapping $g\mapsto G$ where $G=g^{2^n}|_{U_G}.$
Let $b_1=|\underline b_1|.$

Since $\mathcal R$ is a composition of affine
rescalings and composition of analytic mappings,
$\mathcal R$ is analytic.
% QQ
% By Lemma~\ref{K_n}, $\mathcal R(\mathcal U)$
% contains an open set $\mathcal U'$
% containing $F$.
% QQ
By Proposition~\ref{prop:bcpl}, for any $b_1$-dimensional transverse family
$F_{\lambda},\ \lambda\in\mathbb D^{b_1}_{\varepsilon}$ and $\varepsilon>0,$
with
$F_{0}=F,$
there exists a real polynomial-like mapping $G\colon U_G\to V_G\in F_{\lambda},$
arbitrarily close to $F$ so that
$G|_{U_G\cap\mathbb R}$ has positive topological entropy.
By Proposition~\ref{prop:transdim} and Theorem~\ref{thm:R is trans nonsing},
we have that there exist $b_1$ vectors $\{w_1,\dots, w_{b_1}\}$, transverse to the topological conjugacy
class of $f$, so that for $\varepsilon>0,$ small, the family
$\{\lambda_1w_1+\dots+\lambda_{b_1}w_{b_1}:\lambda_i\in\mathbb
D_{\varepsilon}\mbox{ for }i\in\{1,\dots, b_1\}\}$ maps to a $b_1$-dimensional
family transverse to $\mathcal H_f.$ We may assume that $G$ is contained in this transverse family.
Thus, by continuity of $\mathcal R$, there exists an analytic mapping 
$g\in\mathcal U,$ which is a preimage of $G$ under $\mathcal R.$
The mapping $g$ has positive topological entropy,
since its renormalization $G|_{U_G\cap\mathbb R}$ has positive topological entropy.

Showing that there is a sequence with zero topological entropy 
is a little harder, we need to ensure that the preimage under
$\mathcal R$ still has zero entropy, and we need to consider all
turning points at which $f$ is infinitely renormalizable.

Let $c_i,$ $1\leq i\leq m,$ denote the critical points of $f$ such that
$\omega(c_i)$ is a solenoidal attractor. 
Let $m'$ be the number of distinct such solenoidal attractors.
For each distinct $\omega(c_i)$ 
choose a critical point $c_{i,0},$ $1\leq i\leq m',$ of even order so that 
$\omega(c_{i,0})=\omega(c_i)$ and $f$ is infinitely renormalizable
at $c_{i,0}$.
Since $f$ has at most $|\underline b|$ critical points,
and at any critical point at which $f$ is renormalizable,
the period of the restrictive interval is a power of 2,
by Theorem~\ref{thm:complex bounds},
there exists a neighbourhood
$\mathcal U\subset\mathcal{A}_{\underline b}(I)\cap\mathcal{B}_{\Omega_a}$
of $f$ and $N\in\mathbb N$,
so that the following
holds: 
For each $c_{i,0},$ there exists a $b_i$-tuple, $\underline b_i$,
depending on $\underline b$ and the combinatorics of the
renormalization, so that
the mapping $f^{2^N}\colon J_N^i\rightarrow J_N^i$ extends
to a polynomial-like mapping of type
$\underline b_i$, $F_i\colon U_{F_i}\to V_{F_i}$, which persists
over $\mathcal U$, where $J_N^i\owns c_{i,0}$ is the restrictive interval of
period $2^N$ containing $c_{i,0}$.
Let $\mathcal R_i\colon \mathcal U\to\mathcal{PG}_{\underline b_i},$
and let $\hat{\mathcal{R}}\colon \mathcal{U}\to \mathcal{PG}_{\underline
  b_1}\times\dots\times\mathcal{PG}_{\underline b_{m'}}$
be the mapping defined by 
$\hat{\mathcal{R}}(f)=(\mathcal R_1(f),\dots,\mathcal R_{m'}(f)).$ 
We have that $\hat{\mathcal{R}}$ is continuous 
(it is a composition of iteration and rescaling in each coordinate).
%and $\hat{\mathcal{R}}(\mathcal U)$ is open, see \cite[Remark 2.7]{ALdM}.
By Proposition~\ref{prop:transdim}, for each $i$, 
there exist normal vectors $v_{i,1},\dots, v_{i,b_i},$ so that for $1\leq
j<j'\leq b_i,\  \|v_{i,j}-v_{i,j'}\|$ is bounded away from zero  in
$T_{F_i}\mathcal{PL}_{\underline b_i}$, which are transverse to the topological
conjugacy class of $F_i$.
% $an open neighbourhood $\mathcal V_i$ of 
% $F_i$ in $\mathcal{PG}_{\underline b_i}$,
Since $D\mathcal R_i$ is transversally non-singular (Theorem~\ref{thm:R is trans nonsing}),
we have that for each $1\leq  j\leq b_i,$ there exists a sequence of vectors
$w_{i,j}^k\in T_f\mathcal A_{\underline b},$ so that $D\mathcal R_i(f)w_{i,j}^k$
converges to $v_{i,j}$ as $k\to\infty.$ Thus for $k$ sufficiently large,
we have that $D\mathcal R_i(f)w_{i,j}^k$ is transverse to the topological conjugacy class
of $F_i,$ and the family $\{D\mathcal R_i(f)w_{i,1}^k,\dots,D\mathcal R_i(f)w_{i,b_i}^k\}$ spans a
$b_i$-dimensional space transverse to the topological conjugacy class of $F_i.$
By  Proposition~\ref{prop:bcpl} we have that any analytic
family of mappings in $\mathcal{PL}_{b_i}$ that is transverse to $\mathcal H_{F_i}$
contains polynomial-like germs in the
interior of the set of mappings with zero entropy.
Let $f_\lambda,\ \lambda\in\mathbb D^{b_1+b_2+\dots +b_{m'}}_{\varepsilon},$
be the family
$f+\lambda_1w^k_{1,1}+\lambda_2w_{1,2}^k+
\dots+\lambda_{b_1+b_2+\dots+b_{m'}}w^k_{m',b_{m'}}.$
Now, from the choice of the $w_{i,j}^k,$ there exist $\lambda\in \mathbb D^{b_1+b_2+\dots +b_{m'}}_{\varepsilon}$ with
all coordinates arbitrarily close to zero so that
the mapping $\mathcal R_if_\lambda$ is in the interior of mappings with zero entropy in
$\mathcal{PL}_{\underline b_i}.$ Let $N\in\mathbb N,$ and
%which intersects the interior of the set of mappings
%with zero entropy in $\mathcal{PG}_{\underline b_i}$.
%Moreover, we can choose $\mathcal U$ small enough that
%Lemma~\ref{lem:perturb} holds.
%Let $\mathcal U'$ be the preimage of
%$(\mathcal V_1,\dots ,\mathcal V_{m'})$ under $\hat{\mathcal{R}}$.
%Let $g\in\mathcal U\cap\mathcal U'$
%be so that if $\hat{\mathcal{R}}(g)=
%(G_1,\dots, G_{m'})$, then each 
%$G_i,1\leq i\leq m',$ has finitely many periods.
let $K_N$ be the forward invariant set
from Theorem~\ref{omega limit}.
It is the union of restrictive intervals of period $2^N$ 
for $f_\lambda.$
By Lemma~\ref{lem:perturb}, we have that
the set of periodic points of $f_{\lambda}$ in 
$I\setminus K_N$ has finitely many periods, and
we have constructed $f_\lambda$ so that its set of 
periodic points in $K_N$ also has finitely many periods.
Thus, since $K_N$ is forward invariant,
for some $\lambda\in \mathbb D_{\varepsilon}^{b_1+\dots + b_{m'}},$
arbitrarily small,
$f_\lambda$ has finitely many periods.
% there is an analytic mapping $g$
% with entropy zero arbitrarily close to $f$,
% so that each critical point $c$ of $g$
% is periodic.
% So there is an analytic mapping $g_1$ arbitrarily close to $f$,
% so that the corresponding critical point $c_{i,0}(g)$ of $g$
% is periodic and for every critical point $c$ with
% $\omega(c)=\omega(c_{i,0})$ is we have that 
% for some iterate $k_c,$ $g^{k_c}(c(g))=c(g)$.
% Let $C_1$ denote this subset of critical points.
% Thus there exists a neighbourhood $\mathcal U_1$ of 
% $f$ so that for every $g\in \mathcal U_1$ for each critical point
% $c\in C_1,$ $c(g)$ is in the basin of a hyperbolic or super-attracting 
% attractor. Repeating this for each of the finitely many 
% distinct solenoidal attractors $\omega(c_{i,0})$ we conclude the proof.
\end{pf}

% \begin{lem}
% \label{lem:transverse non-singular}
% Suppose that $f\in\Gamma_{\mathcal A_{\underline b}(I)}$,
% and let $\mathcal R(f):U\to V$ be a polynomial-like renormalization 
% of $f$. Suppose that $v\in T_{R(f)}\mathcal{PG}_{\underline b}$,
% which is transverse to the hybrid class of $R(f)$.
% Then there exists a sequence $w_i\in T_f\mathcal A_{\underline b}(I),$
% so that $DR(f)w_i\rightarrow v$.
% \end{lem}

\subsection{Proof of Theorem~\ref{thm:II}}
Theorems~\ref{thm:B1} and \ref{thm:B2} below
imply
Theorem~\ref{thm:II}.

%\medskip
%\noindent\textit{Proof of Theorem~\ref{thm:II}.}
% By Theorem~\ref{thm:main},
% for any open set 
% $\mathcal U\subset\mathcal{A}_{\underline b}(I)$
% with $\Gamma\cap\mathcal U\neq\emptyset$, 
% $\mathcal U\setminus \Gamma$
% is disconnected.
% Thus $\Gamma$ has codimension-one in 
% $\mathcal{A}_{\underline b}(I)$
% (this is a definition of codimension-one
% from geometric topology, which is suitable in this context since
% $\Gamma$ is not in general a manifold).
% Evenmore, observe that by Proposition~\ref{prop:bcpl} in the space of 
% $\mathcal{PG}_{\underline b},$
% for any $f\in \Gamma_{\mathcal{PG}_{\underline b}}$, there exists
% a vector $v$, transverse to the hybrid class of $f$ so that
% for any $\varepsilon>0$, $f-\varepsilon v$ is in the interior of 
% mappings with zero entropy.
% Thus by Lemma~\ref{lem:transverse non-singular}
% and Lemma~\ref{lem:perturb}
% for any $f\in\Gamma_{\mathcal{A}_{\underline b}(I)},$
% there exists a vector $w$
% transverse to the topological conjugacy class of $f$,
% so that for any $\varepsilon>0$,
% $f-\varepsilon w$ is in the interior of zero entropy.

\begin{thm}
\label{thm:B1}
There exists an open and dense subset of 
$\Gamma_{\mathcal{A}_{\underline b}(I)}$  which is
contained in the basin of a
unimodal, polynomial-like fixed point of 
renormalization. 
\end{thm}
\begin{pf}
Let $\Gamma_1$ denote the subset of $\Gamma$ 
consisting of mappings with exactly one solenoidal attractor.
Let us show that $\Gamma_1$ is open and dense in
$\Gamma$. Suppose that $f\in\Gamma_1$. Then 
$f$ has a single solenoidal attractor and the
critical points whose orbits do not converge to the solenoidal attractor
are asymptotic to periodic points of 
period $2^n$, where $n$ is bounded from above. 
Thus in any sufficiently small neighbourhood of $f$,
each mapping has at most one solenoidal attractor. Thus $\Gamma_1$ is open in $\Gamma.$
Let us now show that $\Gamma_1$ is dense in $\Gamma$.
Suppose that $f\in\Gamma\setminus\Gamma_1.$
We need to show that we can approximate $f$ by
mappings with a single solenoidal attractor.

We can argue as in the proof of Theorem~\ref{thm:AE}.
Let $f$ be an analytic mapping with at least two 
solenoidal attractors. For ease of exposition, assume that 
$f$ has exactly two solenoidal attractors.
Then there exist vectors $\underline b_1$ and $\underline b_2,$ 
and a neighbourhood $\mathcal U$ of $f$ so that
$\hat{\mathcal{R}}\colon \mathcal U\to\mathcal{PG}_{\underline
  b_1}\times\mathcal{PG}_{\underline b_2}.$
Let $\hat{\mathcal{R}}(f)=(G_1,G_2).$
By Proposition~\ref{prop:bcpl},
there exist mappings $G$ arbitrarily close to $G_2$
in the interior of zero entropy.
Thus, since $\hat{\mathcal{R}}$ is a continuous
mapping, we can argue as in the proof of Theorem~\ref{thm:AE},
to see that
there exists an analytic mapping $g$ arbitrarily close to $f$
with exactly one solenoidal attractor and zero entropy.
%Moreover, as in the proof of Theorem~\ref{thm:AE},
%by Lemma~\ref{lem:perturb},
%we can find such a $g$ with zero entropy.

Mappings in $\Gamma_1$ could have several critical points
in their solenoidal attractors.
We will now show that there is an 
open and dense set $\Gamma_2$ of $\Gamma_1$
consisting of mappings such that
there is only one critical point in the 
solenoidal attractor.
The proof that $\Gamma_2$ is open (that is, relatively open in
$\Gamma)$ is the same as the proof that $\Gamma_1$ is open,
and so we omit it.
To prove that
$\Gamma_2$ is dense, we use the strategy used to 
prove Theorem~\ref{thm:AE}.
First, let us show that
in the space of stunted sawtooth mappings  we can approximate
any mapping in $\Gamma_{\mathcal S_b}$
by mappings $T_k$ in $\Gamma_{\mathcal S_b}$
with one recurrent, non-periodic plateau.
%remaining plateaus are either periodic
%all but one plateau 
%periodic. 
If $h(T)=0$, and $T$ is at most finitely renormalizable at each plateau, 
let $T'$ be the last renormalization of $T$.
Then, by \cite[Lemma 7.6]{BruvS} 
the $\omega$-limit set of each point
under $T'$ is a fixed point of $T'.$
Moreover, since this fixed point is necessarily
attracting, it is contained in a fixed plateau of $T'$.
By \cite[Lemma 7.7]{BruvS}, if $T$ is a stunted sawtooth mapping
in the interior of zero entropy, then each point under $T$
is either (pre)periodic or in the basin of one of the
periodic attractors (periodic plateaus) of $T$.
By \cite{TH}, see Theorem~\ref{thm:boundary ssm},
$\Gamma_{\mathcal S_m}$ is the limit of 
stunted sawtooth mappings
with periodic plateaus of period $2^n$ as 
$n\rightarrow\infty$.

\medskip
\noindent\textit{Claim.}
For any $\varepsilon_0>0$, there exists $n\in\mathbb N\cup\{0\}$, so that
if $T=(t_1,t_2,\dots,t_b)$ has a periodic plateau of period $2^n$
and zero entropy,
for some $t_i$ and some $\varepsilon\in(0,\varepsilon_0),$
either
$T_1=(t_1,\dots,t_{i-1},t_i+\varepsilon,t_{i+1},\dots,t_b)$
or $T_2=(t_1,\dots,t_{i-1},t_i-\varepsilon,t_{i+1},\dots,t_b)$
is in $\Gamma$.

\medskip
\noindent\textit{Proof of Claim.}
Observe that the space of stunted sawtooth mappings is
compact and recall that period-doubling bifurcations
occur in each parameter separately.
% We have that for any $\eta>0$, there exists $j\in\mathbb N$,
% so that the set of all stunted sawtooth mappings with zero
% entropy and a periodic
% plateau of period $2^j$ is contained in an 
% $\eta$-neighbourhood of $\Gamma$.
% If not, then there exists $\eta>0$,
% so that for any $j$, there exists a stunted sawtooth mapping
% with zero entropy and a periodic plateau of period $2^n$,
% which is not contained in an $\eta$-neighbourhood of $\Gamma$.
% But this contradicts the fact that $\Gamma$ is set of limits of 
% period-doubling bifurcations in the space of stunted sawtooth
% mappings.
Suppose the claim fails. 
Then there exists $\varepsilon_0>0$, so 
that for any $n\in\mathbb N\cup\{0\},$
there exists $T=(t_1,t_2,\dots,t_b)$ with
a periodic plateau of period $2^n$, zero entropy, 
so that for each $i$, we have that
for each $\varepsilon\in(0,\varepsilon_0),$
$T=(t_1,\dots,t_{i-1},t_i\pm\varepsilon, t_{i+1},\dots,t_b)$
is not in $\Gamma$. 
Since there are at most finitely many plateaus,
this implies that for some $i\in\{1,\dots,b\},$
there is a sequence of
stunted sawtooth mappings
$T_{n}=(t_1^n,\dots,t_b^n)$ with the $i$-th plateau periodic with
period $2^{j_n}$, and no plateau periodic with period 
greater than $2^{j_n}$,
where $j_n\to\infty$ as $n\to\infty$.
Since period-doubling bifurcations occur sequentially in 
the space of stunted sawtooth mappings,
we can assume that the parameters $t^n_i$ are monotone.
Thus they converge to a limit $t_*$.
This limiting parameter is accumulated by periodic points of
period $2^n$. Since $|t_*-t^n_i|\rightarrow 0,$ we 
arrive at a contradiction and the claim follows.
\checkmark

Now, by Theorem~\ref{thm:main},
by taking $N$ large 
we can approximate $T$ 
arbitrarily well by mappings with 
$\mathcal Per(T)=\{2^n:0\leq n\leq N\}.$
But now, by the claim, we can perturb such a mapping
by moving just one plateau up or down to obtain
a mapping in $\Gamma$, moreover, the size of
this perturbation tends to zero as $N\rightarrow\infty$.

To conclude the proof,
we can argue as in the 
proof of Theorem~\ref{thm:AE}. 
Let $f$ be an analytic mapping with exactly one solenoidal
attractor, and let $c$ be a critical point in the 
solenoidal attractor.
Then for some $b'\in\mathbb N,$
there is a $b'$-tuple, $\underline b'$, so that
$f$ has a polynomial-like renormalization 
of type $\underline b'$,
$F\colon U\to V,$ about $c$. Let $P=\chi(F)$ be its straightening.
Then by Theorem~\ref{thm:boundary ssm},
$\Psi(P)$ is a stunted sawtooth mapping in 
the boundary of mappings with finitely many periods.
Recall the definition of $\Psi$ on page~\pageref{def:psi}.
By the claim, we can approximate
$\Psi(P)$ by stunted sawtooth 
mappings $T_j\in\Gamma$ with exactly one plateau 
in a solenoidal attractor.
Thus arguing as in the proof
of Theorem~\ref{thm:boundary poly},
we can approximate
$P$ by polynomials $P_j$ of type $\underline b'$
with exactly one solenoidal attractor, which contains exactly one
critical point. So as in the proof of
Proposition~\ref{prop:bcpl},
there exist polynomial-like germs
converging to $F$, which are hybrid conjugate to
the $P_j$, and finally, as in the proof of
Theorem~\ref{thm:AE}, we can approximate $f$, by
analytic mappings in $\Gamma_2.$
\end{pf}

\begin{thm}\label{thm:B2}
Then $\Gamma_{\mathcal{A}_{\underline b}(I)}$ admits a
cellular decomposition.
\end{thm}

\begin{pf}
  Let $\Gamma_2$ be the (relatively) open and dense subset of $\Gamma$ given by
  Theorem~\ref{thm:B1}. Let $X$ denote a connected component of $\Gamma_2.$
  We need to show that there is a relatively open and dense subset of
$\partial X$ consisting of codimension-2 cells.
Let $X_1$ denote the subset of $\partial X$ consisting of mappings 
with a single solenoidal attractor containing exactly 2 critical
points and let $X_2$ denote the subset of $\partial X$ consisting of
mappings with exactly two solenoidal attractors each containing
exactly one critical point.

\medskip
\noindent\textit{Claim 1.}
 $X_1\cup X_2$ is open and dense in $\partial X$.

\medskip
\noindent\textit{Proof of Claim 1.}
First we show that $X_1$ and $X_2$ are open in $\partial X$.
Suppose that $f\in X_1$. Then, relabelling the critical points of $f$
if necessary, we can assume that $f$ has a solenoidal attractor
which contains $c_1$ and $c_2$, but not $c_3,\dots,c_b$.
Let $\mathcal{J}_n$ denote the cycle of the restrictive interval of period 
$2^n$. For $n$ sufficiently large,
$\mathcal J_n\cap \{c_3,\dots,c_b\}=\emptyset.$
Thus, by Lemma~\ref{renor window} for $n$ sufficiently large,
there is an open set of mappings $\mathcal U$ containing $f$,
such that for all $g\in\mathcal U,$ each $g$ has a restrictive interval
of
period $2^n$ and the orbit of this interval contains exactly two 
critical points of $g$.
For mappings $g\in\mathcal U\cap\partial X$,
the number of critical points
in the solenoidal attractor cannot drop to one, since the condition
that $f$ have a solenoidal attractor containing exactly one critical
point is relatively open in $\Gamma$. Thus we have that
$X_1$ is relatively open in $\partial X.$
The proof that $X_2$ is relatively open is similar,
just consider two disjoint cycles of restrictive intervals with
sufficiently high period.

Let us now explain how to see that
$X_1\cup X_2$ is dense in $\partial X$.
Suppose that $f\in\partial X.$ If $f$ has 
exactly one solenoidal attractor (which must contain at 
least two critical points) then we show that
we can approximate $f$ by mappings with a single
solenoidal attractor, which contains exactly two critical points.
If $f$ has more than one solenoidal attractor, then
we show that we can approximate $f$ by 
mappings with two solenoidal attractors,
each containing exactly one critical point.
The strategy for carrying out these approximations
is no different than in the proof that
codimension-one cells (consisting of mappings
with a solenoidal attractor containing exactly one
critical point) are dense in $\Gamma$, and so we omit the 
details. One first proves the corresponding 
statement in the space of stunted sawtooth mappings, and then
transfers it successively to polynomials, polynomial-like
germs and finally to analytic mappings.\checkmark

\medskip
\noindent\textit{Claim 2.} Each of
 $X_1$ and $X_2$ have codimension-two
in $\mathcal{A}(I)$. 

\medskip
\noindent\textit{Proof of Claim 2.}
Suppose that $f\in X_1$.
Then $f$ has a renormalization $R(f)=f^s|_U=F\colon U\to V$
that is contained in a space of polynomial-like germs
$\mathcal{PG}_{\underline b'}$
with exactly two critical points, and indeed
there is an open set $\mathcal U\owns f,$
such that if $g\in\mathcal U$, then $g$ has a 
polynomial-like renormalization
$g^s|_{U_g}=G\colon U_G\to V_G$ in 
$\mathcal{PL}_{\underline b'}.$
The codimension of the hybrid class of $F$ in the space
of polynomial-like germs is two.
Thus we have that there are vectors
$v_1, v_2\in T_{F}\mathcal{PG}_{\underline b'},$
which are transverse to the hybrid class of $F$,
and since $F\in\Gamma_{\mathcal{PG}_{\underline b'}},$
we have that we can choose these vectors so that 
for $t>0$ and small, $F-t v_1$ is 
in the interior of zero entropy and
$F-t v_2$ in $R(\mathcal U\cap X)$.

Suppose that $v\in T_{R(f)}\mathcal{PG}_{\underline b}$,
which is transverse to the hybrid class of $R(f)$.
Then by Theorem~\ref{thm:R is trans nonsing},
% there exists a sequence $w_i\in T_f\mathcal A_{\underline b}(I),$
% so that $DR(f)w_i\rightarrow v$,
% see \cite[Lemma 4.8]{ALdM} and \cite[Theorem 3]{Smania-shy} 
% for a generalization to the multimodal case.
% So the
% renormalization operator is transversally non-singular,
% and
we have that there exist
vectors $w_1,w_2\in T_f\mathcal{A}(I),$
so that
$w_1$ and $w_2$ are transverse to $\partial X$.
If $f\in X_2$, the proof is similar - consider the renormalizations
about each of the critical points separately. \checkmark

Proceeding inductively we see that the union of codimension-$j$ cells in
$\Gamma_{\mathcal{A}(I)}$, 
where $j$ runs from $1$ though to $b$ exhausts $\Gamma$.
\end{pf}

\noindent
\textit{Remark.}
Let us describe the finer structure of the set $X.$
Let $\Gamma_3$ be the subset of $\Gamma_2$ that consists of mappings $f$
with critical points $\{c_1,\dots,c_{b}\},$ so that exactly one critical
point, say $c_1,$ is recurrent and the remaining $b$ critical points are
asymptotic to periodic points.
Arguing just as we did to see that $X_1$ is open, we can show that
$\Gamma_3$ and $\Gamma_2\setminus\Gamma_3$ are relatively open in $\Gamma.$
If $X$ is a connected component of $\Gamma_3,$ then we can decompose $X$ into a
countable union of $\cup Y_i,$ where each $Y_i$ consists of mappings with exactly one
solenoidal attractor, which only attracts $c_1,$ and where the remaining critical points
converge to hyperbolic attractors that persist over $Y_i$. Each $Y_i$ is a
codimension-one set in $\mathcal A_{\underline b}(I)$, and 
$X\setminus (\cup
Y_i)$ consists of mappings with parabolic cycles, which have codimension at
least 2. %One can construct a transverse vector for a mapping with a parabolic
%cycle using the argument of \cite[Lemma 7.1]{ALdM}.

Now, let $X$ be a component of $\Gamma_2\setminus\Gamma_3.$
Then if $f\in X$ with critical points $\{c_1,\dots,c_{b}\},$ we can assume that
$c_1\in\omega(c_1)$ where $\omega(c_1)$ is a solenoidal attractor, which contains
no other critical points, and that there exists $b_1\in\{2,\dots,b\}$ so that for $i\in\{2,\dots b_1\},$ $\omega(c_i)=\omega(c_1)$
and for $i\in\{b_1+1,\dots,b\},$ $c_i$ tend to a periodic orbit.
Using the same argument as in the previous paragraph, we can assume that
 $f\in Y,$ a subset of $X,$ over which periodic orbits that attract $\{c_{b_1+1},\dots,
c_b\}$ do not bifurcate.
For any small perturbation $g\in Y$ of $f,$ let
$c_i(g)$ denote the critical point of $g$ corresponding to $c_i.$
Then, since $g\in Y\subset \Gamma,$
we have that $\omega(c_1(g))\owns c_1(g)$ is a solenoidal attractor.
Moreover, for $n$ sufficiently big $J_n(g),$ contains no attracting cycle, where
$J_n(g)$ is the restrictive interval of period $2^n$ for $g.$
Since for $i\in\{2,\dots, b_1\},$ the orbit of $c_i$ enters every periodic interval
about $c_1,$
we have that either $c_i(g)$ eventually lands on $\partial J_{n'}(g)$ for some
  $n'$ large or  $\omega(c_i(g))=\omega(c_1(g))$.
  Each of these defines a codimension 2 condition.

  \medskip

\subsection{Boundary of chaos for smooth mappings.}

In this section, we prove Theorems~\ref{thm:main smooth}
and \ref{thm:II smooth},
which extend Theorems~\ref{thm:main} and \ref{thm:II} 
to smooth mappings.

% Suppose that $\Lambda$ is a $\mathcal C^k$-manifold with
% base point $\lambda_0\in\Lambda.$
% Suppose that $X\subset \mathbb C,$ we say
% that $h_{\lambda}:X\to\mathbb C,\lambda\in\Lambda,$
% is a $\mathcal C^k$-
% \emph{motion} if 
% \begin{itemize}
% \item $h_{\lambda_0}=\mathrm{id},$
% \item $h_{\lambda}$ is an injection for each $\lambda\in\Lambda,$ and
% \item $z\mapsto h_{\lambda}(z)$ is $\mathcal C^k$ in $\lambda$.
% \end{itemize}

Suppose that $f\in\mathcal{A}_{\underline b}^k(I)$
and
let $\mathcal W$ be a neighbourhood of $f$ in
$\mathcal{A}_{\underline b}^k(I)$.
%We say that 
%a periodic interval $K=K^f\subset [-1,1],$ 
%of period $s$
%\emph{persists} in $\mathcal W$
%if for each $g\in\mathcal W$, there is a $\mathcal C^k$-motion
%$h^g:K^f\to K^g$, and $K^g$ is a periodic interval of period $s$
%and $h^g\circ f^s(z)=g^s\circ h^g(z)$ for $z\in\partial K$.
%We call $K^g$ the \emph{continuation of} $K^f$ to $g.$
If $f^s|_{U}=F\colon U\to V$ is an 
asymptotically holomorphic polynomial-like mapping,
we say that 
$F\colon U\to V$ persists over $\mathcal W$ if for each 
$g$ in $\mathcal W$ there is a $\mathcal C^k$-motion
$h^{g}\colon (U,V)\to (U_G,V_G),$ 
an asymptotically holomorphic 
polynomial-like mapping, $g^s|_{U_G}=G\colon U_G\to V_G,$
and
$h^g\circ F(z)=G\circ h^g(z)$ for $z\in\partial U.$

Before proving Theorem~\ref{thm:main smooth}, 
let us collect some general tools.

\begin{lem}\cite[Proposition 5.5]{GdM}
\label{lem:ahpl approx}
For any bounded domain $U$ in the 
complex plane, there exists $C=C(U)>0$,
with $C(U)\leq C(W)$ if $U\subset W$,
such that the following holds:
Let $\{G_n\colon U\to G_n(U)\}_{n\in\mathbb N}$ be  sequence
of quasiconformal homeomorphisms such that
\begin{itemize}
\item the $G_n(U)$ are uniformly bounded; that is, there
exists $R>0$ so that $G_n(U)\subset B_R(0)$ for all $n$; and
\item $\mu_n\rightarrow 0$ in $L^\infty,$ where $\mu_n$ is
the Beltrami coefficient of $G_n$ in $U$.
\end{itemize}
Then given any domain $U'\Subset U,$ there exists
$n_0\in\mathbb N$ and a sequence 
$\{H_n\colon U'\to H_n(U')\}_{n\geq n_0}$ of biholomorphisms such that
$$\|H_n-G_n\|_{\mathcal C^0(U')}\leq C(U)\Bigg(\frac{R}{d(\partial
  U,\partial U')}\Bigg)\|\mu_n\|_{\infty},$$
where $d$ is the Euclidean distance between the disjoint sets
$\partial U$ and $\partial U'$.
\end{lem}

\begin{lem}\cite[Proposition 11.2]{GMdM}
\label{lem:rbound}
Let $I$ be a compact interval in the real line
and let $U$ be an open set in the complex plane 
containing $I$. Fix $M>0$ and consider the family
$$\mathcal F=\{f\colon U\to\mathbb C, \mathrm{ holomorphic }:
\|f\|_{\mathcal C^0}\leq M\}.$$
Then for any $k\in\mathbb N$ and any $\alpha\in(0,1)$
there exists $L=L(k,\alpha,M)>0$ such that
$$\|f\|_{\mathcal C^k(I)}\leq L(\|f\|_{\mathcal C^0})^{\alpha}.$$
\end{lem}
Combining Lemmas~\ref{lem:ahpl approx} and
\ref{lem:rbound} we obtain a bound on a $\mathcal C^k$
norm from a bound on the dilatation of
Beltrami differential.

We say that a diffeomorphism $\phi\colon I\to I$
is in the \emph{Epstein class}, $\mathcal E_{\beta},$
if there exists $\beta>0$, so that
$\phi^{-1}$ extends to a holomorphic, univalent mapping
from the slit complex plane $$\mathbb C_{(-1-\beta, 1+\beta)}
=\mathbb C\setminus((-\infty, -1-\beta]\cup [1+\beta,\infty))$$
into $\mathbb C.$
Given a set $P=\{p_1,\dots,p_b\}$ of $b$
real unimodal polynomials which 
preserve the interval $[-1,1]$,
we say that a (multimodal) mapping $f\in\mathcal A(I)$
of the interval is in the
\emph{Epstein class}, $\mathcal E_{\beta,P}$
if it can be expressed in the form
$$f=\phi_j\circ p_j\circ\phi_{j-1}\circ p_{j-1}\circ\dots\circ\phi_1\circ p_1,$$
with $j\leq b,$ where 
each 
$\phi_j$ is in $\mathcal E_{\beta}.$

\begin{lem}\cite[Theorem 2]{ShT}\label{lem:Shen-Todd}
Suppose that $f\in\mathcal{A}^k(I),$ $k\geq 2.$
let $T$ be an open interval such that $f^s\colon T\to f^s(T)$
is a diffeomorphism. Then for any 
$S, \alpha, \varepsilon>0,$
there exists $\delta=\delta(S,\alpha,\varepsilon)>0$
and $\beta=\beta(\alpha)>0$
satisfying the following.
Suppose that $\sum_{j=0}^{s-1}|f^j(T)|\leq S$
and that $J$ is a closed subinterval of $T$
such that 
\begin{itemize}
\item $f^s(T)$ is a $\alpha$-scaled neighbourhood of $f^s(J)$, and
\item $|f^j(J)|<\delta$ for $0\leq j< s.$
\end{itemize}
Then letting $\psi_0\colon J\to I$ and $\psi_s\colon f^s(J)\to I$ be affine
diffeomorphisms, there exists a mapping $G\colon I\to I$ in the 
Epstein class $\mathcal E_{\beta}$ such that
$$\|\psi_s\circ f^s\circ\psi_0^{-1}-G\|_{\mathcal C^k}<\varepsilon.$$ 
\end{lem}

To simplify the statements of the following two results about 
Epstein mappings, let us fix a set of $b$ real unimodal polynomials, 
$P=\{p_1,\dots, p_b\}$, with the property that each $p_i$ preserves the interval
$[-1,1]$. Let $\hat{\mathcal{E}}_{\beta}=\mathcal E_{\beta,P}.$ 

By Lemma~\ref{lem:Shen-Todd} and Theorem~\ref{thm:real bounds}, we have
\begin{lem}
\label{lem:Epstein class}
There exist $\beta\in(0,1)$ so that the following holds.
Let $\varepsilon>0$.
Given any mapping
$f\in\Gamma_{\mathcal{A}^k(I)},$ which is infinitely renormalizable at a
critical point $c_0,$
there exists $j_0\in\mathbb N$
and a sequence of Epstein mappings $H_j$ in 
$\hat{\mathcal{E}}_{\beta}$
with the same domain as $\mathcal R_{c_0}^j(f),$
such that for $j\geq j_0$,
$$\|\mathcal R_{c_0}^j(f)-H_j\|_{\mathcal C^{k}(I)}\leq \varepsilon.$$
\end{lem}

\begin{lem}\label{lem:epdom}
For any $\beta\in(0,1)$ and $b\in\mathbb N$,
there exists a Jordan domain $U_{\beta}$ containing $I=[-1,1]$
and a positive constant $M_{\beta}$
so that for any Epstein mapping 
$g\in\hat{\mathcal{E}}_{\beta}$
of $I$
the holomorphic extension of $I$
is well-defined in $U_\beta$
and satisfies $|g(z)|\leq M_{\beta}$
for all $z\in U_{\beta}$.
\end{lem}
\begin{pf}
Since each mapping in the Epstein class $\hat{\mathcal{E}}_{\beta}$
can be expressed as a composition:
$$h_j\circ p_j\circ h_{j-1}\circ p_{j-1}\circ\dots\circ h_1\circ p_1,$$
where $p_i$ is a polynomial and $h_i$ is a diffeomorphism in
$\mathcal{E}_{\beta}$ for $1\leq i\leq k$ and $1\leq k\leq n$,
the result follows from \cite[Proposition 11.5]{GMdM}.
\end{pf}

Stoilow Factorization together with compactness of 
the spaces of
(holomorphic)
polynomial-like mappings \cite[Theorem 5.8]{McM}
and $K$-qc mappings implies:

\begin{lem}
\label{lem:ahpl cpct}
For any $\delta>0,$ there exists $K_0\geq 1,$ so that for any $b\in\mathbb N,$
and $1\leq K\leq K_0,$ we have the following. The space of real $K$-quasiregular asymptotically holomorphic
polynomial-like mappings 
$f\colon U\to V$ of degree $b\geq 2$ with connected Julia sets, critical values in $(1+\delta)^{-1}I$, and
$\mod(V\setminus U)\geq \delta$ is compact up to affine conjugation. More
precisely, if $f_n\colon U_n\to V_n$ is a sequence of such mappings, then there exists
restrictions $f_n\colon U_n'\to V_n'$ (which automatically satisfy the conditions of the lemma
with $\delta$ replaced by $\delta'=\delta'(b,K)$), so that $\{f_n|_{U'_n}\}$ has a convergent subsequence.
\end{lem}

\begin{pf}
% Recall that the space of polynomial-like mappings is endowed with the
% Carath\'eodory topology: Assume that $f_n:U_n\to V_n$ is a sequence of
% polynomial-like mappings, and moreover, assume that for each $n,$ we choose a
% base point $u_n\in U_n.$
% The Carath\'eodory topology on the domains is defined as follows:
% We say that $(U_n,u_n)\to (U,u)$ if $u_n\to u$ and if for any subsequence $n_k\to\infty$ so that
% $(\hat{\mathbb{C}}\setminus U_{n_k})\to K$ in the Hausdorff topology on compact
% subsets of the sphere, we have that $U$ is equal to the component of
% $\hat{\mathbb{C}}\setminus K$ that contains $u.$ We say that $f_n:(U_n,u_n)\to
% (V_n,f_n(u_n))$ converges to $f:(U,u)\to (V,v)$ if $(U_n,u_n)$ converges to
% $(U,u)$ in the Carath\'eodory topology, and if for all $n$ sufficiently large,
% $f_n\to f$ uniformly on compact subsets of $U.$

% Suppose now that $f_n:U_n\to V_n $ is as in the statement of the Lemma.
Let
$\gamma_n$ denote the core curve of the annulus $V_n\setminus\overline{U}_n$.
Let $V_n'$ denote the region bounded by $\gamma_n$ in $\mathbb C.$
Since $\mod(V_n\setminus U_n)$ is bounded from below, $V_n'$ is a
$\kappa=\kappa(\delta)$-quasidisk. Let $U_n'$ denote $f_n^{-1}(V_n').$ Then
$U_n'$ is a $\kappa'=\kappa'(\kappa,\delta,b)$-quasidisk.  

By Stoilow Factorization, we can express
$f_n=g_n\circ\phi_n$ where $g_n\colon U_n'\to V_n'$ is analytic and $\phi_n\colon U_n'\to U_n'$
is $K$-quasiconformal. Moreover, since the proof of Stoilow factorization can be
carried out real-symmetrically, we can assume that $\phi_n$ is a real map, and
that the critical values of $g_n$ are real, and provided that $K$ is
sufficiently small, they are contained in the invariant real interval for $f_n$. 
By  \cite[Theorem 5.8]{McM}, we have that the family
$g_n\colon U_n\to V_n$ is compact in the Carath\'eodory topology. The $\phi_n$
belong to a compact family since they extend to $K'=K'(K,\kappa')$-qc mappings
of the plane, which we can assume are normalized to fix $M,-M,$ and $\infty$, as long as $M>0$
is chosen sufficiently large (such an $M$ exists by Theorem~\ref{thm:complex bounds}). Thus the family of
mappings 
$f_n\colon U_n'\to V_n'$ is compact too.
\end{pf}

\begin{prop}[{\it c.f.} \cite{GMdM}, Theorem 11.1]
\label{prop:smooth attractor}
There exists a compact set $\mathcal K$ of 
polynomial-like germs of type 
$\underline b$ with the following property:
Let $k\geq 3.$
For any $\varepsilon$ and
$f\in \Gamma_{\mathcal A^k_{\underline b}(I)},$
which is infinitely renormalizable at a critical point $c_0,$
there exists a sequence
$\{f_n\}\subset\mathcal K,$ and
$n_0\in\mathbb N$
such that for all $n\geq n_0$,
$$\|\mathcal R_{c_0}^n(f)-f_n\|_{\mathcal C^k(I)}\leq \varepsilon,$$
and $f_n$ is infinitely renormalizable with
the same combinatorics as $\mathcal R_{c_0}^n(f)$.
\end{prop}

\noindent\textit{Remark.}
In Proposition~\ref{prop:smooth attractor},
we have convergence of renormalization to a limit set in the
$\mathcal C^k$ topology; whereas,
\cite[Theorem 11.1]{GMdM} implies exponential convergence in
the $\mathcal C^{k-1}$ topology. 

\begin{pf}
We start with the following claim.

\medskip
\noindent\textit{Claim.}
There exists a compact set of polynomial-like germs
$\mathcal K$ such that given $f\in\Gamma_{\mathcal A^k(I)},$
which is infinitely renormalizable at a critical point $c_0,$
there exists a sequence $g_n\in\mathcal K$ so that
$\|g_n-\mathcal R_{c_0}^n(f)\|_{\mathcal C^0(I)}\rightarrow 0$ as $n\rightarrow\infty$,
and $g_n$ has the same combinatorics as $\mathcal R_{c_0}^n(f)$.

\medskip

\medskip
\noindent\textit{Proof of Claim.}
By Theorem~\ref{thm:complex bounds},
for each $n$ sufficiently big, there exist 
a $b$-tuple, $\underline b,$ and
an
asymptotically holomorphic polynomial-like renormalization,
$F_n\colon U_n\to V_n$ of type $\underline b$ of $f$
 with dilatation bounded
by $\diam(U_n)$. 
By Lemma~\ref{lem:Stoilow}, for each $n$ we can 
express $F_n$ as the composition 
$h_n\circ\phi_n\colon U_n\to V_n,$ where $\phi_n\colon U_n\to U_n$
is quasiconformal with dilatation bounded
by $\diam(U_n)$ and $h_n\colon U_n\to V_n$ is a
real polynomial-like mapping.
Let $U_n'=F_n^{-1}(U_n).$
By Lemma~\ref{lem:ahpl approx}
and Theorem~\ref{thm:complex bounds},
there exist a constant $C_0>0$ and  a real conformal mapping
$\psi_n$ so that
$$\|\phi_n-\psi_n \|_{\mathcal{C^0}(U'_n)}\leq C_0 C(U'_n)\diam(U_n),$$
% Since $\phi_n$ is a $\mathrm{diam}(U_n)$-qc  mapping on
% a $K$-quasidisk, we have that its restriction to $U_n'$ is almost affine,
% hence it is very close to the identity on the identity on $U'_n$.
% QQ
%the mappings $\phi_n$ converge to the identity in
%$\mathcal C^0(U')$,
%where $U'=F_n^{-1}(U)$.
and from \cite[page 53]{GdM},
we see that 
$$ C(U'_n)=\frac{4}{\pi} \sup_{z\in U'_n}\int\int_{U'_n}\Big|\frac{z(z-1)}{w(w-1)(w-z)}\Big|dxdy,
$$
which is uniformly bounded over $n.$
Thus, since $\mathcal R^n_{c_0}(f)$ has bounded geometry, by Theorem~\ref{thm:real bounds}
the mappings $h_n\circ\psi_n\colon U_n'\to U_n$
are polynomial-like mappings with connected Julia sets
and $\mod(U_n\setminus U_n')$ bounded from below.
 Thus
any limit of the $h_n\circ \psi_n$ is contained in a
 compact set of infinitely 
renormalizable polynomial-like germs
$\mathcal K$. For each $j,$ let $\mathcal V_j\subset\mathcal{PG}_{\underline b}$
be the neighbourhood of
$\mathcal K$ consisting of mappings so that 
their $j$-th polynomial-like renormalization persists over
$\mathcal V_j$. Since for any $j\in\mathbb N$, 
$h_n\circ \psi_n$ eventually enters $\mathcal V_j$,
we have that $h_n\circ \psi_n\colon U_n'\to U_n$ is
$j_n$ times renormalizable, where $j_n\rightarrow\infty$
as $n\rightarrow\infty$.

% since there 
% cannot be a sequence of mappings
% that are at most a bounded number of times renormalizable 
% which converges to an infinitely renormalizable mapping.
% Arguing by contradiction, if there exists $M>0$ 
% so that $m_n<M$ for infinitely many $n$, then there
% exists a sequence $h_{n_i}:U_{n_i}\to V_{n_i}$ of
% polynomial-like mappings converging to 
% $R^{n_i}F$
% (in the $\mathcal C^k$ topology on the real line),
% but each $h_{n_i}$ is at most $M$-times renormalizable.
% By compactness of the space of $K$-qr polynomial-like mappings 
% with moduli bounded from below, we obtain a subsequence of 
% the $h_{n_i}:U_{n_i}\to V_{n_i}$ which converges to the
% limit set of the renormalization operator, so 
% some subsequence of the
%  $h_{n_i}$ converges to an infinitely renormalizable mapping,
% which is impossible by the rigidity of infinitely
% renormalizable polynomial-like mappings,
% \cite{KSS-rigidity}.

% Each $h_n,$ with $n$ large, is an analytic polynomial-like 
% mapping,
% which is $j_n$ times
% renormalizable, with $j_n$ very large.
% By Theorem~\ref{thm:complex bounds}
% and \cite[Theorem 5.8]{McM} for $n$ sufficiently large, 
% the $h_n\circ \psi_n$ are contained in a compact family of analytic mappings.
Let $\delta>0$ be the universal constant
so that for all $n$ sufficiently large,
$\mod(U_n\setminus U_n')\geq\delta$, and let
$\Gamma'$ be the intersection of 
$\Gamma_{\mathcal{PG}_{\underline b}}$ with the 
set of polynomial-like germs with 
moduli bounded from below by $\delta$.
Let $\mathcal V_n'$ be the set of all polynomial-like germs
with moduli bounded from below by $\delta$ 
which are at least $n$-times renormalizable.
Then $\Gamma'=\cap_{n=0}^\infty\mathcal V_n'$.
Moreover, by Theorem~\ref{thm:lamination}
and Proposition~\ref{prop:contstr}
for any $\varepsilon>0$,
$\mathcal V_n$ is eventually contained in a
$\varepsilon$-neighbourhood of $\Gamma'$.
Thus for any $\delta>0$, if $n$ is sufficiently large,
there exists a polynomial-like mapping $g_n$ in
the topological conjugacy class of $\mathcal R_{c_0}^n(f)$
within distance $\delta$ from $h_n$
 in $\mathcal C^0(I)$.
% $h_n$ are eventually close in the Carath\'eodory topology
%  to a polynomial-like germ in the topological conjugacy class of
% $\mathcal R^nf$.
% By Theorem~\ref{thm:lamination},
% with a small $\mathcal C^0$ perturbation of $h_n$,
% which tends to zero as $n\rightarrow\infty$, we can
% obtain a polynomial-like mapping $g_n,$ 
% $\mathcal C^0$ close to $R^nf$ and in 
% $\mathcal H_{\mathcal R^n f}.$
\checkmark

\medskip

Associated to each $\underline b$, there
exist a family of polynomials $P$ and $\beta>0$, so that 
by Lemma~\ref{lem:Epstein class},
we have that there exists mapping $H_n$ in the
Epstein class, $\mathcal{E}_{\beta,P}$,
which is arbitrarily close to $\mathcal R_{c_0}^n(f)$
in the $\mathcal C^{k}$-topology.
Thus we have that $\|g_n-H_n\|_{\mathcal C^0(I)}$ is small.
Hence, since $g_n$ and $H_n$ are both analytic,
by Lemmas~\ref{lem:rbound} and \ref{lem:epdom},
we have that $\|g_n-H_n\|_{\mathcal C^{k}(I)}$ is small.
Thus
we have that $\|R^n(f)-g_n\|_{\mathcal C^{k}(I)}$ is small.
\end{pf}

\begin{lem}[\cite{ALdM}, Remark 2.7]
  Suppose that $f\in\mathcal A_{\underline b}^3(I)$ is an infinitely
  renormalizable mapping with a renormalization
  $\mathcal R(f)\in\mathcal A_{\underline b_1}^3(I)$ with the property that all
  the critical points of $\mathcal R(f)$ have the same $\omega$-limit set.
  Then $\mathcal R$ maps any
  sufficiently small open neighbourhood of $f$ to an open neighbourhood of
  $\mathcal R(f).$
\end{lem}
\begin{pf}
  Let us give the proof for a period-two renormalization. The proof in the
  general case is easier, since the closures of the orbit of the restrictive
  intervals are disjoint. Suppose that $f$ is renormalizable with period-two at
  a critical point $c_0.$ Let $J_0\owns c_0$ be the periodic interval of period
  two containing $c_0,$ and let $J_1=\mathrm{Comp}_{f(c_0)}f^{-1}(J_0).$ Let
  $\mathcal W\subset\mathcal A^3_{\underline b}(I)$ be an
  open neighbourhood of $f.$
  Provided that $\mathcal W$ is sufficiently small,
  it is contained in the domain of $\mathcal R.$ It is
  easy to see that the rescaling is open, we need only check that the
  composition is open. To simplify the notation, let us consider the
  renormalization of $f$  just as the composition
  $F=f|_{J_1}\circ f|_{J_0}\colon J_0\to J_0.$
  Choose $g_0\in\mathcal R(\mathcal W)\subset \mathcal A_{\underline b_1}^3(J_0),$
  and let $f_0\in \mathcal W$ be so that $\mathcal R(f_0)=g_0.$
  Fix $\varepsilon_0>0$ to be chosen later,
 and assume that $\|g-g_0\|_{\mathcal C^3}<\varepsilon$ with
$\varepsilon\in(0,\varepsilon_0).$ Let $\alpha$ be the common boundary point of
$J_0$ and $J_1$. Let $Q$ denote a small neighbourhood of $\alpha,$ so that
$f(Q\cap J_0)=Q\cap J_1.$
A straightforward calculation using the Taylor series of $f$ on $Q$
shows that there exists a constant $C_0>0$
and $f_1\in\mathcal C^3(Q)$, so that $\|f-f_1\|_{\mathcal C^3}<C_0\varepsilon,$
and $f_1|_{J_1\cap Q}\circ f_1|_{J_0\cap Q}=g|_{J_0\cap Q}.$
Let $\delta>0$ be so small that the interval $(1+\delta)^{-1}Q=\{x\in I:
|x-a|<(1+\delta/2)^{-1}|Q|/2\},$ where $a$ is the midpoint of $Q$, is a neighbourhood
of $\alpha.$  
Now, approximate $f$ by
$f_2\in\mathcal C^3(I),$ so that $\|f-f_2\|_{\mathcal C^3}<C_0\varepsilon,$ so
that $f_2=f_1$ on $(1+\delta)^{-1}Q.$ Finally, let $x_0$ be a point in
$(1+\delta)^{-1}Q\cap J_1,$ and let $x_1\in \partial J_1\setminus \{\alpha\}.$
Let $X=(f_2|_{J_0})^{-1}(x_0,x_1).$
There exists a constant $C_1>0,$ so that we can approximate $f_2$ on $(x_0,x_1)$ by $f_3\in\mathcal C^3((x_0,x_1))$ so that
$f_3\circ f_2|_X=g|_X.$ Notice that since the critical values of $f$ are not close to $x_1,$ we do not need to change $f$ in a small neighbourhood of $x_1.$
By construction $f_3$ extends to a $\mathcal C^3(I)$ mapping
$$
f_4(x):=\left\{
  \begin{array}{l}
    f_2(x),\quad x\in I\setminus (x_0,x_1)\\
    f_3(x), \quad x\in(x_0,x_1),
    \end{array}
      \right.
      $$
$f_4\circ f_4|J_0=g.$
and, provided that $\varepsilon_0$ was chosen sufficiently small, $f_4\in \mathcal W.$
  \end{pf}

\medskip
\noindent\textit{Proof of Theorem~\ref{thm:main smooth}.}
Suppose that $f\in \mathcal A^k_{\underline b}(I),$ 
$k\geq 3$,
is a mapping with zero topological entropy, and
which is infinitely renormalizable at a critical point $c$.
%where $\underline b$ is $b$-tuple with all even entries.
As usual, throughout the proof, $\mathcal R$ denotes the 
renormalization operator with period-doubling
combinatorics determined by the combinatorics
of restrictive intervals for $f$ about $c$.
By Theorem~\ref{thm:complex bounds}, for $N\in\mathbb N$ sufficiently large,
there exists $\mathcal W\subset\mathcal A^k_{\underline b}(I),$ 
a small open neighbourhood of $f$ chosen
so that each $g\in\mathcal W$ has an 
asymptotically holomorphic polynomial-like renormalization 
$\mathcal R^N g=G\colon U_G\to V_G.$
Let $\mathcal W'=\mathcal R^N (\mathcal W)$
%Each $G\in\mathcal W',$
%is asymptotically holomorphic.
and $F=\mathcal R^Nf$ .

By Proposition~\ref{prop:windows shrink},
to show that there are mappings with
positive entropy and mappings with finitely many 
periods in $\mathcal W'$,
it is enough
to show that there exists an analytic polynomial-like mapping
arbitrarily close to $F\colon U\to V$ in the $\mathcal C^k$-topology
on the real line.
As in the proof of Proposition~\ref{prop:smooth attractor},
it is sufficient to 
prove that we can approximate $F$ by
a polynomial-like mapping in the $\mathcal C^0$ 
topology on a complex neighbourhood of the interval.

% For $j\in\mathbb N$, let $\mathcal T_j$
% denote the connected component 
% containing $F$ of the 
% set of mappings that are $j$ times
% renormalizable,
% and
% for any $m,n$,
% let $\mathcal T_{m,n}$ denote the 
% neighbourhood of $\mathcal R^nF$ 
% consisting of mappings that
% are $m$ times renormalizable at $c$.

Let $n\in\mathbb N$.
There exists a $b_1$-tuple $\underline b_1$
with all entries even so that 
$\mathcal R^n(F)=F_n\colon U_n\to V_n$,
is an asymptotically holomorphic polynomial-like mapping 
of type $\underline b_1$.
Associated to $\underline b_1$, there is a family of
polynomials $P$, and $\beta>0$ so that by 
Lemma~\ref{lem:Epstein class}, 
for any $\varepsilon_1>0$, there exists a mapping 
$G_{n}\colon I\to I$ in the Epstein class, $\mathcal{E}_{\beta,P}$,
so that
$\|G_n-F_n\|_{\mathcal C^k(I)}<\varepsilon_1.$
Moreover, by the claim in the proof
of Proposition~\ref{prop:smooth attractor},
as $n\rightarrow\infty$, 
$F_n\rightarrow\mathcal K$ in $\mathcal C^0(U_n')$,
where $U_n'=F_n^{-1}(U_n)$. The mappings in 
$\mathcal K$ are analytic, so for $n$ sufficiently large,
if $G_n$ is sufficiently close to $F_n$ in $\mathcal C^k(I)$,
then $G_n$ is close to $\mathcal K$ in $\mathcal C^0(X),$
where $X$ is the open neighbourhood of the interval
given by Lemma~\ref{lem:epdom}.
Then, since the mappings in $\mathcal K$ are
polynomial-like mappings, for some $M\in\mathbb N\cup\{0\}$,
uniformly bounded in $\beta$,
$\mathcal R^MG_n\colon U_{\mathcal R^M G_n}\to V_{\mathcal R^M G_n}$ 
is a polynomial-like mapping.
Moreover, since we can take $\epsilon_1$ as small 
as we like, by continuity of the renormalization
operator in the $\mathcal C^k$ topology, see the appendix of 
\cite{AMdM}, we can assume that 
$\mathcal R^M G_n\in\mathcal R^{n+M}(\mathcal W)$, and
now we can conclude the proof as in Theorem~\ref{thm:AE}.

\medskip

\qed

\subsubsection{Proof of Theorem~\ref{thm:II smooth}}
The key step in the proof is
the construction of a codimension-$b$
manifold consisting of mappings
that are infinitely renormalizable at one
critical point, and whose remaining critical points
are periodic.

As usual, we say that a critical point $c$ of $f$ is
{\em non-degenerate} if $D^2f(c)\neq 0$. 
Let $b\in\mathbb N$, and
$\mathcal A^r_{even,b}(I)=\cup_{\underline b}\mathcal
A^r_{\underline b}(I),$ where the union is taken over 
all $b$-tuples $\underline b$ with all entries even. 
\begin{lem}
\label{lem:nondeg}
Let $r= 3+\alpha,$ where $\alpha>0$.
The set of mappings with all critical points non-degenerate is 
open and dense in $\Gamma_{\mathcal{A}_{even, b}^r}(I).$ 
\end{lem}
% \begin{pf}
% Suppose that there is an open set $Y$ of 
% mappings in $\Gamma_{\mathcal{A}_{even,b}(I)},$
% consisting of mappings with a degenerate critical point.
% Since there are at most countably many distinct $b$-tuples,
% $\underline b$, we can assume that 
% for some fixed
% $\underline b$
% each mapping in $Y$ is in
% $\mathcal A^r_{\underline b}(I)$.
% Let $f\in Y$
% From the proof of 
% Theorem~\ref{thm:main smooth},
% we have that there exists a neighbourhood
% $\mathcal U$ of $f$ so that 
% for some fixed $b_1$-tuple,
% $\underline b_1$,
% every mapping
% in $\mathcal U$ has as AHPL renormalization
% of type $\underline b_1$.
% Moreover, choosing $\mathcal U$
% smaller if necessary, by Lemma~\ref{lem:Epstein class},
% we have each mapping
% $g\in\mathcal U$ can be approximated
% by a mapping in the Epstein class, $\mathcal E_{\beta,P}$,
% where $P$ is a sequence of polynomials depending on 
% $\underline b_1$. 
% \end{pf}

\begin{pf}
Let $\underline 2$ denote the $b$-tuple
where every entry is a two.
It is well-known that
the set  $\mathcal{A}_{\underline 2}^r(I)$
of mappings with
all critical points non-degenerate
is open and dense in the space $\mathcal A^r_{even,b}(I)$,
\cite{Wall}.
Thus the set of mappings with all critical points non-degenerate is
relatively open in $\Gamma_{\mathcal A^r_{even,b}(I)}$.

We will now prove density.
Let us assume that $f$ has exactly one solenoidal 
attractor, the case when it has more than one is similar.
Let $f\in\Gamma_{\mathcal A^r_{even, b}(I)},$
and let $\mathcal U$ be an open neighbourhood of $f$.
Let $c$ be a critical point at which $f$ is infinitely
renormalizable. Let $J_n\owns c$ be 
a restrictive interval of period  $2^n$.
Then
for $n$ sufficiently large,
we have that each interval 
$J_n^i=\mathrm{Comp}_{f^i(c)}f^{-(2^n-i)}(J_n)$,
$ i\in\{0,1,\dots,2^n-1\}$, contains at most one
critical point. Moreover, by Lemma~\ref{K_n},
there exists a neighbourhood 
$\mathcal U\subset \mathcal{A}^r_{even,b}(I)$ of $f$
so that
$\cup_{i=0}^{2^n-1}J_n^i$ persists over
$\mathcal U$. 
By Theorem~\ref{thm:main smooth}, and the fact that
the set of mappings with a non-degenerate critical point 
is open and dense, we can approximate
$f$ by mappings $f_0,f_1\in\mathcal U\cap\mathcal{A}_{\underline 2, b}(I),$
where $f_0$ is in the interior of mappings with zero entropy and
$f_1$ has positive entropy.
Let $P(x)=x^2$.
Since $f_0,f_1\in\mathcal U$, we can express 
$\mathcal R^{2^n}(f_0)=h_{0,b}\circ P\circ h_{0,b-1}\circ P\circ\dots\circ
h_{0,1}\circ P,$ and
$\mathcal R^{2^n}(f_1)=h_{1,b}\circ P\circ h_{1,b-1}\circ P\circ\dots\circ
h_{1,1}\circ P,$ where each $h_{i,j}\colon I\to I,$ $i\in\{0,1\}$
and $j\in\{1,\dots, b\},$ is a $\mathcal C^r$ diffeomorphism of the
interval.
Now, for each $j\in\{1,2,\dots,b\}$, let 
$h_{j}^{\lambda}\colon I\to I,$ $0\leq \lambda\leq 1$ be a path of 
diffeomorphisms between
$h_{0,j}$ and $h_{1,j}$. Thus we obtain a path
$F_{\lambda}=h_{b}^{\lambda}\circ P\circ\dots\circ h_1^{\lambda}\circ P$
of multimodal mappings from 
$\mathcal R^{2^n}(f_0)$ to $\mathcal R^{2^n}(f_1).$
Moreover, we can assume that the diameter of the path is
as small as we like by choosing
$f_0, f_1$ close enough to $f$. Taking the preimage
of the path under $\mathcal R^{2^n},$ we obtain
a path $f_{\lambda}$ from $f_0$ to $f_1$,
which crosses $\Gamma_{\mathcal A^r_{\underline 2, b}}.$
Thus there exists a mapping with 
all critical 
points non-degenerate arbitrarily close to $f$.

When $f$ has more than one solenoidal attractor, we have to choose the 
mapping in the interior of zero entropy as we did at the end
of Theorem~\ref{thm:AE}.

\end{pf}

We will make use of the
period-doubling renormalization
operator acting on unimodal mappings 
with non-degenerate critical points,
\cite{Davie}. Let
$\alpha>0$. We let $\mathcal{A}^{2+\alpha}_{2}(I)$
be the space of unimodal $\mathcal C^{2+\alpha}(I)$ mappings
on the interval with a non-degenerate critical point.
The period-doubling renormalization operator
acting on the $\mathcal C^{2+\alpha}(I)$ has a unique fixed point,
$f_*$. By Sullivan's complex bounds, 
\cite{Sullivan}, we
can regard $f_*$ as a quadratic-like germ.
Moreover, at $f_*$ the renormalization
operator is hyperbolic. 
Let $\mathbf u_*$ denote the unstable vector
at $f_*$. The next proposition
describes the stable manifold.

\begin{prop}\cite{Davie}
\label{prop:Davie}
Let $\alpha>0$.
The local stable set of 
$f_*,$ $W^{s,2+\alpha}_{\varepsilon}\subset
\mathcal{A}^{2+\alpha}_{2}(I)$ is
a codimension-one, $\mathcal C^1$-submanifold.
\end{prop}

Let us say that a 
multimodal mapping of type $\underline b,$
$f,$
with critical points $\{c_1,c_2,\dots,c_{b}\}$
has combinatorics
$\sigma^0_*$ if $b-1$ of its critical points are contained
in a periodic cycle and at the remaining critical point,
say $c_0$, $f$ is infinitely renormalizable with period-doubling
combinatorics.

\begin{lem}\cite[Proposition 8.7]{dFdMP}
\label{lem:515}
For real numbers $r>s+1\geq 2$, the composition
operator from $\mathcal C^r\times\mathcal C^{s}\to\mathcal{C}^s$
is a $\mathcal C^1$ mapping.
\end{lem}

\begin{prop}[{\it c.f.} \cite{dFdMP}, Theorem 9.1]
\label{prop:C1mfld}
For every 
$r> 3$,
if $f$ has combinatorics 
$\sigma^0_*$, then the 
connected component containing $f$
of the topological conjugacy class of $f$
is an embedded,
codimension-$b,$ $\mathcal C^1,$ Banach submanifold 
of the space of smooth multimodal mappings.
\end{prop}

\begin{pf}
Let $\alpha>0$ be so that $r=3+\alpha$,
and choose $0<\alpha'<\alpha$.
By Proposition~\ref{prop:Davie}, the local
stable manifold through $f_*$ in the space
$\mathcal C^{2+\alpha'}$ is a  codimension-one 
$\mathcal C^1$-submanifold. Let us denote it by
$W^{s,2+\alpha'}_{\varepsilon}.$
We may assume that $\varepsilon>0$ so small that the
vector $\mathbf{u}_{*}\in T_{f_0}\mathcal A^r_{\underline b}(I)$ is transversal to 
the local stable set $W_{\varepsilon}^{s,2+\alpha'}(f_*)$
at each $f\in W_{\varepsilon}^{s,2+\alpha'}(f_*).$

Let $g\in W^{s,3+\alpha}(f_*),$ the stable set of $f_*$
in $\mathcal{A}^{3+\alpha}_{\underline b}(I).$ 
There exists $N=N(g)>0$ so large that 
$$\mathcal R^N(g)\in W^{s,3+\alpha}_{\varepsilon}(f_*)\subset
W^{s,2+\alpha'}_{\varepsilon}(f_*).$$

Since $v=\mathbf u_{*}$ is transversal at $\mathcal R^N(g)$
to $W^{s,2+\alpha'}_{\varepsilon}(f_*),$
there exist a small open set 
$\mathbb O_0\subset \mathcal{A}^{2+\alpha'}_{\underline b}(I)$
containing $\mathcal R^N(f_0)$ and a $\mathcal C^1$ function
$\Phi\colon \mathbb O_0\to\mathbb R$
such that 
$\Phi^{-1}(0)=W_{\varepsilon}^{s,2+\alpha'}(f_*)
\subset \mathbb O_0$ for which $0\in\mathbb R$ is a
regular value and 
$D\Phi(\mathcal R^N(g))v\neq 0.$

By Lemma~\ref{lem:515},
the operator $\mathcal R^N$ is a $\mathcal C^1$ map from
$\mathcal{A}^{3+\alpha}_{\underline b}(I)$
into $\mathcal{A}^{2+\alpha'}_{\underline b}(I)$.
Let $\mathbb O_1\subset \mathcal{A}^{3+\alpha}_{\underline b}(I)$ 
be an open set 
containing $g$ such that 
$\mathcal R^N(\mathbb O_1)\subset\mathbb O_0$.
We want to show that $0\in\mathbb R$ is a
regular value for
$\Phi\circ \mathcal R^N\colon \mathbb O_1\to\mathbb R.$
Defining $g_t=\mathcal R^N(g)+tv$, with $|t|$ small,
we get a $\mathcal C^1$ family $\{g_t\}$ of mappings in
$\mathcal{A}^{3+\alpha}_{\underline b}(I)$, which is transversal to 
$W^{s,2+\alpha'}_{\varepsilon}(f_*)$ at $g_0=\mathcal R^N(g).$

\medskip
\noindent\textit{Claim.} There exists a 
$\mathcal C^1$ family $\{G_t\}\subset\mathcal{A}^{3+\alpha}_{\underline b}(I)$
such that for all small $t$ we have
$\mathcal R^N(G_t)=g_t$. Moreover,
for each of the $b-1$ critical points 
$c_t$ of $G_t$ which do not correspond to the 
critical point $c_0$ of $g$, the 
itinerary of $c_t$ is the same as the 
itinerary of $c$,
where
$c_t\in\mathrm{Crit}(G_t)$ naturally corresponds to
$c\in\mathrm{Crit}(G)$.
\medskip

\noindent\textit{Proof of Claim.}
First note that $g_t=h_t\circ g_0$ 
where each $h_t\in \mathcal C^k(I)$
is a diffeomorphism.
Since $\mathcal R^N(g)=g_0$,
there exist $p>0$ and closed, pairwise disjoint intervals
$0\in\Delta_0,\Delta_1,\dots,\Delta_{p-1}\subset I$
with $G(\Delta_i)\subset\Delta_{i+1}$ for 
$0\leq i\leq p-1,$ and $G(\Delta_{p-1})\subset\Delta_0$,
such that
$$g_0=\mathcal R^N(G)=\Lambda_G^{-1}\circ G^p\circ \Lambda_G,$$
where $\Lambda_G\colon I\to\Delta_0$ is an affine mapping.

Let $\bar h_t\colon \Delta_0\to\Delta_0$ be the $\mathcal C^k$
diffeomorphism given by 
$\bar h_t=\Lambda_G\circ h_t\circ \Lambda_G^{-1}$.
Let $\Delta_{-1}$ denote the union of immediate basins of
attraction of super-attracting cycles of $G.$
Consider a $\mathcal C^k$ extension of
$\bar h_t$ to a diffeomorphism $H_t\colon I\to I$ with the
property that $H_t|_{\Delta_i}$ is the identity for all 
$i\neq 0$. Then let $G_t\in \mathcal{A}^{3+\alpha}_{\underline b}(I)$
be the map $G_t=H_t\circ f$.
Note that $G_t^i(0)=G^i(0)$ for all
$0\leq i\leq p$, that $G_t$ is $N$-times
renormalizable under $\mathcal R$
and that $\mathcal R^N(G_t)=h_t\circ g_0=g_t$.
\checkmark

\medskip

Now let us show that the claim proves
the proposition.
Observe that the condition that 
$b-1$ of the critical points of $g$ lie in
an periodic cycle defines a codimension
$b-1$ subspace of $\mathcal{A}^{3+\alpha}_{\underline b}(I)$.
Setting
$$w=\frac{d}{dt}\Big|_{t=0}G_t,$$
we obtain that
$$D(\Phi\circ \mathcal R^N)(G)w=D\Phi(g_0)v\neq 0.$$
Therefore $\Phi\circ \mathcal R^N$ is a $\mathcal C^1$ local submersion 
at $G$. By the Implicit Function Theorem,
$(\Phi\circ \mathcal R^N)^{-1}(0)$ is a codimension-one,
$\mathcal C^1$ Banach submanifold of 
$\mathbb O_1$ an open subset of
$\mathcal{A}^{3+\alpha}_{\underline b}(I)$. 
Furthermore, if 
$h\in(\Phi\circ \mathcal R^N)^{-1}(0)$, then $\mathcal R^N(h)\in
W^{s,2+\alpha'}_{\varepsilon}(F_*),$
and so $h$ belongs to the 
global stable set $W^{s,2+\alpha'}(F_*)$.
By Proposition~\ref{prop:smooth attractor},
we have that $h$ in fact belongs to $W^{s,r}(g)$. 
The proposition follows.
\end{pf}

% Repeating the same proof, using only the continuity of 
% $\mathcal R$, we have:
% \begin{cor}
% If $F\in\mathbb V^k,$ $k\in\{3,4,\dots,\infty\}$ 
% has combinatorics 
% $\sigma^0_*$, then the 
% connected component containing $F$
% of the set of mappings
% with the same combinatorics as $F$
% is a 
% codimension-$b$ submanifold 
% of $\mathbb V^3$.
% \end{cor}
% \begin{pf}
% Let us point out how to modify the proof of
% Proposition~\ref{prop:C1mfld}. It no
% longer makes sense to prove that
% $0$ is a regular value of the composition
% $\Phi\circ \mathcal R^N:\mathbb O_1\to\mathbb R$.
% However, we still have that
% there exists a $\mathcal C^1$ family
% $\{G_t\}\subset\mathbb V^3$ so that
% for all $t$ small we have
% $\mathcal R^N(G_t)=g_t,$
% as in the claim in the proof of
% Proposition~\ref{prop:C1mfld}.

% \end{pf}

\medskip
\noindent\textit{Proof of
Theorem~\ref{thm:II smooth}.}
By Lemma~\ref{lem:nondeg},
we have that the family of 
mappings with all critical points 
non-degenerate is open and
dense in $\Gamma=\Gamma_{\mathcal A^r_{even, b}(I)}$.
Let $\Gamma_1$ be the subset of mappings in 
$\Gamma $ that has exactly one solenoidal attractor.
Arguing as in the proof of 
Theorem~\ref{thm:B1}, using 
Theorem~\ref{thm:main smooth} in place of 
Theorem~\ref{thm:main}, we have that 
$\Gamma_1$ is open and dense in $\Gamma$.
Let $\mathcal X\subset\Gamma$ 
denote the open, dense set of mappings with all critical points
non-degenerate and exactly one solenoidal attractor.
We need to show that any mapping 
$f\in\mathcal X$
can be approximated by mappings 
in $\Gamma$ with exactly one solenoidal 
attractor containing
exactly one critical point, which is non-degenerate.

Let $c$ be a recurrent critical point of 
$f$ such that $\omega(c)$ is a solenoidal 
attractor, and
let $F\colon U\to V$ be an asymptotically
holomorphic polynomial-like renormalization
of $f$ at $c$.
By Proposition~\ref{prop:smooth attractor},
there exists a compact set of polynomial-like germs 
$\{f_n\}$ such that for all $\varepsilon>0$ and 
all $n$ sufficiently large
$$\|\mathcal R^n(F) -f_n\|_{\mathcal C^r(I)}<\varepsilon,$$
which are infinitely renormalizable with the 
same combinatorics as $\mathcal R^n(F)$.
By Theorem~\ref{thm:II}
we have that we can approximate each 
$f_n$ by polynomial-like germs $g_m$
which are infinitely renormalizable at 
one critical point and with $b-1$ periodic
critical points.
By Proposition~\ref{prop:C1mfld}
for each $g_m$,
there is a codimension-$b$,
submanifold, $\mathcal{H}_{g_m},$ 
of $\mathcal{A}^{r}_{\underline 2, b}(I)$ 
consisting of mappings 
topologically conjugate to $g_m$,
and
any sequence of mappings,
$\hat g_m\in \mathcal{H}_{g_m},$
accumulates on the topological conjugacy class
of $f_n$ in $\mathcal{A}_{\underline 2,b}(I)$.
Arguing as in the proof of 
Theorem~\ref{thm:main smooth}
any mapping in the topological conjugacy class
of $f_n$ in $\mathcal{A}_{\underline 2,b}(I)$
can be approximated by such mappings, $\hat g_m$.
By Theorem~\ref{thm:qs rigidity},
for $\varepsilon>0$, sufficiently small,
these manifolds laminate $B_{\varepsilon}(f_n)
\subset\mathcal{A}^{r}_{even, b}(I).$ 
So for any neighbourhood 
$\mathcal U'\subset\mathcal A^r_{even, b}(I)$
of $F$, there exists $n$ so that
$\mathcal R^n(\mathcal U'\cap \mathcal T_n)$ intersects
such a topological conjugacy class,
where $\mathcal T_n$ is the set of mappings
which are $n$ times renormalizable. 
We can conclude by arguing as in 
the proof
of Theorem~\ref{thm:AE}.
\qed

\medskip

Finally we obtain:

\begin{thm}
Let $r>3$ and $b\in\mathbb N$.
Let $\underline b$ be a $b$-tuple consisting of even integers.
Each connected component of
$\Gamma_{\mathcal A_{\underline b}^r(I)}$ is locally connected.
\end{thm}

\begin{pf}
Let $\Gamma$ denote a connected component of
$\Gamma_{\mathcal A^r(I)},$ and suppose
that there exists 
$f\in\Gamma$, so that
$\Gamma$ is not locally connected at $f$.
Then there is an arbitrarily small
open set $\mathcal V\subset\mathcal{A}^r(I)$, 
with $f\in \mathcal V$, such that for
every open set $\mathcal U\subset \mathcal V,$ with
$\mathcal U\owns f$,
we have that $\mathcal U\cap\Gamma$ is not connected.
Take $\varepsilon>0$ small enough so that 
$B_{\varepsilon}(f)\subset \mathcal V,$ and set 
$\mathcal U=B_{\varepsilon}(f).$
Since $\Gamma$ is a closed set that is not
locally connected at $f$,
$\Gamma\cap \mathcal U$ has infinitely many components:
If $\Gamma\cap \mathcal U$ contained only
finitely many components,
$\Gamma_0\cup\dots\cup\Gamma_k$, with $f\in \Gamma_0$, then 
$\Gamma_1\cup\dots\cup \Gamma_k$ is a relatively closed subset of 
$\mathcal U$, but now, there is an open set 
$\mathcal U'\subset \mathcal U$ so that
$\mathcal U'\cap \Gamma=\Gamma_0$ 
is connected, which 
contradicts the choice of $\mathcal V$.
Thus we have that $\Gamma\cap \mathcal U$ consists of
infinitely many
 connected components, 
which must accumulate on $f$.

Since $\Gamma$ is connected, by Theorem~\ref{thm:II smooth},
there are codimension-one components $\Gamma_n$ of $\Gamma\cap \mathcal U$, 
so that $\dist(\Gamma_n, f)$ is arbitrarily small, and with 
$\diam(\Gamma_n)\geq \varepsilon/2,$ since they must connect points close to $f$
with points outside $B_{\varepsilon}(f).$
Even more, since
by Theorem~\ref{thm:main smooth},
 $\mathcal U\setminus\Gamma$
 consists of two open sets, one in the interior of mappings with zero entropy,
 and one consisting of mappings with positive entropy,
we have that for each $n$ sufficiently big  
$\partial(B_{\varepsilon/4}(f)\cap\Gamma_n)\subset\partial
B_{\varepsilon/4}(f).$

Let $\mathcal Z\subset \Gamma$, denote the set of 
mappings with exactly one solenoidal attractor, which
contains exactly one non-degenerate critical point and no 
others. By Theorem~\ref{thm:II smooth}, we have
that $\mathcal Z$ is a union of codimension-one
open sets, which is dense in $\Gamma$.

Suppose first that $f\in\mathcal Z.$ 
Then there is a neighbourhood $\mathcal U_1$ of $f$,
and a renormalization $\mathcal R^n$, so that 
$\mathcal R^n(\mathcal U_1)$ is contained in the 
space of asymptotically holomorphic
polynomial-like mappings with non-degenerate critical points.
Moreover, taking a deeper renormalization if necessary,
we can assume that $\mathcal R^{n'}(\mathcal U_1)$
is contained in an arbitrarily small neighbourhood 
of the quadratic-like fixed point of renormalization.

By the claim in the proof of Proposition~\ref{prop:C1mfld},
there exists a
$\varepsilon'>0$ and a neighbourhood $N$ of $f$ in $\mathcal Z$
so that for each $g\in N$, there is a transverse family
$\{g_t\}_{|t|<\varepsilon'}$ to $g=g_0$, 
so that $\mathcal R^{n'}(g_t),
|t|<\varepsilon',$
is a transverse family to the 
local stable manifold of renormalization,
and $\mathcal R^{n'}$ is injective on 
$\{g_t\}_{|t|<\varepsilon'}.$
But now, since each $\Gamma_n$ has codimension-one, 
and they accumulate on $f$, there exist
arbitrarily large $n$ so that for
$g_0$ close to $f$,
$\Gamma_n\cap\{g_t\}_{|t|<\varepsilon'}\neq\emptyset.$
So, since having exactly one solenoidal attractor is an open property in
$\Gamma,$
there exist $g_n\in\Gamma_n\cap\{g_t\}_{|t|<\varepsilon'},$
converging to $f$, so that $\mathcal R^{n'}(g_n)$
is an infinitely renormalizable quadratic-like mapping.
This contradicts the injectivity of
$\mathcal R^{n'}$ on each transverse family.
Thus $\Gamma$ is locally connected at
$f\in\mathcal Z$.

Now assume that $f$ is an arbitrary mapping in $\Gamma$.
In each $\Gamma_n$, there is a dense 
set of relatively
open manifolds consisting of mappings in
$\mathcal Z$.
Since each 
$\Gamma_{n}$ has the property that 
$\partial(B_{\varepsilon/4}(f)\cap\Gamma_n)\subset\partial 
B_{\varepsilon/4}(f)$, we have that
the set $\mathcal Y$
of all limit points of the $\Gamma_{n}$ contains a 
codimension-1 connected submanifold of $\Gamma$, 
contained in 
$\overline{\mathcal U\cap\Gamma}$.
Thus $\mathcal Z$ is dense in $\mathcal Y,$
and points in $\mathcal Z\cap\mathcal Y$ are accumulated by
points in $\Gamma_n$. But this contradicts
the fact that
$\Gamma$ is locally connected at
$f\in\mathcal Z$.
\end{pf}

\medskip\textit{Acknowledgements.}
The authors would like to thank 
Henk Bruin for several helpful comments.
The authors are grateful to Sebastian van Strien 
for suggesting this problem and for
numerous inspiring conversations while this work was being
carried out.

%%% Local Variables:
%%% mode: latex
%%% TeX-master: t
%%% End:

\end{document}